\documentclass{amsart}
\newtheorem{theorem}{Theorem}
\newtheorem{lemma}{Lemma}
\newtheorem{corollary}{Corollary}

\numberwithin{equation}{section}

\newcommand{\scr}{\scriptstyle}
\def\star{\raise .5ex \hbox{*}}
\def\sumstar_#1{\setbox0=\hbox{$\scriptstyle{#1}$}
\setbox2=\hbox{$\displaystyle{\sum}$}
\setbox4=\hbox{${}\star\mathsurround=0pt$}
\dimen0=.5\wd0 \advance\dimen0 by-.5\wd2
\ifdim\dimen0>0pt
\ifdim\dimen0>\wd4 \kern\wd4 \else\kern\dimen0\fi\fi
\mathop{{\sum}\star}_{\kern-\wd4 #1}}
\def\sumprime_#1{\setbox0=\hbox{$\scriptstyle{#1}$}
\setbox2=\hbox{$\displaystyle{\sum}$}
\setbox4=\hbox{${}'\mathsurround=0pt$}
\dimen0=.5\wd0 \advance\dimen0 by-.5\wd2
\ifdim\dimen0>0pt
\ifdim\dimen0>\wd4 \kern\wd4 \else\kern\dimen0\fi\fi
\mathop{{\sum}'}_{\kern-\wd4 #1}}
\font\ger=eufm10
\newcommand{\gs}{\hbox{\ger S}}
\newcommand{\frU}{\hbox{\ger U}}
\def\ndiv{\not \hskip .03in \mid}
\def\lr{\lambda_{R}}
\def\e{\epsilon}
\def\calm{\mathcal{M}}
\def\modulo{\mathrm{mod}\,}
\begin{document}

\title{Higher Correlations of Divisor Sums Related to Primes
II}
\author{D. A. Goldston}
\address{Department of Mathematics, San Jose State University, San Jose, CA
95192, USA} \email{goldston@mathcs.sjsu.edu}
\thanks{Goldston was supported by NSF grants and the American Institute of Mathematics}
\thanks{ Y{\i}ld{\i}r{\i}m was supported
by T\"{U}B\.{I}TAK}
\author{ C. Y. Yildirim}
\address{ Department of Mathematics, Bo\~{g}azi\c{c}i University,
Bebek, Istanbul, 34342 Turkey \& \newline \hspace*{.3cm}
Feza G\"{u}rsey Enstit\"{u}s\"{u}, \c{C}engelk\"{o}y,
Istanbul, P.K. 6, 81220 Turkey}
\email{yalciny@boun.edu.tr}
\date{\today}
\begin{abstract}
We calculate the triple correlations for the truncated divisor sum
$\lambda_{R}(n)$. The $\lambda_{R}(n)$'s behave over certain averages just as
the prime counting von Mangoldt function $\Lambda(n)$ does or is conjectured
to do.
We also calculate the mixed (with a factor of $\Lambda(n)$) correlations.
The results for the moments up to the third degree, and therefore
the implications for the distribution of primes in short intervals,
are the same as those we obtained (in the first paper with this title) by using 
the simpler approximation $\Lambda_{R}(n)$.  However, when $\lambda_{R}(n)$ is 
used the error in the singular
series approximation is often much smaller than what $\Lambda_{R}(n)$ allows. 
Assuming the Generalized Riemann Hypothesis for Dirichlet $L$-functions, we obtain an
$\Omega_{\pm}$-result for the variation of the error term in the prime 
number theorem. 
Formerly, our knowledge under GRH was restricted to $\Omega$-results for the 
absolute value of this variation. 
An important ingredient in the last part of
this work is a recent result due to Montgomery and Soundararajan which 
makes it possible for us to dispense with a large error term in the
evaluation of a certain singular series average.
We believe that our results on $\lambda_{R}(n)$'s and $\Lambda_{R}(n)$'s
can be employed in diverse problems concerning primes.
\end{abstract}

\maketitle

\section{  Introduction}

In this paper we calculate the triple correlations of the short divisor sum
defined by
\begin{equation}
\lambda_R(n) = \sum_{r\leq R} \frac{\mu^2(r)}{\phi(r)}
\sum_{d|(r, n)} d\mu(d), \qquad \mathrm{ for} \ n\ge 1, 
\end{equation}
and $\lambda_R(n)=0$ if $n\le 0$. In the previous paper of
this series \cite{GY2} we gave the calculation of the correlations of
\begin{equation}
\Lambda_R(n)= \sum_{\stackrel{\scr d |n}{\scr d\le R}}\mu(d) \log (R/d),
\qquad \mathrm{ for} \ n\ge 1, 
\end{equation}
and $\Lambda_R(n) = 0$ if $n\le 0$. As can be seen from our results,
these divisor sums tend to behave similarly
to the prime counting von Mangoldt function $\Lambda(n)$,
and thus they may sometimes be used in
place of $\Lambda(n)$ when it is not possible to work directly with
$\Lambda(n)$ itself. Since
\[ \Lambda(n) = \sum_{d|n}\mu(d)\log(R/d), \qquad \mathrm{ for} \ n > 1, \]
$\Lambda_{R}(n)$ comes about as a surrogate for $\Lambda(n)$ by truncation.
We can relate $\lambda_{R}(n)$ to $\Lambda_{R}(n)$ by
interchanging the order of the summations in (1.1), 
thereupon the new inner sum can
be evaluated (eq. (2.15) below) and the contribution of its main term
gives $\Lambda_{R}(n)$. 

Goldston \cite{G2} found $\lambda_{R}(n)$ while remedying the
failure of the circle method in an application to the related problems
of twin primes and short gaps between primes for which a starting point is
the observation that
\begin{equation}
\sum_{n\leq N}\Lambda(n)\Lambda(n+k) = \int_{0}^{1}|S(\alpha)|^{2}e(-k\alpha)\,
d\alpha \: + O(k\log^2 N), 
\end{equation}
where
\[ S(\alpha) = \sum_{n\leq N}\Lambda(n)e(n\alpha), \qquad e(u) = e^{2\pi iu}.\]
For $\alpha$ close to the rational number $a/r$, we write
$\alpha = a/r + \beta$, and approximate $S(\alpha)$ throughout $[0, 1]$
by a sum of local approximations
\[ \sum_{r\leq R}\sum_{\stackrel{\scr 1\leq a \leq r}{\scr (a,r)=1}}
{\mu(r)\over \phi(r)}I({a\over r} + \beta) ,
\qquad \mathrm{where} \qquad I(u) =
\sum_{n\leq N}e(nu) . \]
But the last expression is equal to
\[ \sum_{n\leq N}\lambda_{R}(n)e(n\beta), \] suggesting we replace
$\Lambda(m)$ by $\lambda_{R}(m)$
in sums such as (1.3). Furthermore, Goldston \cite{G2}
showed that among sums of the form
\[ \sum_{\stackrel{\scr r \leq R}{\scr r|n}}a(R,r) \qquad \mathrm{with} \qquad
a(R,1) = 1,\;\; a(R,r) \in \mathbb{R} , \]
$\lambda_{R}(n)$ is the best approximation to $\Lambda(n)$ in an $L^2$ sense.
The proof involves a minimization which was solved in a more general setting
by Selberg \cite{S1} for his upper bound sieve. Hooley's recent use of
$\lambda_{R}(n)$ in \cite{H12}, \cite{H13} leans much on its origin in the
Selberg sieve. It should further be mentioned that, as far as we know,
Heath-Brown \cite{HB} was the first to use $\lambda_{R}(n)$ in additive
prime number theory.

The correlations we are interested in evaluating are
\begin{equation}
\mathcal{ S}_k(N,\mbox{\boldmath$j$}, \mbox{\boldmath $a$})
=\sum_{n=1}^N \lambda_R(n+j_1)^{a_1}\lambda_R(n+j_2)^{a_2}\cdots
\lambda_R(n+j_r)^{a_r} 
\end{equation}
and
\begin{equation}
\tilde{\mathcal{S}}_k(N,\mbox{\boldmath$j$}, \mbox{\boldmath$a$}) =
\sum_{n=1}^N \lambda_R (n+j_1)^{a_1}\lambda_R (n+j_2)^{a_2}\cdots
\lambda_R(n+j_{r-1})^{a_{r-1}}\Lambda(n+j_r) 
\end{equation}
where $\mbox{\boldmath$j$} = (j_1,j_2, \ldots , j_r)$ and $\mbox{\boldmath$a$}
= (a_1,a_2, \ldots a_r)$, the $j_i$'s are distinct integers, $a_i\geq 1$ and
$\sum_{i=1}^r a_i = k$. In (1.5) we assume that $r\ge 2$ and take
$a_r=1$. For later convenience we define
\begin{equation}
\tilde{\mathcal{ S}}_1(N, \mbox{\boldmath$j$},
\mbox{\boldmath$a$}) =\sum_{n=1}^N\Lambda(n+j_1) = \psi(N) +
O(|j_{1}|\log N) \sim N 
\end{equation}
for $\displaystyle|j_1| = o({N\over \log N})$ by the prime number theorem
(as usual $\displaystyle \psi(x) = \sum_{n\leq x}\Lambda(n)$).

For $k=1$ and $k=2$ these correlations have been evaluated before
(\cite{G1}, \cite{H12}, \cite{H13}), and the more general cases of
$n$ running through arithmetic progressions were also worked out
(\cite{GY1}, \cite{HB}, \cite{H13}). 

Correlations which include in their summands factors such as
$\!\Lambda(n)\Lambda(n+j), j\neq 0,$
cannot be evaluated unconditionally; they are the subject of the
Hardy-Littlewood prime $r$-tuple conjecture \cite{HL}.
This conjecture states that for $\mbox{\boldmath$j$} = (j_1,j_2, \ldots, j_r)$
with the $j_i$'s distinct integers,
\begin{equation}\psi_{\mbox{\boldmath$j$}}(N) = \sum_{n=1}^N \Lambda(n+j_1)
\Lambda(n+j_2)\cdots \Lambda(n+j_r) \sim \gs(\mbox{\boldmath$j$}) N
\end{equation}
when $\gs(\mbox{\boldmath$j$})\neq 0$, where
\begin{equation}
\gs(\mbox{\boldmath$j$}) = \prod_p\left(1- \frac{1}{p}\right)^{-r}
\left(1-\frac{\nu_p(\mbox{\boldmath$j$})} {p}\right)
\end{equation}
and $\nu_p(\mbox{\boldmath$j$})$ is the number of distinct residue classes
modulo $p$ that the $j_i$'s occupy. If $r=1$ we see $\gs(\mbox{\boldmath$j$})
=1$, and for $|j_1|\le N$ (1.7) reduces to (1.6), which is the only case where
(1.7) has been proved. The cases $r=2, 3$
will be of particular interest to us in this paper, the explicit
expressions have been shown in \cite{GY2} to be
\begin{eqnarray}
\mbox{} \qquad \qquad \gs((0,j)) & = & \gs_{2}(j), \qquad (j\neq 0) \\
\mbox{} \qquad \qquad \gs((0,j_{1},j_{2})) & =
& \gs_{2}((j_{1},j_{2}))\gs_{3}(j_{1}j_{2}(j_{1}-j_{2})),
\;\; (j_{1}\neq j_{2}, \, j_{1}j_{2} \neq 0),
\end{eqnarray}
where writing
\begin{equation}
        p(n) =
                \left\{
                \begin{array}{ll}
        n, &\mbox{if $n$ is a prime,}\\
                  1, & \mbox{otherwise},
                \end{array}
                \right.
\end{equation}
the singular series for $n\ge 1$ and $j\neq 0$ are defined as
\begin{equation}
        \gs_n(j) =
                \left\{
                \begin{array}{ll}
        C_{n}G_{n}(j)H_{n}(j), &\mbox{if $p(n) |j$,}\\
                  0, & \mbox{otherwise},
                \end{array}
                \right.
\end{equation}
in which
\begin{equation}
C_{n} =  \prod_{\stackrel{\scr p}{ \scr  p\neq n-1,\ p\neq n}}\left( 1 -
\frac{n-1}{(p-1)(p-n+1)}\right),
\end{equation}
\begin{equation}
G_{n}(j) = \prod_{\stackrel{\scr p|j}{ \scr  p= n-1\ {\rm or}\ p=n}}
\left(\frac{p}{p-1}\right),
\end{equation}
\begin{equation}
H_{n}(j) = \prod_{\stackrel{\scr p|j}{ \scr p\neq n-1,\ p\neq n}}
\left(1+\frac{1}{p-n}\right).
\end{equation}
Note that since $\gs(\mbox{\boldmath$j$})= \gs(\mbox{\boldmath$j$}-
\mbox{\boldmath$j_1$})$ for $\mbox{\boldmath$j_1$}$ a vector with $j_1$ in
every component, no loss of generality is incurred when the first components of
the vectors in the arguments of $\gs$ in (1.9) and (1.10) are taken to be $0$.

Gallagher \cite{GA} proved that the moments
\begin{equation}
M_k(N,h,\psi) = \sum_{n=1}^N(\psi(n+h)-\psi(n))^k \qquad (k\in \Bbb{Z}^{+})
\end{equation}
can be calculated from the prime $r$-tuple conjecture (1.7)
for $h \sim \lambda\log N$ as $N \to \infty$, 
with $\lambda$ a positive constant. For this purpose Gallagher showed that
\begin{equation}
\sum_{\stackrel{\scr 1\le j_1 , j_2 , \cdots ,j_r \le h}
{\scr  \mathrm{ distinct}}}\gs(\mbox{\boldmath$j$})\sim h^r, \qquad
(h\to\infty).  
\end{equation}
The calculation of the moments (1.16) was carried out in \cite{GY2} via
expressing them in terms of the quantities (1.7) for which
the prime $r$-tuple conjecture is assumed, with the result that
\begin{equation}
M_k(N,h,\psi) \sim N(\log N)^k \sum_{r=1}^k
 \left\{
                \begin{array}{c} k \\ r
                \end{array}
        \right\}
\lambda^r   \qquad (N \to \infty, \: h \sim \lambda\log N, \: \lambda \ll 1),
\end{equation}
where
\( \left\{
                \begin{array}{c} k \\ r
                \end{array}
        \right\}\)
denotes the Stirling numbers of the second type. 

For larger $h$ the appropriate moments to study are
\begin{equation}
\mu_{k}(N,h,\psi) = \sum_{n=1}^N(\psi(n+h)-\psi(n)-h)^k \, .
\end{equation}
Assuming the Hardy-Littlewood conjecture in the strong form
\begin{equation}
\psi_{\mbox{\boldmath$j$}}(x) = \gs(\mbox{\boldmath$j$}) x + O(N^{{1\over
2}+ \epsilon})
\end{equation}
uniformly for $1\leq r \leq k , \, 1\leq x \leq N$ and distinct $j_{i}$ 
satisfying $1\leq j_{i}\leq h$, 
Montgomery and Soundararajan \cite{MS} proved that
\begin{eqnarray}
\mu_{k}(N,h,\psi) & \sim & (1\cdot 3 \cdots (k-1)) N(h\log 
{N\over h})^{{k\over 2}} \qquad \mbox{if $k$ is even,} \\
\mu_{k}(N,h,\psi) & \ll & N(h\log N)^{{k\over 2}}({h\over \log N})^{-{1\over
8k}} + h^{k}N^{{1\over 2}+\epsilon}
\qquad \mbox{if $k$ is odd,}
\end{eqnarray}
uniformly for $(\log N)^{1+\delta} \leq h \leq N^{{1\over k}-\epsilon}$ 
(with any fixed $\delta >0$).
They also conjectured upon heuristics that $\mu_{k}(N,h,\psi) =
([2|k](1\cdot 3 \cdots (k-1)) +o(1))N(h\log {N\over h})^{{k\over 2}}$
holds uniformly for
$(\log N)^{1+\delta} \leq h \leq N^{1-\delta}$ for each fixed $k$ 
(see (1.49) below for the notation $[2|k]$).
Their proof depends on the estimation of the quantities 
\begin{equation}
R_{r}(h) = \sum_{\stackrel{\scr 1\le j_1 , j_2 , \cdots ,j_r \le h}
{\scr  \mathrm{ distinct}}}\frU((j_1 , \cdots, j_r )) ,
\end{equation}
where
\begin{equation}
\frU((j_1 , \cdots, j_r )) = \sum_{\mathcal{J} \subset \{j_1, \cdots, j_r
\}} (-1)^{r-|\mathcal{J}|}\gs(\mathcal{J})
\end{equation}
($R_{0}(h)$ and $\gs(\emptyset)$ are taken to be $1$), as
\begin{eqnarray}
\qquad \;\;
R_{r}(h) & = & (1\cdot 3 \cdots (r-1))(-h\log h + Ah)^{{r\over 2}} 
+O_{r}(h^{{r\over 2} - {1\over 7r} + \epsilon}), \quad
\mbox{$r :$ even,} \\
\qquad \;\;
R_{r}(h) & \ll & h^{{r\over 2} - {1\over 7r} + \epsilon}, \qquad
\mbox{$r :$ odd}
\end{eqnarray}
($A = 2-\gamma - \log 2\pi$, and $\gamma$ denotes Euler's constant).
Gallagher's result (1.17) can be deduced from these.
Note that it is easy to see $R_{1}(h)=0$, and 
for $r=2$ we know from Goldston \cite{FG} that (1.25) holds with the much
smaller error term $O(h^{{1\over 2}+\epsilon})$.

Only the first moment is known unconditionally as a simple
consequence of the prime number theorem. 
The work of Goldston and Montgomery \cite{GM} reveals, upon
assuming the Riemann Hypothesis, an equivalence
between the asymptotic formulae for the second moment and the pair
correlation conjecture for the zeros of the Riemann zeta-function. 

From the surrogate prime-counting function $\lambda_R(n)$, we write
\begin{equation} \psi_R(x) = \sum_{n\le x}\lambda_R(n),
\end{equation}
and we wish to examine the moments $M(N,h,\psi_R)$ defined as in (1.16).
We have
\begin{eqnarray*}
M_k(N,h,\psi_{R}) &=& \sum_{n=1}^N \left(\sum_{1\le m\le h}
\lambda_{R}(n+m)\right)^k \\
&=& \sum_{\stackrel{\scr 1\le m_i\le h }{\scr 1\le i\le k}}\sum_{n=1}^N
\lambda_{R}(n+m_1)\lambda_{R}(n+m_2)\cdots \lambda_{R}(n+m_k).
\end{eqnarray*}
Now suppose that the $k$ numbers $m_1, m_2, \ldots , m_k$ take on $r$ distinct
values
$j_1, j_2, \ldots , j_r$ with $j_i$ having multiplicity $a_i$, so that
$\sum_{1\le i \le r}a_i =k$. Grouping the terms leads to the expression
\begin{equation}
M_k(N,h,\psi_R) =  \sum_{r=1}^k \sum_{\stackrel{\scr a_1 , a_2, \ldots, a_r}
{\scr a_i \ge 1, \sum a_i =k}}
        \left(
                \begin{array}{c} k \\  a_1 , a_2 , \ldots , a_r
\end{array}
        \right)
\sum_{1\le j_1 < j_2 < \cdots < j_r \le h}\mathcal{S}_k(N,\mbox{\boldmath$j$},
\mbox{\boldmath$a$}),
\end{equation}
where $\mathcal{S}_k(N,\mbox{\boldmath$j$}, \mbox{\boldmath$a$})$ is the
correlation given in (1.4).
Our main result on these correlations is the following theorem.

\begin{theorem}  Given $1\le k\le 3$, let $\mbox{\boldmath$j$}
= (j_1,j_2, \ldots , j_r)$ and $\mbox{\boldmath$a$} = (a_1,a_2, \ldots a_r)$,
where the $j_i$'s are distinct integers, and $a_i\geq 1$ with $\sum_{i=1}^r
a_i = k$. Assume $\max_i |j_i| \leq N^{1-\epsilon}$ and 
$R\gg N^\epsilon$. Then we have
\begin{equation} \mathcal{S}_k(N,\mbox{\boldmath$j$},\mbox{\boldmath$a$}) =
\big(\mathcal{ C}_k(\mbox{\boldmath$a$})\gs(\mbox{\boldmath$j$})+o(1)\big)
N(\log R)^{k-r} +O(R^k),
\end{equation}
where $\mathcal{ C}_k(\mbox{\boldmath$a$})$ has the values
\begin{eqnarray}  \mathcal{ C}_1(1)&=&1, \quad
\mathcal{ C}_2(2) = 1, \quad  \mathcal{ C}_2(1,1)=1, \nonumber \\
 \mathcal{ C}_3(3)&=& \frac{3}{ 4}, \quad \mathcal{ C}_3(2,1)=1, \quad
\mathcal{ C}_3(1,1,1)=1 . \nonumber \end{eqnarray}
\end{theorem}

\noindent
(As a notational convention extra parentheses have been dropped,
so for example $\mathcal{ C}_2((1,1))=\mathcal{C}_2(1,1)$).
The method of proof used in this paper may be carried out for $k > 3$, but
the calculation of the constants $\mathcal{ C}_k
(\mbox{\boldmath$a$})$ and controlling the error terms
become extremely complicated even for $k=4$. 
In the third paper in this series it is shown by a different method that
Theorem 1 for $\Lambda_{R}(n)$ holds for all $k$, and also a way is
found to calculate the constants $\mathcal{ C}_{k}(\mbox{\boldmath$a$})$
for small $k$.
We also believe the error term $O(R^k)$ can be reduced in size.
In Hooley's method \cite{H13} for the special case $S_2(N,(0),(2))$,
the error term $O(R^2)$ doesn't arise at all.

Letting $m=n +\min_{i}j_{i}$ in the sum of (1.4), and then shifting
the summation range to extend from $1$ to $N$ again, we pick up an
error $O(|\min_{i}j_{i}|N^{\e})$ since $\lambda_{R}(n)\ll n^{\e}$.
This error is absorbed in the error term $o(N)$ under the 
conditions of the theorem.
Also, as was remarked after (1.15), $\gs(\mbox{\boldmath$j$})$ is
not affected by this shift. Hence
in proving Theorem 1 we may take $j_1=0$, and $j_{2},\ldots,j_{r}$
all positive.
To see the upper bound for $\lr(n)$, note that with
$\displaystyle n' = \prod_{\stackrel{\scr p|n}{\scr p\leq R}}p$,
one has $\lambda_R(n) = \lambda_R(n')$, so 
\begin{equation}
\lambda_R(n) = \sum_{r\leq R}
{\mu^{2}(r)\mu((r,n'))\phi((r,n'))\over \phi(r)} =
\sum_{\stackrel{\scr t\leq R}{\scr (t,n')
= 1}}{\mu^{2}(t)\over \phi(t)}\sum_{\stackrel{\scr s|n'}{\scr s\leq R/t}}\mu(s)
\ll d(n')\log 2R.
\end{equation}

We now apply Theorem 1 in (1.28), and obtain upon using (1.17) that
for $h\ll N^{1-\epsilon}$ and $h\to\infty$, $R=N^{\theta_{k}}$ with
$0 < \theta_k < {1\over k}$ for $M_k$,
\begin{eqnarray} 
M_1(N,h,\psi_R) & \sim & Nh , \qquad M_2(N,h,\psi_R) \sim Nh^2 + Nh\log R ,
\nonumber \\
M_3(N,h,\psi_R) & \sim & Nh^3 + 3Nh^2 \log R + {3\over 4}Nh\log^{2}R .
\end{eqnarray}
The choice $h = \lambda\log N$ renders full meaning to (1.31) allowing us to
state
\begin{corollary}
For $h\sim \lambda \log N$, $\lambda \ll 1$, and
$R=N^{\theta_k}$,  where $\theta_k$ is fixed and $0<\theta_k<\frac{1}{k}$ for
$ 1\le k \le 3$, we have
\begin{eqnarray} M_1(N,h,\psi_R)& \sim& \lambda N\log N,  \qquad
M_2(N,h,\psi_R) \sim (\theta_2 \lambda +\lambda^2)N\log^2 N, \nonumber \\
\qquad  \qquad
M_3(N,h,\psi_R)& \sim &(\frac{3}{4}{\theta_3}^2 \lambda + 3\theta_3 \lambda^2
+\lambda^3)N\log^3 N.
\end{eqnarray}
\end{corollary}

We next consider the mixed moments
\begin{equation} \tilde{M}_k(N,h,\psi_R) =
\sum_{n=1}^N(\psi_R(n+h)-\psi_R(n))^{k-1}
(\psi(n+h)-\psi(n))
\end{equation}
for $k\ge 2$, while for $k=1$ we take
$\tilde{M}_1(N,h,\psi_R)= M_1(N,h,\psi) $. Writing 
\begin{equation}
\psi(x) = x + E(x), 
\end{equation}
we have for $1\leq h\leq N$,
\begin{eqnarray}
& & M_1(N,h,\psi) = \sum_{n=1}^{N}\sum_{n<k\leq n+h}\Lambda(k) =
\sum_{k=2}^{N+h}\Lambda(k)\sum_{n=\max(1,k-h)}^{\min(k-1,N)}1
\nonumber \\
& = & \sum_{k=2}^{h}(k-1)\Lambda(k) + \sum_{k=h+1}^{N}h\Lambda(k)
+ \sum_{k=N+1}^{N+h}(N+h-k+1)\Lambda(k) \nonumber \\
& = & \psi(N+h)-\psi(N)-\psi(h) - \int_{2}^{h}\psi(t)\, dt 
+\int_{N}^{N+h}\psi(t)\, dt \nonumber \\
& = & Nh + E(N+h) - E(N) - E(h) - \int_{2}^{h}E(t)\, dt +
\int_{N}^{N+h}E(t)\, dt + O(1), \nonumber \\
& & \mbox{} 
\end{eqnarray}
where partial summation on (1.34) has also been used. The prime 
number theorem says $E(x) = o(x)$, so that we obtain for $1 \leq h\leq N$ as
$N\to\infty$,
\begin{equation}
\tilde{M}_1(N,h,\psi_R)= M_1(N,h,\psi) \sim Nh.
\end{equation}
If the Riemann Hypothesis is assumed, then it is known that
$E(x) \ll x^{1/2}\log^{2}x$, giving
\begin{equation}
\tilde{M}_1(N,h,\psi_R)= M_1(N,h,\psi) = Nh + O(N^{{1\over 2}}h\log^{2}N).
\end{equation}

For $k \ge 2$,
leaving out the details that were included in \cite{GY2}, we have 
\begin{equation} \tilde{M}_k(N,h,\psi_R) =  \sum_{r=2}^{k}
\sum_{\stackrel{\scr
a_1 , a_2, \ldots, a_{r-1}}{\scr a_i \ge 1, \sum a_i =k-1}}
\frac{1}{(r-1)!}
        \left(
                \begin{array}{c} k-1 \\  a_1 , a_2 , \ldots , a_{r-1}
      \end{array}
        \right)
W_r(N,\mbox{\boldmath$j$}, \mbox{\boldmath$a$}) + O( RN^\epsilon),
\end{equation}
where
\begin{displaymath} W_r(N,\mbox{\boldmath$j$}, \mbox{\boldmath$a$})=
\sum_{i=1}^{r-1} (\mathcal{L}_{1}(R))^{a_i}
\sum_{\stackrel{\scr 1\le j_1 ,j_2 ,\cdots ,
j_{r-1} \le h}{ \scr  \mathrm{ distinct}}} \tilde{\mathcal{ S}}_{k-a_{i}}
(N,\mbox{\boldmath$j_i$}, \mbox{\boldmath$a_i$}) +
\sum_{\stackrel{\scr 1\le j_1 ,j_2 ,\cdots , j_{r} \le h}{ \scr
\mathrm{ distinct}}}
\tilde{\mathcal{ S}}_{k}(N,\mbox{\boldmath$j$}, \mbox{\boldmath$a$}),
\end{displaymath}
\begin{displaymath} 
\mbox{\boldmath$j_i$}= (j_1,j_2, \ldots ,j_{i-1},j_{i+1},
\ldots , j_{r-1}, j_i),     \quad \mbox{\boldmath$a_i$}= (a_1,a_2,\ldots ,
a_{i-1},a_{i+1}, \ldots , a_{r-1}, 1),
\end{displaymath}
\begin{equation}
\mathcal{L}_{k}(R) :=
\sum_{\stackrel{\scr r\le R}{\scr (r,k)=1}} \frac{\mu^2(r)}{\phi(r)} .
\end{equation}
(It is easily seen that $\mathcal{L}_{1}(R) \ll \log 2R$ for $R\geq 1$.
A precise estimation of this sum due to Hildebrand \cite{Hi} is given in
(2.15) below). 

Formula (1.38) reduces the calculation of the mixed moments to the calculation
of mixed correlations. Our results depend on the extent of uniformity in the
distribution of primes in arithmetic progressions. We let
\begin{equation} \psi(x;q,a) = \sum_{\stackrel{\scr n\le x}
{ \scr n\equiv a (q)} }\Lambda(n),
\end{equation}
and on taking
\begin{equation} E(x;q,a)= \psi(x;q,a) - [(a,q)=1]\frac{x}{ \phi(q)},
\end{equation}
(where we have used the Iverson notation (1.49) below)
the estimate we need is, for some fixed $0<\vartheta \le 1$,
\begin{equation} \sum_{1\le q\le x^{\vartheta - \epsilon}}
\max_{\stackrel{\scr a }{ \scr (a,q)=1}}|E(x;q,a)|  \ll
\frac{x}{\log ^\mathcal{A}x},
\end{equation}
for any $\epsilon >0$, any $\mathcal{A}=\mathcal{A}(\epsilon)>0$, and $x$
sufficiently large (see \cite{D}, Chapter 28).
This is a weakened form of the Bombieri-Vinogradov
theorem if $\vartheta=\frac{1}{2}$, and therefore (1.42) holds unconditionally
if $\vartheta \le \frac{1}{2}$.  Elliott and Halberstam conjectured (1.42)
is true with $\vartheta=1$. The range of $R$ where our results on mixed
correlations hold depends on $\vartheta$ in (1.42). We prove

\begin{theorem}  Given $2\le k\le 3$,  let
$\mbox{\boldmath$j$} = (j_1,j_2, \ldots , j_r)$ and $\mbox{\boldmath$a$} =
(a_1,a_2, \ldots a_r)$, where $r\ge 2$, $a_r=1$, and where the $j_i$'s are
distinct integers, and $a_i\geq 1$ with $\sum_{i=1}^r a_i = k$.
Then we have, for $N^{\epsilon}\ll R \ll N^{\frac{\vartheta}{k-1} -\epsilon} $
where (1.42) holds with $\vartheta$, and 
$\max_{i}|j_{i}| \ll N^{{1\over k-1} - \epsilon}$
\begin{equation} \tilde{\mathcal{S}}_{k}(N,\mbox{\boldmath$j$},
\mbox{\boldmath$a$}) = \big(\gs(\mbox{\boldmath$j$}) +o(1)\big)N(\log R)^{k-r}.
\end{equation}
\end{theorem}

In proving Theorem 2 we may take the argument of $\Lambda$ to be $n$,
the error arising from arranging this by shifts of the range of summation
being $O(|j_{r}|N^{\e}) + O((\max |j_{i}|)^{1+\e})$.

Using (1.43) in (1.38), and (1.17), we find that 
\begin{eqnarray} 
M_2(N,h,\psi_R) & \sim & Nh^2 + Nh\log R, \nonumber \\
M_3(N,h,\psi_R) & \sim & Nh^3 + 3Nh^2 \log R + Nh\log^{2}R ,
\end{eqnarray}
and similar to Corollary 1 we have
\begin{corollary} For $h\sim \lambda \log N$, $\lambda \ll 1$,
and $R=N^{\theta_k}$,  where $\theta_k$ is fixed and
$0<\theta_k<\frac{\vartheta}{k-1}$ for $k =2\,{\mathrm or}\, 3$, we have,
\begin{eqnarray} \tilde{M}_1(N,h,\psi_R)& \sim &\lambda N\log N, \qquad
\tilde{M}_2(N,h,\psi_R) \sim (\theta_2 \lambda +\lambda^2)N\log^2 N,
\nonumber \\
\qquad \qquad
\tilde{M}_3(N,h,\psi_R)& \sim& ({\theta_3}^2 \lambda + 3\theta_3 \lambda^2
+\lambda^3)N\log^3 N.
\end{eqnarray}
\end{corollary}

The results for the correlations
up to and including the third order of $\lambda_R(n)$ and $\Lambda_R(n)$
coincide asymptotically, thereby implying the results (Theorems 1.6 and 1.7)
of \cite{GY2} on primes in short intervals, in particular
\begin{equation}
\liminf_{n\to \infty}({p_{n+r}-p_{n}\over \log p_{n}}) \leq r-{\sqrt{r}\over 2}.
\end{equation}

For longer intervals, instead of (1.42) we shall have recourse to Hooley's 
\cite{H6} bound depending on GRH that for all $q\leq x$
\begin{equation}
\sum_{\stackrel{\scr 1\leq a\leq q}{(a,q)=1}}\max_{u\leq x}|E(u;q,a)|^2
\ll x(\log x)^4  .
\end{equation}
It is easy to see that the same bound holds when the sum is taken over all 
$1\leq a \leq q$. For $(a,q)>1$ we have
\begin{displaymath}
E(u;q,a) = \psi(u;q,a) \leq \sum_{\stackrel{\scr n\leq x}{\scr n\equiv
a(\modulo q)}}\Lambda(n) ,
\end{displaymath}
where it is seen that only those $n$'s which are powers of a prime divisor
of $(a,q)$ contribute. The sum is not void only if $(a,q)$ has just one
prime factor, say $p$, in which case its value is $\leq (\log p)\lfloor
{\log x\over \log p}\rfloor \leq \log x$. So the addition of $(a,q)>1$
terms in the sum of (1.47) brings in $\leq q\log^{2}x = o(x\log^{4}x)$.

We prove
\begin{theorem} Assume the Generalized Riemann Hypothesis. For any
arbitrarily small but fixed $\eta > 0$, and for sufficiently large $N$,
with $\log^{14}N \ll h \ll N^{{1\over 7}-\epsilon}$ and writing
$h=N^{\alpha}$, there exists $n_1 ,\, n_2 \in [N+1,2N]$ such that
\begin{eqnarray}
\psi(n_1 +h) -\psi(n_1) - h & > & ({\sqrt{1-5\alpha}\over 2}-\eta) 
(h\log N)^{{1\over 2}} \nonumber \\
\psi(n_2 +h) -\psi(n_2) - h & < & -({\sqrt{1-5\alpha}\over 2}-\eta) 
(h\log N)^{{1\over 2}} .
\end{eqnarray}
\end{theorem}
This is a new development in the sense that formerly our knowledge
under GRH was restricted to lower-bound estimates for the absolute value of
the variation of the error term in the prime number theorem. The
strongest of such results were attained in \cite{GY1} in the more general
case of primes in an arithmetic progression which yielded as a special
case $\max_{x\leq y \leq 2x}|\psi(y+h)-\psi(y)-h| \gg_{\e} 
(h\log x)^{{1\over 2}}$ for $1 \leq h \leq N^{{1\over 3}-\epsilon}$. 
A proof of this is included in \S 10. In fact the general case
was also obtained by using the correlations of $\lambda_{R}(n)$. There only
the first and second level correlations were employed, nevertheless in the
more general case of $n\in [N+1, 2N]$ running through an arithmetic
progresssion $n\equiv a (\bmod \, q)$.

\bigskip

\textit{Notation.}
In this paper $N$ is always a large natural number, $p$ is a prime
number. The largest squarefree positive integer
divisor of a nonzero integer $j$ will be denoted by $j^{*}$.
If a lower limit is unspecified in a summation it will be understood
that the sum starts at $1$. When a sum is
denoted with a dash as $\sumprime_{} $ this always indicates we will sum over
all variables expressed by inequalities in the conditions of summation and
these variables will all be pairwise coprime. We will
always take the value of a void sum to be zero and the value of a void product
to be 1. The letter $\epsilon$ will denote a small positive number which may
change each time it occurs.
We will also use the Iverson notation of putting brackets around
a truth-valued statement $P(x)$ which means
\begin{equation}
        [P(x)]=
                \left\{
                \begin{array}{ll}
        1, &\mbox{if $P(x)$ is true,}\\
                  0, & \mbox{if $P(x)$ is false}.
                \end{array}
                \right.
\end{equation}
As usual, $(a,b)$ denotes the greatest common divisor of $a$ and $b$ and
$[a_1,a_2, \cdots , a_n]$ denotes the least common multiple of $a_1,a_2,
\ldots , a_n$. If $k=0$, the condition $d|k$ means $d$ can be any positive
integer; and we will take $(0,a) = 0$ for $a\neq 0$. We define
$\phi_{2}(p) = p-2$ on the primes, $\phi_{2}(1)=1$, and extend the
definition to squarefree integers multiplicatively.
For arithmetical functions
$\alpha,\, \beta$, we will sometimes write $\alpha\cdot\beta(n)$ for the
product $\alpha(n)\beta(n)$, and $\displaystyle{\alpha\over \beta}(n)$ for 
the quotient $\displaystyle{\alpha(n)\over \beta(n)}$.

\section{Lemmas}

Let us recall some well-known facts to be used in this paper.

We shall need the elementary estimates (see \cite{GY1}), for an integer
$k\neq 0$,
\begin{eqnarray}
\sum_{p|k}{\log p\over p} & \ll & \log\log 3|k| \\ 
m(k) := \sum_{d|k}\frac{\mu^2(d)}{ \sqrt{d}} & = & \prod_{p|k}
\left( 1 + \frac{1}{\sqrt{p}}\right) \ll 
\exp\left(\frac{c\sqrt{\log |k|}}{\log\log 3|k|}\right) \\
\prod_{p|k}(1+{1\over p}) & \ll & \log\log 3|k|,  \quad
\gs_{n}(k) \ll \log\log 3|k|
\end{eqnarray}
which follow from the prime number theorem 

For a multiplicative function $f(n)$ we have
\begin{equation}
\sum_{d|n}\mu^{2}(d)f(d)\log d = (\sum_{p|n}{f(p)\log p\over 1+f(p)})
\prod_{p|n}(1+f(p)) \qquad (n\neq 0).
\end{equation}

If $f\, : \mathbb{N} \to \mathbb{R}$ is a multiplicative
function satisfying $0 \leq f(p^{m}) \leq \alpha_{1}\,\alpha_{2}^{m}$
at all prime powers, with constants $\alpha_{1} > 0,\,
0 < \alpha_{2} < 2$, then we have uniformly for $x\geq 2$
\begin{equation}
\sum_{n\leq x}f(n) \ll_{\alpha_{1},\, \alpha_{2}}
{x\over \log x}\exp\sum_{p\leq x}{f(p)\over p}.
\end{equation}
This result, quoted from \cite{Hi} (which refers to \cite{HR} for the proof
of a sharper version), helps us see that for monic polynomials $P_{i}$,
$\displaystyle \prod_{p|n}{P_{1}(p)\over P_{2}(p)}$ behaves on average the
same as
$n^{\deg P_{1} - \deg P_{2}}$. In particular, we have
\begin{eqnarray}
\sum_{n\leq x}\prod_{p|n}{P_{1}(p)\over P_{2}(p)} & \ll & x
\qquad (\deg P_{1} = \deg P_{2})
\\
\sum_{n\leq x}\prod_{p|n}{P_{1}(p)\over P_{2}(p)} & \ll &
\log x \qquad (1 + \deg P_{1} = \deg P_{2})
\\
\sum_{x_{1} < n\leq x_{2}}\prod_{p|n}{P_{1}(p)\over P_{2}(p)} &
\ll & 1 + \log ({x_{2}\over x_{1}})
\qquad (1 + \deg P_{1} = \deg P_{2})
\\
\sum_{n > x}{1\over n^{\alpha}} \prod_{p|n}{P_{1}(p)\over
P_{2}(p)} & \ll_{\alpha} & {1\over x^{\alpha}}
\qquad (1 + \deg P_{1} = \deg P_{2},\; \mathrm{ fixed} \
\alpha > 0)
\\
\sum_{n > x}\log n \prod_{p|n}{1\over P_{2}(p)} & \ll &
{\log x\over x} \qquad (\deg P_{2} = 2) .
\end{eqnarray}
Here (2.6) follows from a direct application of (2.5), (2.7) and
(2.8) can be obtained by partial summation on (2.6), and to get (2.9)
one may split the sum into ranges $(x,2x],\, (2x,4x], \ldots $ and apply
(2.8) to each part. Then (2.10) is shown by partial summation on
(2.9) with $\alpha = 1$. We will also need
\begin{equation}
\sum_{n\leq x}{\mu^{2}(n)m(n)\over \sqrt{n}} \ll \sqrt{x} ,
\end{equation}
to see which we apply (2.5) with $f(n)=m(n)$ and then do partial
summation.

We also quote a Perron-formula type of result
Titchmarsh \cite{T}, \S 3.12. It will be used in proving Lemmas 2 and 3
which are needed in \S 7.
Let $\displaystyle A(s) =
\sum_{n=1}^{\infty}{a_{n}\over
n^{s}}$, for $\sigma = \Re\, s > 1$. Assume that
$a_{n} = O(\mathfrak{a}(n)),\, \mathfrak{a}(n)$ being nondecreasing, and
\begin{equation}
A(s) = \sum_{n=1}^{\infty}{|a_{n}|\over n^{\sigma}} =
O({1\over (\sigma -1)^{\alpha}}) \quad \, (\sigma \to 1^{+}).
\end{equation}
For $c>0,\, c+\sigma > 1$, we have for
$x\not\in \mathbb{N}$ and $X$ is the nearest integer to $x$
\begin{eqnarray}
\sum_{1\leq n < x}{a_{n}\over n^{s}} & = & {1\over
2\pi i}\int_{c-iT}^{c+iT}A(s+w){x^{w}\over w}\, dw
+ O({x^{c}\over T(\sigma + c - 1)^{\alpha}}) \nonumber \\
&  & \qquad + O({\mathfrak{a}(2x)x^{1-
\sigma}\log x\over T}) + O({\mathfrak{a}(X)x^{1-\sigma}\over T|x-X|}) .
\end{eqnarray}

Our first lemma is a generalization of a result of Hildebrand \cite{Hi}.

\begin{lemma} Let $P_{1}$ and $P_{2}$ be monic polynomials such
that $\deg P_{2} = 1 + \deg P_{1}$ and $P_{2}(p) \neq 0$ for prime $p$.
We have for each positive integer $k$, uniformly for $x\geq 1$,
\begin{eqnarray}
& & \sum_{\stackrel{\scr n\le x}{\scr (n,k)=1}}\mu^2(n)
\prod_{p|n}{P_{1}(p)\over P_{2}(p)} = \nonumber \\
& & \prod_{p}(1+{(p-1)P_{1}(p)-P_{2}(p)\over pP_{2}(p)})
\prod_{p|k}({P_{2}(p)\over P_{1}(p)+P_{2}(p)})
[\log x + \gamma +  \\
& & \mbox{} \;\; \sum_{p}{P_{2}(p)-(p-2)P_{1}(p)\over (p-1)(P_{1}(p)+
P_{2}(p))}\log p + \sum_{p|k}{P_{1}(p)\log p\over P_{1}(p)+P_{2}(p)}]
+ O(\frac{m(k)}{\sqrt{x}} ). \nonumber
\end{eqnarray}
\end{lemma}
Hildebrand's result is the special case
\begin{equation}
\mathcal{L}_{k}(x) =
\sum_{\stackrel{\scr n\le x}{\scr (n,k)=1}} \frac{\mu^2(n)}{\phi(n)}
= \frac{\phi(k)}{k}\left( \log x + \gamma + \sum_{p}\frac{\log p}{p(p-1)}
 + \sum_{p|k}{\log p\over p}\right) + O(\frac{m(k)}{\sqrt{x}} ) .
\label{2.15}
\end{equation}
Note for future use that, by partial summation on (2.15),
\begin{equation}
\sum_{n\leq x}{\mu^{2}(n)\over \phi(n)}\log({x\over n}) = {1\over
2}\log^{2}x + O(\log x).
\end{equation}
Another special case we will use is
\begin{eqnarray}
& & \sum_{\stackrel{\scr n\le x}{\scr (n,k)=1}}\mu^{2}(n)\prod_{p|n}
{p^2 - p-1\over (p-1)^3 } =  \nonumber \\
& &  \prod_{p}(1+{p-2\over p(p-1)^2 })
\prod_{p|k}({(p-1)^{3}\over p^3 - 2p^2 + 2p -2})
\{ \log x + \gamma + \nonumber \\
& & \mbox{} \;\; \sum_{p}{(2p-3)\log p\over (p-1)(p^3 -2p^2 +2p-2)} +
\sum_{p|k}{(p^2 -p-1) \log p\over (p^3 -2p^2 +2p-2)} \}
+ O({m(k)\over \sqrt{x}}) .
\end{eqnarray}
\begin{lemma} There exists a constant $C$ such that for all $x\geq 1$,
\begin{equation}
|\sum_{n\leq x}{\mu(n)\phi_{2}(n)\over n\phi(n)}| \leq C.
\end{equation}
\end{lemma}
\begin{lemma} As $x \to \infty$,
\begin{equation}
\sum_{n\leq x}\mu^{2}(n)
\prod_{p\mid n}{3p-4\over (p-1)(\sqrt{p}-1)} =
P(1)x^{{1\over 2}}\log^{2}x + Dx^{{1\over 2}}\log x + (E+o(1))
x^{{1\over 2}},
\end{equation}
where $P(s)$ is defined below in (3.20), and $D,\, E$ are constants
specified in (3.24).
\end{lemma}
In \S 4, \S 5 and \S 8 we will need
\begin{lemma}
For nonzero integers $j$ and $k$, we have uniformly in $x\geq 1$
\begin{equation}
\sum_{\stackrel{\scr n\leq x}{\scr (n,k)=1}}
{\mu(n)\mu\cdot\phi((n,j))\over \phi^{2}(n)} = \{1-[2\ndiv k]\mu((2,j))\}
C_{2}\prod_{\stackrel{\scr p|k}{\scr p>2}}{(p-1)^{2}\over p(p-2)}
\prod_{\stackrel{\stackrel{\scr p|j}{\scr p\ndiv k}}{\scr p>2}}
({p-1\over p-2})
+ O({d(j')j'\over x\phi(j')}) ,
\end{equation}
where $\displaystyle j' = {j^{*}\over (j^{*},k)}$.
\end{lemma}
From (2.20) we derive
\begin{equation} \begin{split}
-\sum_{n\leq x} &
{\mu(n)\,\mu\cdot\phi((n,j))\,\log n\over \phi^{2}(n)} \\
& = \left\{ \begin{array}{ll}
\gs_{2}(j)[\displaystyle\sum_{p\ndiv j}{\log p\over p(p-2)} -
\displaystyle\sum_{p|j}{\log p\over p}]
+ O({j^{*}d(j^{*})\log 2x\over \phi(j^{*})x}), &\mbox{if $2|j$,}\\
\mbox{} & \mbox{} \\
 \gs_{2}(2j){\log 2\over 2}
+ O({j^{*}d(j^{*})\log 2x\over \phi(j^{*})x}), & \mbox{if $2\ndiv j$.}\\
\end{array} \right.
\end{split}
\end{equation}

For use in \S 9 we prove
\begin{lemma}
For even $J\neq 0$ with $k|J$, and $J \ll x^{A}$\, (any fixed $A>0$),
we have uniformly in $x\geq 1$
\begin{eqnarray}
& & \qquad \qquad \sum_{n\leq x}
{\mu(n) d(n)\over \phi(n)\phi_{2}({n\over (n,2)})}
{\mu\over d}((n,J)){\mu\over \phi}((n,k))\phi_{2}(({n\over (n,2)},J)) = \\
& & 2\, [2\ndiv k]
\prod_{p\ndiv J}(1-{2\over (p-1)(p-2)})
\prod_{\stackrel{\stackrel{\scr p>2}{\scr p|J}}{\scr p\ndiv k}}
(1+{1\over p-1})
\prod_{\stackrel{\scr p>2}{\scr p|k}}
(1-{1\over (p-1)^{2}}) +  O(x^{-1+\epsilon }). \nonumber
\end{eqnarray}
\end{lemma}

\section{Proofs of the lemmas}

\emph{ Proof of Lemma 1.} We follow Hildebrand's way \cite{Hi} of
obtaining (2.15). Let, for $n\in \mathbb{N}$,
\begin{equation} f_{k}(n) =  \left\{ \begin{array}{ll}
n\mu^{2}(n)\displaystyle
\prod_{p|n}{P_{1}(p)\over P_{2}(p) }, & \mbox{ if $(n,k)=1$}, \\
\mbox{} & \mbox{} \\
0,   & \mbox{ if $(n,k)>1$ }.
\end{array} \right.
\end{equation}
Also define $g_k$, the M\"{o}bius transform of $f_k$, through
\begin{equation}
f_{k}(n) = \sum_{d|n}g_{k}(d), \qquad (n\in \mathbb{N}).
\end{equation}
We have
\begin{equation}
\sum_{\stackrel{\scr n\le x}{\scr (n,k)=1}}\!\!\mu^{2}(n)\prod_{p|n}
{P_{1}(p)\over P_{2}(p)}\!\! = \!\!\!\sum_{n\leq x}\!{f_{k}(n)\over n}\!\!
=\!\!\sum_{n\leq x}\!{g_{k}(n)\over n}\sum_{m \leq {x\over n}}\!\!{1\over m}
\!\! = \!\!\sum_{n\leq x}{g_{k}(n)\over n}(\log {x\over n} + \gamma +
O({n\over x})).
\end{equation}
The arithmetical functions $f_k$ and $g_k$ are multiplicative, their
values at the prime powers are
\begin{equation} f_{k}(p^{m}) =  \left\{ \begin{array}{ll}
{pP_{1}(p)\over P_{2}(p)}, & \mbox{ for $m=1$ and $p \ndiv k$}, \\
\mbox{} & \mbox{} \\
0,   & \mbox{otherwise },
\end{array} \right.
\end{equation}
and since $g_{k}(p^{m}) = f_{k}(p^{m})- f_{k}(p^{m-1})$,
\begin{equation} g_{k}(p^{m}) =  \left\{ \begin{array}{ll}
{pP_{1}(p)\over P_{2}(p)}-1, & \mbox{ for $m=1$ and $p \ndiv k$}, \\
\mbox{} & \mbox{} \\
{-pP_{1}(p)\over P_{2}(p)}, & \mbox{ for $m=2$ and $p \ndiv k$}, \\
\mbox{} & \mbox{} \\
-1, & \mbox{ for $m=1$ and $p|k$}, \\
\mbox{} & \mbox{} \\
0,   & \mbox{ for $m=2,\, p|k$ or $m>2$}.
\end{array} \right.
\end{equation}
We see that $\displaystyle \sum_{p,\, m\geq 1}{|g_{k}(p^{m})|\over p^{m}}$
is convergent, so that $\displaystyle\sum_{n = 1}^{\infty}{g_{k}(n)\over n}$
is absolutely convergent, giving
\begin{eqnarray}
\sum_{n=1}^{\infty}{g_{k}(n)\over n} & = & \prod_{p}(1+\sum_{m\geq 1}
{g_{k}(p^{m})\over p^m }) \nonumber \\
& = & \prod_{p|k}(1-{1\over p}) \prod_{p\ndiv k}(1+{{pP_{1}(p)\over
P_{2}(p)}-1\over p}-{P_{1}(p)\over pP_{2}(p)}) \nonumber \\
& = & \prod_{p}(1+{(p-1)P_{1}(p)-P_{2}(p)\over pP_{2}(p)})
\prod_{p|k}{P_{2}(p)\over P_{1}(p)+P_{2}(p)}.
\end{eqnarray}
Now $g_{k}(n) \neq 0$ only when $n$ is of the form $n=n_{1}n_{2}n_{3}^2 $,
with pairwise coprime $n_{i}$'s $(i=1, 2, 3)$ satisfying
$\mu^{2}(n_{i}) = 1,\, n_{1}|k,\, n_{2}n_{3}\ndiv k$, in which case
\begin{equation}
g_{k}(n) = \mu(n_{1})\mu^{2}(n_{2})\prod_{p|n_{2}}({pP_{1}(p)-P_{2}(p)\over
P_{2}(p)})\:
\mu(n_{3})n_{3}\prod_{p|n_{3}}{P_{1}(p)\over P_{2}(p)}.
\end{equation}
Hence, for $t\geq 1$, we have
\begin{eqnarray}
\sum_{n\leq t}|g_{k}(n)| & \leq & \sum_{n_{1}|k}\mu^{2}(n_{1})\sum_{n_{2}\leq
{t\over n_{1}}}\prod_{p|n_{2}}({pP_{1}(p)-P_{2}(p)\over P_{2}(p)})
\sum_{n_{3}\leq \sqrt{{t\over n_{1}n_{2}}}}n_{3}\prod_{p|n_{3}}
{P_{1}(p)\over P_{2}(p)} \nonumber \\
& \ll & \sqrt{t}\sum_{n_{1}|k}{\mu^{2}(n_{1})\over \sqrt{n_{1}}}
\sum_{n_{2}\leq {t\over n_{1}}}
{1\over \sqrt{n_{2}}}\prod_{p|n_{2}}{pP_{1}(p)-P_{2}(p)\over P_{2}(p)} \ll
\sqrt{t}\, m(k),
\end{eqnarray}
where we have made use of (2.6) (resp. (2.9)) for the sum
over $n_{3}$ (resp. $n_{2}$). Next, for $u\geq 1$, we have
\begin{equation}
\left| \sum_{n > u}{g_{k}(n)\over n}\right| = \left| \int_{u}^{\infty}
{1\over t^2 }\sum_{u<n \leq t}g_{k}(n)\, dt \right| \ll m(k)\int_{u}^{\infty}
t^{-3/2}\, dt \ll {m(k)\over \sqrt{u}},
\end{equation}
and therefore
\begin{equation}
\sum_{n \leq x}{g_{k}(n)\over n} = \sum_{n=1}^{\infty}{g_{k}(n)\over n} +
O({m(k)\over \sqrt{x}}) .
\end{equation}
For the main term observe that
\begin{eqnarray}
& & \sum_{n\leq x}{g_{k}(n)\over n}\log {x\over n}  =  \int_{1}^{x}
{1\over u}\sum_{n\leq u}{g_{k}(n)\over n}\, du \nonumber \\
& & \mbox{} =
\sum_{n=1}^{\infty}{g_{k}(n)\over n}\int_{1}^{x}{du\over
u}-\int_{1}^{\infty}
{1\over u}\sum_{n > u}{g_{k}(n)\over n}\, du + \int_{x}^{\infty}{1\over u}
\sum_{n > u}{g_{k}(n)\over n}\, du \nonumber \\
& & \mbox{} = \sum_{n=1}^{\infty}{g_{k}(n)\over n}\log x
-\int_{1}^{\infty}
{1\over u}\sum_{n > u}{g_{k}(n)\over n}\, du + O({m(k)\over \sqrt{x}}) .
\end{eqnarray}
The last integral here is
\begin{eqnarray}
& & \int_{1}^{\infty}{1\over u}\sum_{n > u}{g_{k}(n)\over n}\, du  =
\sum_{n=1}^{\infty}{g_{k}(n)\log n\over n} =
\sum_{n=1}^{\infty}{g_{k}(n)\over n}
\sum_{p^{m}||n}\log p^m \nonumber \\
& & \mbox{} =
\sum_{p,\, m\geq 1}\log p^m \sum_{\stackrel{\scr n=1}{\scr
p^{m}||n}}^{\infty}{g_{k}(n)\over n} =
\sum_{p,\, m\geq 1}{g_{k}(p^{m})\log p^m \over p^m }
\sum_{\stackrel{\scr n=1}{\scr p\ndiv n}}^{\infty}{g_{k}(n)\over n}
\nonumber \\
& & \mbox{} = \sum_{n=1}^{\infty}{g_{k}(n)\over n}\sum_{p,\, m\geq 1}
{g_{k}(p^{m})\log p^m \over p^m }(1+\sum_{l=1}^{\infty}{g_{k}
(p^{l})\over p^{l}})^{-1} \nonumber \\
& & \mbox{} =
\sum_{n=1}^{\infty}{g_{k}(n)\over n}\{
\sum_{p\ndiv k}{[(p-2)P_{1}(p)-P_{2}(p)]\log p\over
(p-1)(P_{1}(p)+P_{2}(p))} - \sum_{p|k}{\log p\over p-1}\} \nonumber \\
& & \mbox{} = \sum_{n=1}^{\infty}{g_{k}(n)\over n}
\{ \sum_{p}{[(p-2)P_{1}(p)-P_{2}(p)]\log p\over (p-1)(P_{1}(p)+P_{2}(p))} 
-\sum_{p|k}{P_{1}(p)\log p\over (P_{1}(p)+P_{2}(p))} \} .
\end{eqnarray}
Plugging (3.12) in (3.11), and then using (3.11), (3.10), (3.8) and (3.6) in
(3.3), we complete the proof of Lemma 1.

\emph{ Proof of Lemma 2. }  Let
\begin{equation}
A(s) := \sum_{n=1}^{\infty}{\mu(n)\phi_{2}(n)\over \phi(n)}{1\over n^{s}},
\end{equation}
the series being absolutely convergent for $\Re\, s > 1$. From (2.13)
with $x$ half an odd integer, we have
\begin{equation}
\sum_{n\leq x}{\mu(n)\phi_{2}(n)\over n\phi(n)}= {1\over 2\pi i}
\int_{c-iT}^{c+iT}A(1+w){x^{w}\over w}\, dw
\, + O({x^{c}\over cT}) + O({\log x\over T}).
\end{equation}
Taking $\displaystyle c={1\over \log x}$, which minimizes
$\displaystyle {x^{c}\over c}$, the error of (3.14) is
$\displaystyle O({\log x\over T})$. Next note that
\begin{equation}
A(s) = \prod_{p}(1-{p-2\over (p-1)p^s}) = {1\over \zeta(s)}\prod_{p}(1+{1\over
(p-1)(p^{s}-1)}),
\end{equation}
where the last product is expressible as a Dirichlet series which is
absolutely convergent for $\Re\, s > 0$. So we have
\begin{equation}
\sum_{n\leq x}{\mu(n)\phi_{2}(n)\over n\phi(n)} =
{1\over 2\pi i}
\int_{c-iT}^{c+iT}{1\over \zeta(w+1)}\prod_{p}(1+{1\over (p-1)(p^{1+w}-1)})
{x^{w}\over w}\, dw
+O({\log x\over T}).
\end{equation}
Now we pull the line of integration to $w = -{K\over \log T}+ it,\,
-T \leq t \leq T$ in accordance with the well-known zero-free region
for $\zeta(s)$, so that the integrand has no poles in the region thus formed.
Here
\begin{equation}
{1\over \zeta(s)} = O(\log (|t|+2)), \qquad
(\sigma \geq 1-{K\over \log (|t|+2)})
\end{equation}
holds (see Titchmarsh \cite{T}, Thm. 3.8 and Eq. (3.11.8)), so that
\begin{eqnarray}
\int_{-{K\over \log T} -iT}^{-{K\over \log T} + iT} & &
{1\over \zeta(1+w)}\prod_{p}(1+{1\over (p-1)(p^{w+1}-1)})
{x^{w}\over w}\, dw  \nonumber \\
& & \,\, \ll x^{-{K\over \log T}}\log T \int_{-T}^{T}{dt\over
\sqrt{{K\over \log T}^{2}
+ t^2}} \ll x^{-{K\over \log T}}\log^2 T.
\end{eqnarray}
For the integrals over the horizontal sides of the contour we have
\begin{eqnarray}
\int_{{-K\over \log T}}^{1\over \log x}& &
{1\over \zeta(1+u + iT)}\prod_{p}(1+{1\over (p-1)(p^{1+u + iT}-1)})
{x^{u + iT}\over (u+ iT)}\, du  \nonumber \\
& & \,\, \ll {\log T\over T}\int_{{-K\over \log T}}^{1\over \log x}
x^{u}\, du \ll {\log T\over T}.
\end{eqnarray}
By taking $\log T = \sqrt{\log x}$, all the error terms in
(3.16), (3.18) and (3.19) are made to tend to $0$ as $x\to \infty$. So, as
$x\to \infty$, the sum of (2.18) tends to $0$, and this
completes the proof of Lemma 2.

\emph{ Proof of Lemma 3.} Consider
\begin{eqnarray}
A(s) & := & \sum_{n=1}^{\infty}{\mu^{2}(n)\over n^s }
\prod_{p\mid n}{(3p-4)\sqrt{p}\over (p-1)(\sqrt{p} -1)}
= \prod_{p}(1+{3\over p^s}{(p-{4\over 3})\sqrt{p}\over (p-1)(\sqrt{p}-1)})
\nonumber \\
& = & \zeta^{3}(s)\prod_{p}[1+{3\over p^{s+{1\over 2}}}
({p-{\sqrt{p}\over 3}-1\over p-\sqrt{p}-1+{1\over \sqrt{p}}})
+{3\over p^{2s}}
(1-{(3p-4)\sqrt{p}\over (p-1)(\sqrt{p}-1)})
\nonumber \\
& & \qquad + {1\over p^{3s}}
({(9p-12)\sqrt{p}\over (p-1)(\sqrt{p}-1)}-1)
-{1\over p^{4s}}{(3p-4)\sqrt{p}\over (p-1)(\sqrt{p}-1)}]
\end{eqnarray}
The series $A(s)$ converges absolutely for $\Re\, s > 1$, and
the last product, call it P(s), is absolutely convergent for
$\Re\, s > {1\over 2}$. Eq. (2.13) can be applied with $\alpha = 3$,
$x$ half an odd integer, and $c= {1\over 2}+{1\over \log x}$
to have
\begin{eqnarray}
& & \sum_{n\leq x}\mu^{2}(n)
\prod_{p\mid n}{(3p-4)\over (p-1)(\sqrt{p}-1)} \nonumber \\
& & \qquad =
{1\over 2\pi i}\int_{c-iT}^{c+iT}A({1\over 2}+w){x^{w}\over w}\, dw
\, + O({x^{{1\over 2}}\log^{3} x\over T}) +
O({e^{{K\log x\over \log\log x}}x^{{1\over 2}}\log x\over T}).
\end{eqnarray}
In writing the very last error term (with an appropriate constant $K$)
we have used
\begin{equation}
\prod_{p\mid n}{(3p-4)\sqrt{p}\over (p-1)(\sqrt{p}-1)}
\ll 4^{\omega(n)}m(n),
\end{equation}
and elementary deductions from the prime number theorem, namely (2.2) and
\begin{equation}
\omega(n) \leq c_{1}{\log n\over \log\log n} .
\end{equation}
Now we pull the line of integration to $\Re w = \delta,\,
0 < \delta < {1\over 2}$. In doing so we pass the triple
pole of the integrand at $w={1\over 2}$, where the residue is
\begin{equation}
P(1)x^{{1\over 2}}\log^{2}x + Dx^{{1\over 2}}\log x + Ex^{{1\over 2}}
\end{equation}
(with the constants $D,\, E$ made up of the Stieltjes constants and
the values of $P(s)$ and its first two derivatives at $s=1$).
On the left vertical side of the contour we will have
\begin{eqnarray}
\int_{\delta -iT}^{\delta + iT}
\zeta^{3}(w+{1\over 2})P(w+{1\over 2}){x^{w}\over w}\, dw
& \ll & x^{\delta}\int_{0}^{T}{|\zeta({1\over 2}+\delta + it)|^{3}\over
\sqrt{\delta^{2}+t^2}}\, dt \nonumber \\
& \ll & x^{\delta}({1\over ({1\over 2}-\delta)^{2}}+ \log^{3}T) ,
\end{eqnarray}
if we take
${1\over 2} -{U\over \log^{{2\over 3}}T}< \delta < {1\over 2}$, $U$ being
an appropriate constant.
To see this we employ the estimate
(see Karatsuba and Voronin \cite{KV}, p.116)
\begin{equation}
\zeta(\sigma + it) = O(\log^{{2\over 3}}|t|), \quad (\sigma \geq 1-
{U\over \log^{{2\over 3}}|t|},\, |t|\geq 2),
\end{equation}
which implies
\begin{equation}
x^{\delta}\int_{2}^{T}{|\zeta({1\over 2}+\delta + it)|^{3}\over
\sqrt{\delta^{2}+t^{2}}}\, dt \ll x^{\delta}\int_{2}^{T}{\log^{2}t\over t}\, dt
\ll x^{\delta}\log^{3}T .
\end{equation}
We also have
\begin{equation}
x^{\delta}\int_{0}^{2}{|\zeta({1\over 2}+\delta + it)|^{3}\over
\sqrt{\delta^{2}+t^2}}\, dt \ll x^{\delta}\int_{0}^{2}{1\over |\delta -
{1\over 2} + it|^3}\,dt \ll {x^{\delta}\over ({1\over 2}-\delta)^2}.
\end{equation}
For the horizontal sides of the contour we have, by (3.26),
\begin{equation}
\int_{\delta +iT}^{c+ iT}\zeta^{3}(w +{1\over 2})P(w+{1\over 2}){x^{w}\over w}
\, dw \ll {1\over T}\int_{\delta}^{c}|\zeta({1\over 2}+\sigma + iT)|^{3}
x^{\sigma}\, d\sigma \ll {x^{{1\over 2}}\log^{2}T\over T\log x}.
\end{equation}
Choosing
\begin{equation}
\delta = {1\over 2} - {1\over (\log T)^{{3\over 4}}}, \quad
T= e^{{2K\log x\over \log\log x}}\log x \log\log x
\end{equation}
we make all the error terms in (3.21), (3.25) and (3.29) to be $o(\sqrt{x})$.
This completes the proof of Lemma 3.

\emph{Proof of Lemma 4.} We have
\begin{equation}
\sum_{\stackrel{\scr n = 1}{\scr (n,k)=1}}^{\infty}
{\mu(n)\mu\cdot\phi((n,j))\over \phi^{2}(n)} =
\prod_{\stackrel{\scr p|j}{\scr p\ndiv k}}(1+{1\over p-1})
\prod_{\stackrel{\scr p\ndiv j}{\scr p\ndiv k}}(1-{1\over (p-1)^{2}}) ,
\end{equation}
where we notice that the last product is $0$ if $2\ndiv jk$. If $2|jk$,
the products are re-organized to give the main term of (2.20),
which has been expressed so as to be valid whether or not $2|jk$.
The $O$-term of (2.20) is the tail of the series 
\begin{eqnarray}
\sum_{\stackrel{\scr n >x}{\scr (n,k)=1}}
{\mu(n)\mu\cdot\phi((n,j))\over \phi^{2}(n)} & = &
\sum_{\stackrel{\scr d|j^{*}}{(d,k)=1}}\mu\cdot\phi(d)
\sum_{\stackrel{\stackrel{\scr n>x}{\scr (n,k)=1}}{\scr (n,j)=d}}
{\mu(n)\over \phi^{2}(n)} \nonumber \\
= & & \!\!\!\! \sum_{\stackrel{\scr
d|j^{*}}{(d,k)=1}}{\mu^{2}(d)\over \phi(d)} \sum_{\stackrel{\scr
t>{x\over d}}{(t,kj)=1}}{\mu(t)\over \phi^{2}(t)} \ll {1\over x}
\sum_{\stackrel{\scr d|j^{*}}{(d,k)=1}}{d\over \phi(d)} \ll
{d(j')j'\over x\phi(j')} .
\end{eqnarray}
To obtain (2.21), let
\begin{equation}
B(s) := \sum_{n = 1}^{\infty}
{\mu(n)\mu\cdot\phi((n,j))\over \phi^{2}(n)n^{s}} =
\prod_{p|j}(1+{1\over (p-1)p^{s}})
\prod_{p\ndiv j}(1-{1\over (p-1)^{2}p^{s}}) ,
\end{equation}
where for $\Re\, s > -1$ the product is absolutely convergent. Then
\begin{equation}
        B'(0) =
                \left\{
                \begin{array}{ll}
        \gs_{2}(j)[\displaystyle\sum_{p\ndiv j}{\log p\over p(p-2)} -
\displaystyle\sum_{p|j}{\log p\over p}], &\mbox{if $2|j$,}\\
\mbox{} & \mbox{} \\
                  \gs_{2}(2j){\log 2\over 2}, & \mbox{if $2\ndiv j$ .}
                \end{array}
                \right.
\end{equation}
But
\begin{displaymath}
\sum_{n \leq x}
{\mu(n)\mu\cdot\phi((n,j))\over \phi^{2}(n)}(-\log n)
= B'(0) +
\sum_{n > x}
{\mu(n)\mu\cdot\phi((n,j))\over \phi^{2}(n)}(\log n) ,
\end{displaymath}
and the very last sum is shown to be
$\ll {d(j^{*})j^{*}\log 2x\over \phi(j^{*})x}$ similarly to (3.32).

\emph{Proof of Lemma 5.} 
We extend the sum to all $n$, introducing as error the
tail of the series for $n>x$. We have
\begin{eqnarray}
& & \qquad \sum_{n=1}^{\infty}
{\mu(n) d(n)\over \phi(n)\phi_{2}({n\over (n,2)})}
{\mu\over d}((n,J)){\mu\over \phi}((n,k))\phi_{2}(({n\over (n,2)},J)) \\
& & \qquad = (1+\mu((2,k)))
\prod_{p\ndiv J}(1-{2\over (p-1)(p-2)})
\prod_{\stackrel{\stackrel{\scr p>2}{\scr p|J}}{\scr p\ndiv k}}
(1+{1\over p-1})
\prod_{\stackrel{\scr p>2}{\scr p|k}}
(1-{1\over (p-1)^{2}}), \nonumber
\end{eqnarray}
and the introduced error is bounded as
\begin{eqnarray}
& \ll & \sum_{n >x}
{\mu^{2}(n) d(n)\phi_{2}(({n\over (n,2)},J))\over 
\phi(n)\phi_{2}({n\over (n,2)})d((n,J))}
\ll \sum_{m|J}{\phi_{2}\over d}(m) 
\sum_{\stackrel{\scr n>x}{\scr (n,{J\over 2})=m}}
{\mu^{2}(n) d(n)\over \phi(n)\phi_{2}({n\over (n,2)})} \nonumber \\
& \ll & \sum_{m|J}{1\over \phi(m)}
\sum_{\stackrel{\scr t\geq {x\over m}}{\scr (t,2J)=1}}
{\mu^{2}(t) d(t)\over \phi(t)\phi_{2}(t)} \ll
\sum_{m|J}{1\over \phi(m)}\sum_{t\geq {x\over m}}{1\over t^{2-\epsilon}}
\nonumber \\
& \ll & {1\over x^{1-\epsilon}}\sum_{m|J}{m^{1-\epsilon}\over \phi(m)}
\ll {1\over x^{1-\epsilon}}.
\end{eqnarray}

\section{ Pair Correlations of $\lambda_R(n)$}
In the case $k=1$ of Theorem 1 we have
\begin{eqnarray} \mathcal{ S}_1(N,(j_{1}),(1)) &=& \sum_{n\le
N}\lambda_R(n+j_{1})=
\sum_{r\le R}{\mu^{2}(r)\over \phi(r)}\sum_{d|r}d\mu(d)\sum_{\stackrel{\scr
\max(1,1+j_{1})\leq n \leq N+j_{1}}{\scr d|n}}1 \nonumber \\
&=& \sum_{r\le R}{\mu^{2}(r)\over \phi(r)}\sum_{d|r}d\mu(d)
({N+ j_{1}-\max(0,j_{1})\over d}+O(1))
\nonumber \\
&=& N +\min(0,j_{1})
+ O(\sum_{r\le R}{\mu^{2}(r)\sigma(r)\over \phi(r)}) = N(1+o(1)) + O(R),
\end{eqnarray}
where we refer to (2.6) for
\begin{equation}
\sum_{r\le R}{\mu^{2}(r)\sigma(r)\over \phi(r)} \ll R .
\end{equation}

To examine the case $k=2$ of Theorem 1 we need to consider
\begin{equation} \mathcal{ S}_2(j)= \sum_{n=1}^N \lambda_R(n)\lambda_R(n+j).
\end{equation}
In our earlier notation, $\mathcal{ S}_2(j)= \mathcal{ S}_2(N, (0,j),(1,1))$
if $j\neq 0$, and $\mathcal{ S}_2(0) = \mathcal{ S}_2(N,(0),(2))$. We have
for any $j$,
\begin{equation}
\mathcal{ S}_2(j) = \sum_{r_{1}, r_{2} \leq R}{\mu^{2}(r_1)
\mu^{2}(r_2)\over \phi(r_1)\phi(r_2)}\sum_{\stackrel{\scr d|r_1 }{\scr
e|r_2 }}d\mu(d)e\mu(e)
\sum_{\stackrel{\stackrel{\scr n\leq N}{\scr d|n}}{\scr e|n+j}}1 .
\label{4.4}
\end{equation}
The innermost sum is over $n$'s in a unique residue class modulo
$[d,e]$ whenever $(d,e)|j$, in which case its value is $N/[d,e] + O(1)$,
otherwise the innermost sum is void. By (4.2) the last $O(1)$ leads to a
contribution of $O(R^2)$ in (4.4). Hence
\begin{equation}
\mathcal{ S}_2(j) = N\sum_{r_{1}, r_{2} \leq R}{\mu^{2}(r_1)
\mu^{2}(r_2)\over \phi(r_1)\phi(r_2)}
\sum_{\stackrel{\stackrel{\scr d|r_1 }{\scr e|r_2 }}{\scr (d,e)|j}}
\mu(d)\mu(e)(d,e) + O(R^2). \label{4.5}
\end{equation}
Let $(d,e) = \delta,\, d=d'\delta,\, e=e'\delta$ so that $(d',e')=1$. Then the
inner sums over $d$ and $e$ become
\begin{equation*}
\sum_{\stackrel{\scr \delta|(r_{1},r_{2})}{\scr \delta|j}}\delta
\sum_{d'|{r{_1}\over \delta}}\mu(d')
\sum_{\stackrel{\scr e'|{r_{2}\over \delta}}{\scr (e',d')=1}}\mu(e').
\end{equation*}
Here the innermost sum is
\begin{equation*}
\sum_{\stackrel{e'|{r_{2}\over \delta}}{(e',d')=1}}\mu(e') =
\prod_{\substack{p|{r_{2}\over \delta}\\ p\ndiv d'}}(1+\mu(p)) =
                \left\{
                \begin{array}{ll}
        1, &\mbox{if ${r_{2}\over \delta}|d'$ ,}\\
                  0, & \mbox{otherwise}.
                \end{array}
                \right.
\end{equation*}
Next the sum over $d'$ becomes
\begin{equation*}
\sum_{\stackrel{\scr d'|{r_{1}\over \delta}}{\scr {r_{2}\over \delta}|d'}}
\mu(d') =
                \left\{
                \begin{array}{ll}
        \mu({r_{1}\over \delta}), &\mbox{if $r_{1} = r_{2}$ ,}\\
                  0, & \mbox{otherwise}.
                \end{array}
                \right.
\end{equation*}
Hence we have, for any $j$,
\begin{eqnarray}
\mathcal{ S}_2(j) & = & N\sum_{r_{1}\leq R}{\mu(r_1)\over \phi^{2}(r_1)}
\sum_{\delta|(r_{1},j)}\delta\mu(\delta) + O(R^2) \nonumber \\
& = & N\sum_{r_{1}\leq R}{\mu(r_1)\mu((j,r_{1}))\phi((j,r_{1}))\over
\phi^{2}(r_1)} + O(R^2) . 
\end{eqnarray}
Now, if $j=0$, (4.6) reduces to
\begin{equation}
\sum_{n\leq N}(\lr(n))^2 = N\mathcal{L}_{1}(R) + O(R^2), 
\end{equation}
and by (2.15) this proves Theorem 1 for the case $\mathcal{S}_2(N,(0),(2))$.
If $j \neq 0$, then by Lemma 4
\begin{equation}
\sum_{r_{1}\leq R}{\mu(r_1)\mu((j,r_{1}))\phi((j,r_{1}))\over \phi^{2}(r_1)} 
= \gs_{2}(j) + O({j^{*}d(j^{*})\over R\phi(j^{*})}) ,
\end{equation}
and therefore
\begin{equation}
\sum_{n\leq N}\lr(n)\lr(n+j)  =  N\gs_{2}(j) +
O\left({Nj^{*}d(j^{*})\over R\phi(j^{*})}\right) + O(R^2), \qquad (j\neq 0),
\end{equation}
which completes the proof of Theorem 1 for the case
$\mathcal{ S}_2(N, (0,j),(1,1))$.

\section{The Mixed Correlations}
We now turn our attention to the mixed correlations
$\tilde{\mathcal{S}}_k(N,\mbox{\boldmath$j$}, \mbox{\boldmath$a$})$ defined in
(1.5). The $k=1$ case was noted in (1.6), and the first mixed moment
was treated in (1.35)-(1.37). We consider
\begin{equation}
\tilde{\mathcal{S}}_2(j) = \sum_{n=1}^N \lambda_R(n+j)\Lambda(n), \qquad
(j\neq 0)
\end{equation}
in the case of mixed second level correlations.
In the notation of (1.5), $\tilde{\mathcal{S}}_2(j) =
\tilde{\mathcal{S}}_2(N, (j,0), (1,1))$. We have
\begin{equation}
\tilde{\mathcal{S}}_2(j) =
\sum_{r\le R}{\mu^{2}(r)\over \phi(r)}\sum_{d|r}d\mu(d)\sum_{\stackrel{\scr
\max(2,1-j) \leq n \leq N}{\scr d|n+j}}\Lambda(n), \qquad (j\neq 0) .
\end{equation}
The innermost sum of (5.2) is 
\begin{eqnarray}
& = & \psi(N;d,-j) - \psi(\max(2,1-j);d,-j) \\
& = & [(d,j)=1]{N\over \phi(d)}+
E(N;d,-j) + [j<0]\cdot O((1+{|j|\over d})\log |j|) , \nonumber
\end{eqnarray}
by (1.41). Hence (5.2) becomes
\begin{eqnarray}
\tilde{\mathcal{S}}_2(j) & = & N
\sum_{r\le R}{\mu^{2}(r)\over \phi(r)}\sum_{\stackrel{\scr d|r}{\scr (d,j)=1}}
{d\mu(d)\over \phi(d)} 
+ O\bigl(\sum_{d\leq R}{d\mu^{2}(d)\over \phi(d)}|E(N;d,-j)| 
\sum_{\stackrel{\scr s\le {R\over d}}{\scr (s,d)=1}}{\mu^{2}(s)\over \phi(s)}
\bigr) \nonumber \\
& & \mbox{} + [j<0]\cdot O\bigl(\log |j|\sum_{r\leq R}{\mu^{2}(r)\over
\phi(r)}\sum_{d|r}d(1+{|j|\over d})\bigr).
\end{eqnarray}
Since
\begin{equation}
\sum_{\stackrel{\scr d|r}{\scr (d,j)=1}}{d\mu(d)\over \phi(d)} =
{\mu({r\over (r,j)})\over \phi({r\over (r,j)})} , 
\end{equation}
the main term of (5.4) is the same as that of (4.6), so (4.8) settles it.
The last error term is easily bounded as 
$O((R+|j|\log^{2}R)\log |j|)$. As for the
error term with $E$ in (5.4), it is
\begin{displaymath}
\ll  \log R \log\log R (\sum_{\stackrel{\scr d\leq R}{\scr (d,j)=1}}
|E(N;d,-j)| + \sum_{\stackrel{\scr d\leq R}{\scr (d,j)>1}}\mu^{2}(d)
\psi(N;d,-j) ) .
\end{displaymath}
Here the first sum is estimated by the Bombieri-Vinogradov theorem (1.42),
provided that $R \ll N^{1/2 - \epsilon}$. In the second sum only 
$n$'s which are powers of those primes that are divisors of
$(d,j)$ contribute to $\psi$. Hence the second sum is
\begin{displaymath}
\ll \sum_{p|j}\sum_{\stackrel{\scr d\leq R}{\scr p|d}}\sum_{\stackrel
{\scr p^{a}\leq N}{\scr d|p^{a}+j}}\log p  \ll \log N\sum_{p|j}
\sum_{\stackrel{\scr d\leq R}{\scr p|d}}1 \ll R\log N \log\log 3|j| .
\end{displaymath}
Thus we obtain, for $j\neq 0$,
\begin{equation}
\tilde{\mathcal{S}}_2(j) = \sum_{n=1}^N \lambda_R(n+j)\Lambda(n) =
N\gs_{2}(j) + O\left({Nj^{*}d(j^{*})\over R\phi(j^{*})}\right)
+ O(\frac{N}{\log^{\mathcal{A}}N}) ,
\end{equation}
completing the proof of Theorem 2 in the case $k=2$.

For mixed correlations of the third level we may begin with
\begin{equation}
\tilde{\mathcal{S}}_3(j_{1},j_{2},0) = \sum_{n=1}^N \lr(n+j_{1})\lr(n+j_{2})
\Lambda(n), \qquad (j_{1}j_{2} \neq 0). 
\end{equation}
In the beginning
the two cases, $j_1 = j_2$ or not, will undergo a common treatment.
Using the definition of $\lr(n)$ we can rewrite this as
\begin{equation}
\tilde{\mathcal{S}}_3(j_{1},j_{2},0) =
\sum_{r_{1}, r_{2}\le R}{\mu^{2}(r_{1})\mu^{2}(r_{2})
\over \phi(r_{1})\phi(r_{2})}
\sum_{\stackrel{\scr d|r_{1}}{\scr e|r_{2}}}d\mu(d)e\mu(e)
\sum_{\stackrel{\stackrel{\scr \max(2,1-j_{1},1-j_{2}) \leq n \leq N}
{\scr n\equiv -j_{1}(\mathrm{mod}\, d)}}
{\scr n\equiv -j_{2}(\mathrm{mod}\, e)}}
\Lambda(n) . 
\end{equation}
The congruence conditions in the innermost sum are compatible if and only if
$(d,e)|j_1 - j_2$, in which case there is a unique residue class $j$ such
that $n\equiv j(\modulo [d,e])$, and similar to (5.3) the innermost sum is
\begin{equation}
= [([d,e],j)=1]{N\over \phi([d,e])}
+ E(N;[d,e],j) + O(\max_{i}|j_{i}|\log N).
\end{equation}
Placed in (5.8), the error part of (5.9) contributes, (calling the 
last $O$-term $M$)
\begin{eqnarray}
& \ll &
\sum_{\stackrel{\scr d,e \leq R}{\scr (d,e)|j_1 - j_2 }}
{de\over \phi(d)\phi(e)}
[|E(N;[d,e],j)| + O(M)]
\sum_{\stackrel{\scr s_1 \le {R\over d}}{\scr (s_{1},d)=1}}
{\mu^{2}(s_{1})\over \phi(s_{1})}
\sum_{\stackrel{\scr s_2 \le {R\over e}}{\scr (s_{2},d)=1}}
{\mu^{2}(s_{2})\over \phi(s_{2})} \nonumber \\
& \ll & (\log R)^2 (\log\log R)^2
\sum_{\stackrel{\scr d,e \leq R}{\scr (d,e)|j_1 - j_2}}
[|E(N;[d,e],j)| + O(M)]  \nonumber \\
& \ll & (\log R)^3 \sum_{D\leq R^2 }
[\max_{a(\modulo D)}
|E(N;D,a)| + O(M)]\sum_{\stackrel{\scr d,e \leq R}{\scr [d,e]
= D}}1  \nonumber \\
& \ll & (\log R)^3 \sum_{D\leq R^2 }
[\max_{a(\modulo D)}|E(N;D,a)| + O(M)]d_{3}(D) \nonumber \\
& \ll & (\log R)^3
\sqrt{\sum_{D\leq R^2 }{d_{3}(D)^{2}\over D}}\sqrt{\sum_{D\leq R^2 }
D[\max_{a(\modulo D)}|E(N;D,a)| + O(M)]^2} \nonumber \\
& \ll & (\log R)^{{15\over 2}} \sqrt{\sum_{D\leq R^2 }
D[{N\log N\over D}\max_{a(\modulo D)}|E(N;D,a)| + O(M^{2})]}\nonumber \\
& \ll & (\log R)^{{15\over 2}} \sqrt{{N^{2}\over \log^{\mathcal{A} - 1}N} +
NR^{2}\log^{2}N + (\max_{i}|j_{i}|)^{2} R^{4}\log^{2}N }
\nonumber \\
& \ll & {N\over \log^{\mathcal{B}}N}, \nonumber \\
& \, & \mbox{}
\end{eqnarray}
where $\mathcal{B}$ is as large as wished by taking $\mathcal{A}$ large
enough,
provided that $R \ll N^{1/4 - \epsilon}$ for the applicability of the
Bombieri-Vinogradov estimate (1.42), and that $R^{2}\max_{i}|j_{i}| \ll
N^{1-\e}$.
In the above sequence of inequalities
we have used some well-known features of the divisor function $d_{3}(n)$,
and we have employed the trivial estimate $\displaystyle |E(N;D,a)| \ll
{N\log N\over D}$, and the $(a,D)>1$ part of 
$\displaystyle\max_{a(\modulo D)}$
is estimated to be $\ll NR^{2}\log^{2}N$ as was done for (5.6).

From (5.8)--(5.10) we have
\begin{equation}
\tilde{\mathcal{S}}_3(j_{1},j_{2},0) = N
\sum_{r_{1}, r_{2}\le R}
{\mu^{2}(r_{1})\mu^{2}(r_{2})\over \phi(r_{1})\phi(r_{2})}
\sum_{\stackrel{\stackrel{\scr d|r_{1}, e|r_2 }{\scr (d,j_{1})=1, (e,j_{2})=1}}
{\scr (d,e)|j_{1}-j_{2}}}{d\mu(d)e\mu(e)\over \phi([d,e])} +
O({N\over \log^{\mathcal{B}}N}),
\end{equation}
for $(j,[d,e]) = 1$ if and only if $(j_{1},d)=1$ and $(j_{2},e)=1$.
Letting
\begin{equation}
r_{1}' = {r_{1}\over (r_{1},j_{1})}, \;\; r_{2}' = {r_{2}\over (r_{2},j_{2})},
\;\; (d,e)= \delta, \;\; d=d'\delta, \;\; e =e'\delta,
\end{equation}
the inner sums over $d$ and $e$ in (5.11) take the form
\begin{equation*}
\sum_{\stackrel{\scr \delta|(r_{1}',r_{2}')}{\scr \delta|j_{1}-j_{2} }}
{\delta^{2}\over \phi(\delta)}
\sum_{d'|{r{_1}'\over \delta}}{d'\mu(d')\over \phi(d')}
\sum_{\stackrel{\scr e'|{r_{2}'\over \delta}}{\scr (e',d')=1}}
{e'\mu(e')\over \phi(e')}.
\end{equation*}
Here the innermost sum is
\begin{equation*}
\sum_{\stackrel{\scr e'|{r_{2}'\over \delta}}{\scr (e',d')=1}}
{e'\mu(e')\over \phi(e')} =
\prod_{\substack{p|{r_{2}'\over \delta}\\ p\ndiv d'}}(1-{p\over p-1}) =
\prod_{\substack{p|{r_{2}'\over \delta}\\ p\ndiv d'}}{-1\over p-1} =
{\mu({{r_{2}'\over \delta}\over ({r_{2}'\over \delta},d')})
\over \phi({{r_{2}'\over \delta}\over ({r_{2}'\over \delta},d')})}.
\end{equation*}
So the inner sums over $d$ and $e$ in (5.11) become
\begin{equation*}
{\mu(r_{2}')\over \phi(r_{2}')}
\sum_{\stackrel{\scr \delta|(r_{1}',r_{2}')}{\scr \delta|j_{1}-j_{2} }}
\delta^{2}\mu(\delta)\sum_{d'|{r{_1}'\over \delta}}{d'\mu(d')\over \phi(d')}
\mu(({r_{2}'\over \delta},d'))\phi(({r_{2}'\over \delta},d')),
\end{equation*}
in which we can evaluate the sum over $d'$ as
\begin{equation*}
\prod_{p|({r_{1}'\over \delta},{r_{2}'\over \delta})}(1+p)
\prod_{\substack{p|{r_{1}'\over \delta}\\ p\ndiv {r_{2}'\over \delta}}}
{-1\over p-1} = \sigma(({r_{1}'\over \delta},{r_{2}'\over \delta}))
{\mu({{r_{1}'\over \delta}\over ({r_{1}'\over \delta},{r_{2}'\over \delta})})
\over
\phi({{r_{1}'\over \delta}\over ({r_{1}'\over \delta},{r_{2}'\over \delta})})}.
\end{equation*}
Now the inner sums over $d$ and $e$ in (5.11) have been simplified to
\begin{eqnarray*}
& \, & {\mu(r_{1}')\over \phi(r_{1}')}{\mu(r_{2}')\over \phi(r_{2}')}
\mu\cdot\phi\cdot\sigma((r_{1}',r_{2}'))
\sum_{\stackrel{\scr \delta|(r_{1}',r_{2}')}{\scr \delta|j_{1}-j_{2} }}
{\delta^{2}\mu(\delta)\over \sigma(\delta)} \\
& = & {\mu(r_{1}')\over \phi(r_{1}')}{\mu(r_{2}')\over \phi(r_{2}')}
\mu\cdot\phi\cdot\sigma((r_{1}',r_{2}'))
\prod_{p|(r_{1}',r_{2}',j_{1}-j_{2})}(1-{p^{2}\over p+1}).
\end{eqnarray*}
Plugging this into (5.11) the main term of
$\tilde{\mathcal{S}}_3(j_{1},j_{2},0)$ is now expressed as
\begin{equation}
N\sum_{r_{1},r_{2} \leq R}
{\mu^{2}(r_{1})\mu^{2}(r_{2})\mu(r_{1}')\mu(r_{2}')\over
\phi(r_{1})\phi(r_{2})\phi(r_{1}')\phi(r_{2}')}
\mu\cdot\phi\cdot\sigma((r_{1}',r_{2}'))
\prod_{p|(r_{1}',r_{2}',j_{1}-j_{2})}({1+ p -p^{2}\over p+1}).
\end{equation}
We now express $r_1 $ and $r_2 $ as products of coprime factors and
transform (5.13) into a sum over these new variables. Let
\begin{equation}
r_{1} = (r_{1},j_{1})r_{1}' = (r_{1},j_{1})(r_{1}',j_{2})r_{1}'' =
(r_{1},j_{1})(r_{1}',j_{2})(r_{1}'',j_{1}-j_{2})r_{1}''' =
s_{11}s_{12}s_{13}r_{1}''',
\end{equation}
say, so that $s_{11}|j_{1};\, s_{12}|j_{2}, (s_{12},j_{1})=1;\, s_{13}|(j_{1}-
j_{2}), (s_{13},j_{1}j_{2}) = 1; $ and \newline
$ (r_{1}''',j_{1}j_{2}(j_{1}-j_{2}))=1$.
Similarly let
\begin{equation}
r_{2} = (r_{2},j_{2})r_{2}' = (r_{2},j_{2})(r_{2}',j_{1})r_{2}'' =
(r_{2},j_{2})(r_{2}',j_{1})(r_{2}'',j_{1}-j_{2})r_{2}''' =
s_{22}s_{21}s_{23}r_{2}''',
\end{equation}
with $s_{22}|j_{2};\, s_{21}|j_{1}, (s_{21},j_{2})=1;\, s_{23}|(j_{1}-
j_{2}), (s_{23},j_{1}j_{2}) = 1;\, (r_{2}''',j_{1}j_{2}(j_{1}-j_{2}))=1$.
Then we let
\begin{equation}
(s_{13},s_{23}) = s_{3},\; s_{13}=t_{1}s_{3},\; s_{23} = t_{2}s_{3},\;
(r_{1}''',r_{2}''') = r,\; r_{1}'''=s_{14}r,\; r_{2}'''=s_{24}r,
\end{equation}
so that $t_{1}, t_{2}, s_{3}$ are each divisors
of $j_{1}-j_{2}$ which are coprime to $j_{1}j_{2}$, and $s_{14}, s_{24}, r$
are relatively prime to $j_{1}j_{2}(j_{1}-j_{2})$. We now have the
factorizations
\begin{equation}
r_{1}=s_{11}s_{12}t_{1}s_{3}s_{14}r, \qquad
r_{2}=s_{22}s_{21}t_{2}s_{3}s_{24}r,
\end{equation}
with the just stated coprimality and divisibility conditions on these
variables to be specified by a star on the summation signs below.
Now the sum of (5.11) has been transformed into
\begin{eqnarray}
\sumstar_{\stackrel{\scr s_{11}s_{12}t_{1}s_{3}s_{14}r \leq R}{\scr
s_{22}s_{21}t_{2}s_{3}s_{24}r \leq R}} & &
{\mu^{2}(s_{11})\mu(s_{12})\mu(t_{1})\mu(t_{2})\mu^{2}(s_{22})\mu(s_{21})
\over \phi(s_{11})\phi^{2}(s_{12})\phi^{2}(t_{1})\phi^{2}(t_{2})\phi(s_{22})
\phi^{2}(s_{21})}
{\mu^{2}(s_{3})\over \phi^{3}(s_{3})}\prod_{p|s_{3}}(p^2 - p -1) \nonumber \\
& & \mbox{} \times{\mu(s_{14})\mu(s_{24})\mu(r)\sigma(r)\over \phi^{2}(s_{14})
\phi^{2}(s_{24})\phi^{3}(r)}.
\end{eqnarray}
If $j_{1}=j_{2}=j\neq 0$, we have $s_{12} = s_{21} = r_{1}''' = r_{2}''' =1$,
so that (5.18) reduces to
\begin{equation}
\sumstar_{\stackrel{\scr s_{11}t_{1}s_{3} \leq R}
{\scr s_{22}t_{2}s_{3} \leq R}}
{\mu^{2}(s_{11})\mu^{2}(s_{22})\mu(t_{1})\mu(t_{2})
\over \phi(s_{11})\phi(s_{22})\phi^{2}(t_{1})\phi^{2}(t_{2})}
\mu^{2}(s_{3})\prod_{p|s_{3}}{(p^2 - p -1)\over (p-1)^3 } .
\end{equation}

We first deal with (5.19) by starting to sum over $s_{3}$ as
\begin{eqnarray}
& & \sum_{\stackrel{\scr s_{3} \leq
\min({R\over s_{11}t_{1}},{R\over s_{22}t_{2}})}
{\scr (s_{3},t_{1}t_{2}j) = 1}}\mu^{2}(s_{3})
\prod_{p|s_{3}}{(p^2 - p -1)\over (p-1)^3 } = \nonumber \\
& &  \prod_{p}(1+{p-2\over p(p-1)^2 })
\prod_{p|t_{1}t_{2}j}({(p-1)^{3}\over p^3 - 2p^2 + 2p -2})
\{ \log\min({R\over s_{11}t_{1}},{R\over s_{22}t_{2}}) + \gamma + \nonumber \\
& & \mbox{} \;\;
\sum_{p}{(2p-3)\log p\over (p-1)(p^3 -2p^2 +2p-2)} +
\sum_{p|t_{1}t_{2}j}{(p^2 -p-1) \log p\over (p^3 -2p^2 +2p-2)} \} +\nonumber \\
& & \mbox{} \;\;
O({m(t_{1}t_{2}j)\over \sqrt{\min({R\over s_{11}t_{1}},{R\over s_{22}t_{2}})}})
\end{eqnarray}
according to (2.14). The last error term will contribute to (5.19)
\begin{eqnarray}
& \ll &{m(j)\over \sqrt{R}}\sum_{\stackrel{\scr s_{11}|j}{\scr s_{22}|j}}
{\mu^{2}(s_{11})\mu^{2}(s_{22})\over \phi(s_{11})\phi(s_{22})}
\sum_{\stackrel{\scr t_{1}\leq {R\over s_{11}}}{\scr (t_{1},j)=1}}
{\mu^{2}(t_{1})m(t_{1})\over \phi^{2}(t_{1})} \nonumber \\
& & \qquad \qquad
\sum_{\stackrel{\scr t_{2}\leq {R\over s_{22}}}{\scr (t_{2},t_{1}j)=1}}
{\mu^{2}(t_{2})m(t_{2})\over \phi^{2}(t_{2})}
\sqrt{\max (s_{11}t_{1},s_{22}t_{2})}
\nonumber \\
& = & {m(j)\over \sqrt{R}}\sum_{\stackrel{\scr s_{11}|j}{\scr s_{22}|j}}
{\mu^{2}(s_{11})\mu^{2}(s_{22})\over \phi(s_{11})\phi(s_{22})}
\sum_{\stackrel{\scr t_{1}\leq {R\over s_{11}}}{\scr (t_{1},j)=1}}
{\mu^{2}(t_{1})m(t_{1})\over \phi^{2}(t_{1})} \nonumber \\
& & \qquad
\{ \sum_{\stackrel{\scr t_{2}\leq {s_{11}t_{1}\over s_{22}}}
{\scr (t_{2},t_{1}j)=1}}
{\mu^{2}(t_{2})m(t_{2})\over \phi^{2}(t_{2})} \sqrt{s_{11}t_{1}} +
\sum_{\stackrel{\scr {s_{11}t_{1}\over s_{22}} < t_{2} \leq {R\over s_{22}}}
{\scr (t_{2},t_{1}j)=1}}
{\mu^{2}(t_{2})m(t_{2})\over \phi^{2}(t_{2})}
\sqrt{s_{22}t_{2}} \; \} \nonumber \\
& \ll & {m(j)\over \sqrt{R}}\{ \sum_{\stackrel{\scr s_{11}|j}{\scr s_{22}|j}}
{\mu^{2}(s_{11})\sqrt{s_{11}}\over \phi(s_{11})}
{\mu^{2}(s_{22})\over \phi(s_{22})}
\sum_{\stackrel{\scr t_{1}\leq {R\over s_{11}}}{\scr (t_{1},j)=1}}
{\mu^{2}(t_{1})\sqrt{t_{1}}m(t_{1})\over \phi^{2}(t_{1})} \nonumber \\
& & \qquad \qquad +\sum_{\stackrel{\scr s_{11}|j}{\scr s_{22}|j}}
{\mu^{2}(s_{11})\over \phi(s_{11})}
{\mu^{2}(s_{22})\sqrt{s_{22}}\over \phi(s_{22})}
\sum_{\stackrel{\scr t_{1}\leq {R\over s_{11}}}{\scr (t_{1},j)=1}}
{\mu^{2}(t_{1})m(t_{1})\over \phi^{2}(t_{1})} \} \nonumber \\
& & \qquad
\ll {m(j)\over \sqrt{R}} \sum_{\stackrel{\scr s_{11}|j}{\scr s_{22}|j}}
{\mu^{2}(s_{11})\sqrt{s_{11}}\over \phi(s_{11})}
{\mu^{2}(s_{22})\over \phi(s_{22})}
\ll {j^{*}m^{2}(j^{*})\over \phi(j^{*})\sqrt{R}} = O(R^{-{1\over 2} +
\epsilon}).
\end{eqnarray}
Before calculating the main contribution from the $\log\min$ term,
we find the contributions from the other terms of (5.20). Since
\begin{equation}
\gamma + \sum_{p}{(2p-3)\log p\over (p-1)(p^3 -2p^2 +2p-2)} +
\sum_{p|t_{1}t_{2}j}{(p^2 -p-1) \log p\over (p^3 -2p^2 +2p-2)} \ll
\log\log 3t_{1}t_{2}j^{*}
\end{equation}
as in (2.1), the sums over $t_{1},\, t_{2}$ are each $O(1)$, and the
sums over $s_{11},\, s_{22}$ running through the divisors of $j$ are each of
value $j/\phi(j)$, and by (2.3)
\begin{equation}
({j\over \phi(j)})^{2}
\prod_{p|j}{(p-1)^{3}\over p^3 - 2p^2 + 2p -2} \ll \prod_{p|j}(1+{1\over p})
\ll \log\log 3j^{*},
\end{equation}
the secondary terms of (5.20) contribute to (5.19)
\begin{equation}
\ll (\log\log 3j^{*})^{2}.
\end{equation}
As for the main term of the brackets of (5.20), $\log R -
\log\max(s_{11}t_{1},
s_{22}t_{2})$, the $\log\max$ part will contribute little. For the $t_{1},\,
t_{2}$ sums are again $O(1)$, and the
$s_{11},\, s_{22}$ sums are
\begin{equation}
\ll \sum_{\stackrel{\scr s_{11}|j}{\scr s_{22}|j}}
{\mu^{2}(s_{11})\log s_{11}\over \phi(s_{11})}
{\mu^{2}(s_{22})\over \phi(s_{22})} \ll {j\over \phi(j)}\sum_{s_{11}|j}
{\mu^{2}(s_{11})\log s_{11}\over \phi(s_{11})} 
\qquad = ({j\over \phi(j)})^{2} \sum_{p|j}{\log p\over p} ,
\end{equation}
where we have employed (2.4).
Using (5.23), and (2.1) to bound the very last sum,
the contribution from the $\log\max$ term is also majorized as in (5.24).
We note that it is possible to carry out the calculation (5.22)-(5.24)
more precisely, but this wouldn't be significant since we have not been
able to give a better evaluation of the contribution of the $\log\max$ term.

The main term of (5.19) has now been reduced to
\begin{eqnarray}
& & \log R \prod_{p}(1+{p-2\over p(p-1)^2 })
\prod_{p|j}{(p-1)^{3}\over p^3 - 2p^2 + 2p -2}
\sum_{\stackrel{\scr s_{11}|j}{\scr s_{22}|j}}
{\mu^{2}(s_{11})\mu^{2}(s_{22})\over \phi(s_{11})\phi(s_{22})} \nonumber \\
& & \qquad
\sum_{\stackrel{\scr t_{1}\leq {R\over s_{11}}}{\scr (t_{1},j)=1}}
\mu(t_{1})\prod_{p|t_{1}}{p-1\over p^3 -2p^2 + 2p -2} \,
\sum_{\stackrel{\scr t_{2}\leq {R\over s_{22}}}{\scr (t_{2},t_{1}j)=1}}
\mu(t_{2})\prod_{p|t_{2}}{p-1\over p^3 -2p^2 + 2p -2}.
\end{eqnarray}
Now observe that by (2.9)
\begin{eqnarray}
\!\!\!\! & & \sum_{\stackrel{\scr n\leq x}{\scr (n,k)=1}}
\mu(n)\prod_{p|n}{p-1\over p^3 -2p^2 + 2p -2}
 =  \sum_{\stackrel{\scr n = 1}{\scr (n,k)=1}}^{\infty}
\mu(n)\prod_{p|n}({p-1\over p^3 -2p^2 + 2p -2}) + O({1\over x})
\nonumber \\
& & \qquad \qquad
\qquad = \prod_{p\ndiv k}(1-{p-1\over p^3 -2p^2 + 2p -2}) + O({1\over x})
\nonumber \\
& &  \qquad \qquad \qquad = \prod_{p}(1-{p-1\over p^3 -2p^2 +2p -2})
\prod_{p|k}({p^3 -2p^2 +2p -2\over p^3 -2p^2 +p-1}) + O({1\over x}).
\end{eqnarray}
Upon doing the $t_{2}$-sum according to (5.27),
the error of the last kind contributes to (5.26)
\begin{eqnarray}
& \ll & {\log R\over R}
\sum_{\stackrel{\scr s_{11}|j}{\scr s_{22}|j}}
{\mu^{2}(s_{11})\mu^{2}(s_{22})s_{22}\over \phi(s_{11})\phi(s_{22})}
\sum_{\stackrel{\scr t_{1}\leq {R\over s_{11}}}{\scr (t_{1},j)=1}}
\mu^{2}(t_{1})\prod_{p|t_{1}}{p-1\over p^3 -2p^2 + 2p -2} \nonumber \\
& \ll & {\log R\over R}
\sum_{s_{11}|j}{\mu^{2}(s_{11})\over \phi(s_{11})}
\sum_{s_{22}|j}{\mu^{2}(s_{22})s_{22}\over \phi(s_{22})}
\ll {\log R\over R}{j^{*}\over \phi(j^{*})}{j^{*}d(j^{*})\over \phi(j^{*})}
\ll R^{-1+\epsilon}.
\end{eqnarray}
From (5.27) we see that the main term is now
\begin{eqnarray}
& & \log R \prod_{p}(1-{1\over p(p-1)^2 })
\prod_{p|j}{(p-1)^{3}\over p^3 - 2p^2 + p -1}
\sum_{\stackrel{\scr s_{11}|j}{\scr s_{22}|j}}
{\mu^{2}(s_{11})\mu^{2}(s_{22})\over \phi(s_{11})\phi(s_{22})} \nonumber \\
& & \qquad
\sum_{\stackrel{\scr t_{1}\leq {R\over s_{11}}}{\scr (t_{1},j)=1}}
\mu(t_{1})\prod_{p|t_{1}}{p-1\over p^3 -2p^2 + p -1}.
\end{eqnarray}
As in (5.27) we have
\begin{equation}
\sum_{\stackrel{\scr t_{1}\leq {R\over s_{11}}}{\scr (t_{1},j)=1}}
\mu(t_{1})\prod_{p|t_{1}}{p-1\over p^3 -2p^2 + p -1} =
\prod_{p\ndiv j}(1-{p-1\over p^3 -2p^2 +p -1}) + O({s_{11}\over R}),
\end{equation}
and this last error leads to a contribution of $O(R^{-1+\epsilon})$ just as
in (5.28). The main term becomes
\begin{equation}
\log R \prod_{p}(1-{1\over p(p-1)^2 })
\prod_{p|j}{(p-1)^{3}\over p^3 - 2p^2 + p -1}
\prod_{p\ndiv j}(1-{p-1\over p^3 -2p^2 +p -1})
\sum_{\stackrel{\scr s_{11}|j}{\scr s_{22}|j}}
{\mu^{2}(s_{11})\mu^{2}(s_{22})\over \phi(s_{11})\phi(s_{22})}
\end{equation}
Note that if $2\ndiv j$, then the third product is $0$. We simplify (5.31)
and obtain the final expression for the main term as
\begin{equation}
[2|j]\, 2 \prod_{p>2}(1-{1\over (p-1)^2 })\prod_{\stackrel{\scr p|j}{
\scr p>2}}({p-1\over p-2}) \log R  = \gs_{2}(j)\log R.
\end{equation}
Since the largest error term all along this calculation (besides that of
(5.11)) was $O(N(\log\log 3j^{*})^{2})$, this completes the proof of
Theorem 2 in the case $k=3,\, r=2$ in view of (1.9).

\medskip

We now calculate the expression (5.18), with $j_{1} \neq j_{2},\, j_{1}j_{2}
\neq 0$, for the case $k=r=3$ of Theorem 2. Some additional notation will
prove to be convenient. Let $J = [j_{1}j_{2}(j_{1}-j_{2})]^{*}$,
$j_{1}^{*} = (j_{1}^{*})'(j_{1}^{*},j_{2}^{*})$,
$j_{2}^{*} = (j_{2}^{*})'(j_{1}^{*},j_{2}^{*}),\, (j_{1} - j_{2})^{*} =
(j_{1}^{*},j_{2}^{*})j_{3}^{*}$, so that
$J=(j_{1}^{*})'(j_{2}^{*})'(j_{1}^{*},j_{2}^{*})j_{3}^{*}$ as a product
of relatively prime factors. The conditions on the variables are rewritten as
$s_{11}|j_{1}^{*},\, s_{22}|j_{2}^{*},\, s_{12}|(j_{2}^{*})',\,
s_{21}|(j_{1}^{*})',\, s_{3}t_{1}t_{2}|j_{3}^{*},\, (rs_{14}s_{24},J)=1$, with
$s_{3},\, t_{1},\, t_{2}$ being pairwise coprime and the same for
$r,\, s_{14},\, s_{24}$.

We start by summing over $r$, observing that
\begin{equation}
\sum_{\stackrel{\scr n\leq x}{\scr (n,k)=1}}{\mu(r)\sigma(r)\over \phi^{3}(r)}
= \sum_{\stackrel{\scr n=1}{\scr (n,k)=1}}^{\infty}
{\mu(r)\sigma(r)\over \phi^{3}(r)} + O({1\over x})
= \prod_{p\ndiv k}(1-{p+1\over (p-1)^3 }) + O({1\over x}) .
\end{equation}
In our case
\begin{eqnarray}
& & \sum_{\stackrel{\scr r\leq \min({R\over s_{11}s_{12}t_{1}s_{3}s_{14}},
{R\over s_{21}s_{22}t_{2}s_{3}s_{24}})}{\scr (r, Js_{14}s_{24})=1}}
{\mu(r)\sigma(r)\over \phi^{3}(r)}  = \nonumber \\
& & \prod_{p}{p^3 -3p^2 + 2p -2\over (p-1)^3 }\prod_{p|Js_{14}s_{24}}
{(p-1)^{3}\over p^3 -3p^2 + 2p -2} \nonumber \\
& & \qquad
+O({1\over R}\max(s_{11}s_{12}t_{1}s_{3}s_{14}, s_{21}s_{22}t_{2}s_{3}s_{24}))
\end{eqnarray}
is fed into (5.18). The error term of (5.34) brings
\begin{eqnarray}
& & \ll {1\over R}
\sumstar_{\stackrel{\scr s_{11}s_{12}t_{1}s_{3}s_{14} \leq R}{\scr
s_{22}s_{21}t_{2}s_{3}s_{24} \leq R}}
{\mu^{2}(s_{11})\mu^{2}(s_{12})\mu^{2}(t_{1})\mu^{2}(t_{2})\mu^{2}(s_{22})
\mu^{2}(s_{21})
\over \phi(s_{11})\phi^{2}(s_{12})\phi^{2}(t_{1})\phi^{2}(t_{2})\phi(s_{22})
\phi^{2}(s_{21})} \nonumber \\
& & \mbox{} \qquad \times
{\mu^{2}(s_{3})s_{3}\over \phi^{3}(s_{3})}\prod_{p|s_{3}}(p^2 - p -1) \,
{\mu^{2}(s_{14})\mu^{2}(s_{24})\over \phi^{2}(s_{14})
\phi^{2}(s_{24})}\max(s_{11}s_{12}t_{1}s_{14}, s_{21}s_{22}t_{2}s_{24}),
\end{eqnarray}
in which the next summation over $s_{24}$ reads
\begin{equation}
\sum_{\stackrel{\scr s_{24} \leq {s_{11}s_{12}t_{1}s_{14}\over
s_{21}s_{22}t_{2}}}{\scr (s_{24}, Js_{14})=1}}
{\mu^{2}(s_{24})\over \phi^{2}(s_{24})}s_{11}s_{12}t_{1}s_{14} +
\sum_{\stackrel{\scr {s_{11}s_{12}t_{1}s_{14}\over s_{21}s_{22}t_{2}} <
s_{24} \leq {R\over s_{21}s_{22}t_{2}s_{3}}}
{\scr (s_{24}, Js_{14})=1}}
{\mu^{2}(s_{24})\over \phi^{2}(s_{24})}s_{21}s_{22}t_{2}s_{24} .
\end{equation}
This makes (5.35) majorized as
\begin{eqnarray}
& &\ll {1\over R}\, \{
\sumstar_{\stackrel{\scr s_{11}s_{12}t_{1}s_{3}s_{14} \leq R}{\scr
s_{22}s_{21}t_{2}s_{3} \leq R}}
{\mu^{2}(s_{11})s_{11}\over \phi(s_{11})}
{\mu^{2}(s_{12})s_{12}\over \phi^{2}(s_{12})}
{\mu^{2}(t_{1})t_{1}\over \phi^{2}(t_{1})}
{\mu^{2}(t_{2})\over \phi^{2}(t_{2})}
{\mu^{2}(s_{22})\over \phi(s_{22})}
{\mu^{2}(s_{21})\over \phi^{2}(s_{21})} \nonumber \\
& & \mbox{} \qquad \times
{\mu^{2}(s_{3})s_{3}\over \phi^{3}(s_{3})}\prod_{p|s_{3}}(p^2 - p -1) \,
{\mu^{2}(s_{14})s_{14}\over \phi^{2}(s_{14})} \nonumber \\
& & + \sumstar_{\stackrel{\scr s_{11}s_{12}t_{1}s_{3}s_{14} \leq R}{\scr
s_{22}s_{21}t_{2}s_{3} \leq R}}
{\mu^{2}(s_{11})\over \phi(s_{11})}
{\mu^{2}(s_{12})\over \phi^{2}(s_{12})}
{\mu^{2}(t_{1})\over \phi^{2}(t_{1})}
{\mu^{2}(t_{2})t_{2}\over \phi^{2}(t_{2})}
{\mu^{2}(s_{22})s_{22}\over \phi(s_{22})}
{\mu^{2}(s_{21})s_{21}\over \phi^{2}(s_{21})} \nonumber \\
& & \mbox{} \qquad \times
{\mu^{2}(s_{3})s_{3}\over \phi^{3}(s_{3})}\prod_{p|s_{3}}(p^2 - p -1) \,
{\mu^{2}(s_{14})\over \phi^{2}(s_{14})}\log R \, \} ,
\end{eqnarray}
where we have made use of (2.8).
Now we sum over $s_{14}$ and majorize (5.37) as
\begin{eqnarray}
& &\ll {\log R\over R}\, \{
\sumstar_{\stackrel{\scr s_{11}s_{12}t_{1}s_{3} \leq R}{\scr
s_{22}s_{21}t_{2}s_{3} \leq R}}
{\mu^{2}(s_{11})s_{11}\over \phi(s_{11})}
{\mu^{2}(s_{12})s_{12}\over \phi^{2}(s_{12})}
{\mu^{2}(t_{1})t_{1}\over \phi^{2}(t_{1})}
{\mu^{2}(t_{2})\over \phi^{2}(t_{2})}
{\mu^{2}(s_{22})\over \phi(s_{22})}
{\mu^{2}(s_{21})\over \phi^{2}(s_{21})} \nonumber \\
& & \mbox{} \qquad \times
{\mu^{2}(s_{3})s_{3}\over \phi^{3}(s_{3})}\prod_{p|s_{3}}(p^2 - p -1)
\nonumber \\
& & + \sumstar_{\stackrel{\scr s_{11}s_{12}t_{1}s_{3} \leq R}{\scr
s_{22}s_{21}t_{2}s_{3} \leq R}}
{\mu^{2}(s_{11})\over \phi(s_{11})}
{\mu^{2}(s_{12})\over \phi^{2}(s_{12})}
{\mu^{2}(t_{1})\over \phi^{2}(t_{1})}
{\mu^{2}(t_{2})t_{2}\over \phi^{2}(t_{2})}
{\mu^{2}(s_{22})s_{22}\over \phi(s_{22})}
{\mu^{2}(s_{21})s_{21}\over \phi^{2}(s_{21})} \nonumber \\
& & \mbox{} \qquad \times
{\mu^{2}(s_{3})s_{3}\over \phi^{3}(s_{3})}\prod_{p|s_{3}}(p^2 - p -1) \, \} \\
& & \ll R^{-1+\epsilon} ,
\nonumber
\end{eqnarray}
where to see the last line it is enough to observe that all of the
summations are over variables which divide $j_{1}$ or $j_{2}$ or
$j_{1}-j_{2}$, and a more precise (but more complicated looking)
factor than $R^{\epsilon}$ could easily be given.

We revert to (5.18) with the sum over $r$ already performed in (5.34) so
that the main term has been turned into
\begin{eqnarray}
& & \prod_{p}{p^3 -3p^2 +2p -2\over (p-1)^3 }
\prod_{p|J}{(p-1)^{3}\over p^3 -3p^2 +2p -2} \nonumber \\
& & \mbox{}
\times \sumstar_{\stackrel{\scr s_{11}s_{12}t_{1}s_{3}s_{14} \leq R}{\scr
s_{22}s_{21}t_{2}s_{3}s_{24} \leq R}}
{\mu^{2}(s_{11})\mu(s_{12})\mu^{2}(s_{22})\mu(s_{21})\mu(t_{1})\mu(t_{2})
\over \phi(s_{11})\phi^{2}(s_{12})\phi(s_{22})\phi^{2}(s_{21})
\phi^{2}(t_{1})\phi^{2}(t_{2})} \nonumber \\
& & \mbox{} \;\;\; \times
{\mu^{2}(s_{3})\over \phi^{3}(s_{3})}\prod_{p|s_{3}}(p^2 - p -1) \,
\mu(s_{14})\mu(s_{24})\prod_{p|s_{14}s_{24}}{p-1\over p^3 -3p^2 +2p-2}.
\end{eqnarray}
For the sum over $s_{24}$ we have
\begin{eqnarray}
& & \sum_{\stackrel{\scr s_{24}\leq {R\over s_{21}s_{22}t_{2}s_{3}}}
{\scr (s_{24},Js_{14})=1}}
\mu(s_{24})\prod_{p|s_{24}}{p-1\over p^3 -3p^2 + 2p -2} = \nonumber \\
& & \mbox{} \qquad
\prod_{p\ndiv Js_{14}}
(1-{p-1\over p^3 -3p^2 +2p -2}) + O({s_{21}s_{22}t_{2}s_{3}\over R}) .
\end{eqnarray}
The error terms, here and in what follows, can be considered similar to above
ending up with $O(R^{-1 + \epsilon})$ as in (5.38). So from now on we shall
just concentrate on the main term, which upon (5.39) has become
\begin{eqnarray}
& & \prod_{p}{p^3 -3p^2 +p -1\over (p-1)^3 }
\prod_{p|J}{(p-1)^{3}\over p^3 -3p^2 +p -1} \nonumber \\
& & \mbox{}
\times \sumstar_{\stackrel{\scr s_{11}s_{12}t_{1}s_{3}s_{14} \leq R}{\scr
s_{22}s_{21}t_{2}s_{3} \leq R}}
{\mu^{2}(s_{11})\mu(s_{12})\mu^{2}(s_{22})\mu(s_{21})\mu(t_{1})\mu(t_{2})
\over \phi(s_{11})\phi^{2}(s_{12})\phi(s_{22})\phi^{2}(s_{21})
\phi^{2}(t_{1})\phi^{2}(t_{2})} \nonumber \\
& & \mbox{} \;\;\; \times
{\mu^{2}(s_{3})\over \phi^{3}(s_{3})}\prod_{p|s_{3}}(p^2 - p -1) \,
\mu(s_{14})\prod_{p|s_{14}}{p-1\over p^3 -3p^2 +p-1}.
\end{eqnarray}
Next in row is the sum over $s_{14}$,
\begin{eqnarray}
& & \sum_{\stackrel{\scr s_{14}\leq {R\over s_{11}s_{12}t_{1}s_{3}}}
{\scr (s_{14},J)=1}}
\mu(s_{14})\prod_{p|s_{14}}{p-1\over p^3 -3p^2 + p -1} = \nonumber \\
& & \mbox{} \qquad
\prod_{p\ndiv J}
({p^{2}(p-3)\over p^3 -3p^2 +p -1}) + O({s_{11}s_{12}t_{1}s_{3}\over R}) .
\end{eqnarray}
This shows that if $3\ndiv J$, then the main term is $0$. Note that $2|J$
always. So now we can express the main term as
\begin{eqnarray}
& & [3|J]\prod_{p>3}{p^2 (p -3)\over (p-1)^3 }
\prod_{\stackrel{\scr p|J}{\scr p>3}}
{(p-1)^{3}\over p^2 (p -3)} \nonumber \\
& &
\times \sumstar_{\stackrel{\scr s_{11}s_{12}t_{1}s_{3} \leq R}{\scr
s_{22}s_{21}t_{2}s_{3} \leq R}}
{\mu^{2}(s_{11})\mu(s_{12})\mu^{2}(s_{22})\mu(s_{21})\mu(t_{1})\mu(t_{2})
\over \phi(s_{11})\phi^{2}(s_{12})\phi(s_{22})\phi^{2}(s_{21})
\phi^{2}(t_{1})\phi^{2}(t_{2})}
{\mu^{2}(s_{3})\over \phi^{3}(s_{3})}\prod_{p|s_{3}}(p^2 - p -1) . \nonumber \\
& & \mbox{}
\end{eqnarray}
From now on the inequality conditions in $\sumstar_{\mbox{}}$ \,
become superfluous
as all of the remaining variables to be summed over are divisors of $J$,
which is $\ll R^{\epsilon}$, and therefore satisfy these inequalities anyway.
The sum over $s_{3}$ is
\begin{equation}
\sum_{s_{3}|{j_{3}^{*}\over t_{1}t_{2}}}\mu^{2}(s_{3})\prod_{p|s_{3}}
{p^2 -p-1\over (p-1)^3 } = \prod_{p|{j_{3}^{*}\over t_{1}t_{2}}}
{p^3 -2p^2 + 2p-2\over (p-1)^3 },
\end{equation}
turning the main term into
\begin{eqnarray}
& & [3|J]\prod_{p>3}{p^2 (p -3)\over (p-1)^3 }
\prod_{\stackrel{\scr p|J}{\scr p>3}}{(p-1)^{3}\over p^2 (p -3)}
\prod_{p|j_{3}^{*}}{p^3 -2p^2 +2p-2\over (p-1)^3 } \nonumber \\
& & \qquad
\times \sumstar_{\mbox{}}
{\mu^{2}(s_{11})\mu(s_{12})\mu^{2}(s_{22})\mu(s_{21})
\over \phi(s_{11})\phi^{2}(s_{12})\phi(s_{22})\phi^{2}(s_{21})}
\mu(t_{1})\mu(t_{2})\prod_{p|t_{1}t_{2}}{p-1\over p^3 -2p^2 +2p-2}.
\end{eqnarray}
Now the sum over $t_{2}$ is
\begin{equation}
\sum_{t_{2}|{j_{3}^{*}\over t_{1}}}\mu(t_{2})\prod_{p|t_{2}}
{p-1\over p^3 -2p^2 +2p-2 } = \prod_{p|{j_{3}^{*}\over t_{1}}}
{p^3 -2p^2 + p-1\over p^3 -2p^2 + 2p-2},
\end{equation}
and the main term becomes
\begin{eqnarray}
& & [3|J]\prod_{p>3}{p^2 (p -3)\over (p-1)^3 }
\prod_{\stackrel{\scr p|J}{\scr p>3}}{(p-1)^{3}\over p^2 (p -3)}
\prod_{p|j_{3}^{*}}{p^3 -2p^2 +p-1\over (p-1)^3 } \nonumber \\
& & \qquad
\times \sumstar_{\mbox{}}
{\mu^{2}(s_{11})\mu(s_{12})\mu^{2}(s_{22})\mu(s_{21})
\over \phi(s_{11})\phi^{2}(s_{12})\phi(s_{22})\phi^{2}(s_{21})}
\mu(t_{1})\prod_{p|t_{1}}{p-1\over p^3 -2p^2 +p-1}.
\end{eqnarray}
We continue by summing over $t_{1}$,
\begin{equation}
\sum_{t_{1}|j_{3}^{*}}\mu(t_{1})\prod_{p|t_{1}}
{p-1\over p^3 -2p^2 +p-1 } = \prod_{p|j_{3}^{*}}
{p^2 (p-2)\over p^3 -2p^2 + p-1},
\end{equation}
where the very last product is $0$ if $2|j_{3}^{*}$ (which is equivalent to
$2\ndiv j_{1}j_{2}$). Now the main term has been simplified to
\begin{eqnarray}
& & [3|J]\prod_{p>3}{p^2 (p -3)\over (p-1)^3 }
\prod_{\stackrel{\scr p|J}{\scr p>3}}{(p-1)^{3}\over p^2 (p -3)}
\prod_{p|j_{3}^{*}}{p^2 (p-2)\over (p-1)^3 } \nonumber \\
& & \qquad
\times \sumstar_{\mbox{}}
{\mu^{2}(s_{11})\mu(s_{12})\mu^{2}(s_{22})\mu(s_{21})
\over \phi(s_{11})\phi^{2}(s_{12})\phi(s_{22})\phi^{2}(s_{21})}.
\end{eqnarray}
For the final reduction of the main term we have
\begin{eqnarray}
\sum_{s_{11}|j_{1}^{*}}{\mu^{2}(s_{11})\over \phi(s_{11})} & = &
{j_{1}^{*}\over \phi(j_{1}^{*})} = {(j_{1}^{*})'(j_{1}^{*},j_{2}^{*})\over
\phi((j_{1}^{*})')\phi((j_{1}^{*},j_{2}^{*}))} \nonumber \\
\sum_{s_{22}|j_{2}^{*}}{\mu^{2}(s_{22})\over \phi(s_{22})} & = &
{(j_{2}^{*})'(j_{1}^{*},j_{2}^{*})\over
\phi((j_{2}^{*})')\phi((j_{1}^{*},j_{2}^{*}))} \nonumber \\
\sum_{s_{12}|(j_{2}^{*})'}{\mu(s_{12})\over \phi^{2}(s_{12})} & = &
\prod_{p|(j_{2}^{*})'}{p(p-2)\over (p-1)^2 } \nonumber \\
\sum_{s_{21}|(j_{1}^{*})'}{\mu(s_{21})\over \phi^{2}(s_{21})} & = &
\prod_{p|(j_{1}^{*})'}{p(p-2)\over (p-1)^2 } ,
\end{eqnarray}
which shows that the main term is $0$ if $2|(j_{1}^{*})'(j_{2}^{*})'$. Hence,
since $2|J$, in order to have a nonzero main term it must be that
$2|(j_{1}^{*},j_{2}^{*})$. Thus we find that the result of the summation (5.18)
is
\begin{eqnarray}
& & [2|(j_{1}^{*},j_{2}^{*})][3|J]\prod_{p>3}{p^2 (p -3)\over (p-1)^3 }
\prod_{\stackrel{\scr p|J}{\scr p>3}}{(p-1)^{3}\over p^2 (p -3)}
\prod_{p|j_{3}^{*}}{p^2 (p-2)\over (p-1)^3 } \nonumber \\
& & \qquad \times
{(j_{1}^{*})'\over \phi((j_{1}^{*})')}{(j_{2}^{*})'\over \phi((j_{2}^{*})')}
({(j_{1}^{*},j_{2}^{*})\over \phi((j_{1}^{*},j_{2}^{*}))})^{2}
\prod_{p|j_{1}^{*}j_{2}^{*}}{p(p-2)\over (p-1)^2 } .
\end{eqnarray}
To re-organize (5.51), note that the product over all $p>3$ is
$ {4\over 3}C_{2}C_{3}$ (see (1.13)), and recall that for
nonzero main term $2|(j_{1}^{*},j_{2}^{*})$ so that
$2\ndiv (j_{1}^{*})'(j_{2}^{*})'j_{3}^{*}$, and then consider one by one
the four possibilities arising from
$3|(j_{1}^{*})'(j_{2}^{*})'(j_{1}^{*},j_{2}^{*})j_{3}^{*}$
as to which factor $3$ divides. In this way we find that the main term is
\begin{eqnarray}
& &
[2|(j_{1},j_{2})][3|j_{1}j_{2}(j_{1} - j_{2})]\, 6 C_{2} C_{3}\prod_{\stackrel{
\scr p|(j_{1},j_{2})}{\scr p>2}}({p-1\over p-2})\prod_{\stackrel{\scr
p|j_{1}j_{2}(j_{1}- j_{2})}{\scr p>3}}({p-2\over p-3}) \nonumber \\
& & \qquad = \gs_{2}((j_{1},j_{2}))\gs_{3}(j_{1}j_{2}(j_{1} - j_{2})),
\end{eqnarray}
and, by (1.10), this completes the proof of Theorem 2 for the case $k=r=3$.

\section{Re-expression of the pure triple correlations
$\mathcal{S}_{3}(N,\mbox{\boldmath$j$},\mbox{\boldmath$a$}) $}

In this section the sum
\begin{equation}
\mathcal{S}_{3}(N,\mbox{\boldmath$j$},\mbox{\boldmath$a$}) =
\sum_{n=1}^N \lambda_R(n+j_{1})\lambda_{R}(n+j_{2})\lambda_R(n)
\end{equation}
will be reduced to a multiple sum over relatively prime variables.
Later on three cases will be considered: If $j_{1} = j_{2} = 0$, then
the sum is $\mathcal{ S}_{3}(N,(0),(3))$; if $j_{1} = j_{2} = j \neq 0$, then
the sum is $\mathcal{ S}_{3}(N,(0,j),(1,2))$; if $j_{1}\neq j_{2}$ and
$j_{1}j_{2}\neq 0$, then the sum is
$\mathcal{ S}_{3}(N,(0,j_{1},j_{2}),(1,1,1))$.
From the definition of $\lambda_{R}(n)$'s we have
\begin{equation}
\mathcal{S}_{3}(N,\mbox{\boldmath$j$},\mbox{\boldmath$a$}) =
\sum_{r_{1}, r_{2}, r_{3} \leq R}{\mu^{2}(r_1)\mu^{2}(r_2)\mu^{2}(r_{3})\over
\phi(r_1)\phi(r_2)\phi(r_3)}\sum_{\stackrel{\stackrel{\scr d|r_{1}}{\scr
e|r_2 }}{\scr f|r_{3}}}d\mu(d)e\mu(e)f\mu(f)
\sum_{\stackrel{\stackrel{\stackrel{\scr n\leq N}{\scr n\equiv -j_{1}
(\mathrm{mod}\, d)}}{\scr n\equiv -j_{2} (\mathrm{mod}\, e)}}{\scr n\equiv 0
(\mathrm{mod}\, f)}} 1.
\end{equation}
The innermost sum is over $n$'s in a unique residue class modulo
$[d,e,f]$ whenever $(d,e)|j_{1}-j_{2}, \, (d,f)|j_{1},\, (e,f)|j_{2}$
in which case its value is $N/[d,e,f] + O(1)$,
otherwise the innermost sum is void. As in (4.5),
by (4.2) the last $O(1)$ leads to a
contribution of $O(R^3)$ in (6.2). Hence
\begin{equation}
\mathcal{S}_{3}(N,\mbox{\boldmath$j$},\mbox{\boldmath$a$}) =
N\sum_{r_{1}, r_{2}, r_{3} \leq R}{\mu^{2}(r_1)\mu^{2}(r_2)\mu^{2}(r_{3})\over
\phi(r_1)\phi(r_2)\phi(r_3)}
\sum_{\stackrel{\stackrel{\scr d|r_{1}, e|r_{2}, f|r_{3}}{\scr
(d,e)|j_{1}-j_{2}}}{\scr (d,f)|j_{1}, (e,f)|j_{2}}}
{d\mu(d)e\mu(e)f\mu(f)\over [d,e,f]} + O(R^3).
\end{equation}
We can express the summation variables as products of coprime factors (since
the M\"{o}bius function restricts us to squarefree $r_{i}$'s)
\begin{equation}
r_{1} = a_{1}a_{12}a_{13}a_{123} ; \qquad r_{2} = a_{2}a_{12}a_{23}a_{123} ;
\qquad r_{3} = a_{3}a_{13}a_{23}a_{123},
\end{equation}
with the understanding that
for a subscript $\chi$, $a_{\chi}$ is a divisor of those $r_{j}$'s where $j$
occurs in $\chi$. We can now write
\begin{equation}
d = d_{1}d_{12}d_{13}d_{123}; \qquad e = e_{2}e_{12}e_{23}e_{123} ; \qquad
f = f_{3}f_{13}f_{23}f_{123} ,
\end{equation}
where a $d$ or $e$ or $f$ with a certain subscript is a divisor of the $a$
with the same subscript (e.g. $d_{12} \mid a_{12}$). Then we have
\begin{equation}
[d,e,f] = d_{1}e_{2}f_{3}[d_{12},e_{12}][d_{13},f_{13}][e_{23},f_{23}]
[d_{123},e_{123},f_{123}] .
\end{equation}
So the inner sum of (6.3) over $d, e, f$ becomes
\begin{equation}
\sumstar_{\stackrel{\stackrel{\scr d_{1}d_{12}d_{13}d_{123}
\mid a_{1}a_{12}a_{13}a_{123}}
{\scr e_{2}e_{12}e_{23}e_{123}
\mid a_{2}a_{12}a_{23}a_{123}}}
{\scr f_{3}f_{13}f_{23}f_{123} \mid a_{3}a_{13}a_{23}a_{123}}}
{d_{12}e_{12}\over [d_{12},e_{12}]}{d_{13}f_{13}\over [d_{13},f_{13}]}
{e_{23}f_{23}\over [e_{23},f_{23}]}{d_{123}e_{123}f_{123}\over
[d_{123},e_{123},f_{123}]}\mu(d_{1})\cdots \mu(f_{123}),
\end{equation}
where  $\cdots$ indicates that we have the M\"{o}bius functions
of all of the twelve variables coming from (6.5), and the star in
$\sumstar_{\mbox{}}$ \, reminds us that the variables of summation also obey
the conditions
\begin{eqnarray}
(d,e) & = & (d_{12},e_{12})(d_{123},e_{123})|j_{1}-j_{2} \nonumber \\
(d,f) & = & (d_{13},f_{13})(d_{123},f_{123})|j_{1} \nonumber \\
(e,f) & = & (e_{23},f_{23})(e_{123},f_{123})|j_{2} .
\end{eqnarray}
Now (6.7) can be broken into
simpler sums as
\begin{eqnarray}
& & \!\!\! \sum_{\stackrel{\stackrel{\stackrel{\scr d_{123},\, e_{123},\,
f_{123}|a_{123}}
{\scr (d_{123},e_{123})|j_{1}-j_{2}}}{\scr (d_{123},f_{123})|j_{1}}}
{\scr (e_{123},f_{123})|j_{2}}}
{\mu(d_{123})\mu(e_{123})\mu(f_{123})d_{123}e_{123}
f_{123}\over [d_{123},e_{123},f_{123}]} \nonumber \\
& & \sum_{\stackrel{\scr d_{12},\, e_{12}| a_{12}}{\scr
(d_{12},e_{12})| {j_{1}-j_{2}\over (d_{123},e_{123})}}}
\mu(d_{12})\mu(e_{12})(d_{12},e_{12})
\sum_{\stackrel{\scr d_{13},\, f_{13}| a_{13}}{\scr
(d_{13},f_{13})|{j_{1}\over (d_{123},f_{123})}}}
\mu(d_{13})\mu(f_{13})(d_{13},f_{13}) \nonumber \\
& & \sum_{\stackrel{\scr e_{23},\, f_{23}| a_{23}}{\scr
(e_{23},f_{23})|{j_{2}\over (e_{123},f_{123})}}}
\mu(e_{23})\mu(f_{23})(e_{23},f_{23})
\sum_{d_{1}| a_{1}}\mu(d_{1})
\sum_{e_{2}\mid a_{2}}\mu(e_{2})\sum_{f_{3}\mid a_{3}}\mu(f_{3}) .
\end{eqnarray}
The last three sums yield a nonzero contribution only if $a_{1}=a_{2}=a_{3}=1$.
As for the other sums, by multiplicativity, it suffices to evaluate them
when
the $a_{\chi}$ are prime. We have for squarefree $a_{\chi}$'s
\begin{eqnarray}
\sum_{\stackrel{\scr d_{12},\, e_{12}| a_{12}}{\scr
(d_{12},e_{12}) \mid {j_{1}-j_{2}\over (d_{123},e_{123})}}}
\mu(d_{12})\mu(e_{12})(d_{12},e_{12}) & = &
\mu(a_{12})\mu\cdot\phi(a_{12},j_{1}-j_{2}) \nonumber \\
\sum_{\stackrel{\scr d_{13},\, f_{13}| a_{13}}{\scr
(d_{13},f_{13})|{j_{1}\over (d_{123},f_{123})}}}
\mu(d_{13})\mu(f_{13})(d_{13},f_{13}) & = &
\mu(a_{13})\mu\cdot\phi(a_{13},j_{1}) \nonumber \\
\sum_{\stackrel{\scr e_{23},\, f_{23}| a_{23}}{\scr
(e_{23},f_{23})|{j_{2}\over (e_{123},f_{123})}}}
\mu(e_{23})\mu(f_{23})(e_{23},f_{23}) & = &
\mu(a_{23})\mu\cdot\phi(a_{23},j_{2}) ,
\end{eqnarray}
and
\begin{eqnarray}
& & \qquad \qquad
\sum_{\stackrel{\stackrel{\stackrel{\scr d_{123},\, e_{123},\,
f_{123}|a_{123}}
{\scr (d_{123},e_{123})|j_{1}-j_{2}}}{\scr (d_{123},f_{123})|j_{1}}}
{\scr (e_{123},f_{123})|j_{2}}}
{\mu(d_{123})\mu(e_{123})\mu(f_{123})d_{123}e_{123}
f_{123}\over [d_{123},e_{123},f_{123}]} =  \\
&  & \mu\cdot\phi((a_{123},j_{1},j_{2})) \nonumber \\
& & \qquad \times \, \phi_{2}((a_{123},j_{1}))
\phi_{2}({(a_{123},j_{2})\over (a_{123},j_{1},j_{2})})
\phi_{2}({(a_{123},j_{1}-j_{2})\over (a_{123},j_{1}-j_{2},j_{1}j_{2})})\!\!\!\!
\prod_{\stackrel{\scr p|a_{123}}{\scr p\ndiv j_{1}j_{2}(j_{1}-j_{2})}}
\!\!\!\! (-2) . \nonumber
\end{eqnarray}

Thus the sum in the main term of (6.3) has been transformed into a sum
over pairwise coprime variables $a_{\chi}$,
\begin{eqnarray}
& &  \sumprime_{\stackrel{\stackrel{\scr a_{12}a_{13}a_{123} \leq R}
{\scr a_{12}a_{23}a_{123} \leq R}}
{\scr a_{13}a_{23}a_{123} \leq R}}
{\mu(a_{12})\mu\cdot\phi((a_{12},j_{1}-j_{2}))
\mu(a_{13})\mu\cdot\phi((a_{13},j_{1}))
\mu(a_{23})\mu\cdot\phi((a_{23},j_{2}))\over
\phi^{2}(a_{12})\phi^{2}(a_{13})\phi^{2}(a_{23})} \nonumber \\
& & \times
{\mu^{2}(a_{123})\mu\cdot\phi((a_{123},j_{1},j_{2}))
\phi_{2}((a_{123},j_{1}))
\phi_{2}({(a_{123},j_{2})\over (a_{123},j_{1},j_{2})})
\phi_{2}({(a_{123},j_{1}-j_{2})\over (a_{123},j_{1}-j_{2},j_{1}j_{2})})\over
\phi^{3}(a_{123})} \nonumber  \\
& & \qquad \qquad \times \, \mu\cdot d({a_{123}\over
(a_{123},j_{1}j_{2}(j_{1}-j_{2}))}).
\end{eqnarray}

\section{Pure triple correlations: The case $j_1 = j_2 = 0$}

In this case (6.12) reads
\begin{equation}
\sumprime_{\stackrel{\stackrel{\scr uvy \leq R}{\scr uwy \leq R}}
{\scr vwy \leq R}}{\mu^{2}(u)\mu^{2}(v)
\mu^{2}(w)\mu(y)\phi_{2}(y)\over \phi(u)\phi(v)\phi(w)\phi^{2}(y)} .
\end{equation}
We begin by summing over $u$ as
\begin{equation}
\sum_{\stackrel{\scr u \leq \min({R\over vy}, {R\over wy})}{\scr
(u,vwy) =1}}
{\mu^{2}(u)\over \phi(u)} ,
\end{equation}
which is evaluated by (2.15), and plugging into (7.1) we have
\begin{eqnarray}
& & \sumprime_{vwy \leq R}
{\mu^{2}(v)\mu^{2}(w)\mu(y)\phi_{2}(y)\over vwy\phi(y)}
[\log\min({R\over vy}, {R\over wy}) + D_{1}
+ \sum_{p \mid vwy}{\log p\over p}] \nonumber \\
& & \qquad +O(\sumprime_{vwy \leq
R}{\mu^{2}(v)\mu^{2}(w)\mu^{2}(y)\phi_{2}(y)
m(vwy)\over \phi(v)\phi(w)\phi^{2}(y)\sqrt{\min({R\over vy},
{R\over wy})}}),
\end{eqnarray}
where we have written $D_{1}$ for the constant
$ \gamma + \sum_{p}{\log p\over p(p-1)}$.
We will regard the $\log\min$ occurring in (7.3) as
$\log({R\over vwy}) + \log\min(v,w)$. First, upon letting $z=vwy$, we have
\begin{eqnarray}
 \sumprime_{vwy \leq R}
{\mu^{2}(v)\mu^{2}(w)\mu(y)\phi_{2}(y)\over vwy\phi(y)}
\log({R\over vwy}) & \!\!\! = & \!\!\! \sum_{z\leq
R}{\mu^{2}(z)\over z}\log ({R\over z}) \sum_{y\mid
z}{\mu(y)\phi_{2}(y)\over \phi(y)}\sum_{w\mid {z\over y}}1
\nonumber \\
= \!\! \sum_{z\leq R}\!{\mu^{2}(z)d(z)\over z}\log ({R\over z})
\!\sum_{y\mid z}{\mu(y)\phi_{2}(y)\over d(y)\phi(y)} & \!\!\!\! = &
\!\!\!\!\sum_{z\leq R}{\mu^{2}(z)d(z)\over z}\log ({R\over
z})\!\prod_{p\mid
z}(1-{p-2\over 2(p-1)}) \nonumber \\
 = \sum_{z\leq R}{\mu^{2}(z)\over \phi(z)}\log
({R\over z})& = &{1\over 2}\log^{2}R + O(\log R),
\end{eqnarray}
by (2.16). Similarly, we have
\begin{equation}
\sumprime_{vwy \leq R}
{\mu^{2}(v)\mu^{2}(w)\mu(y)\phi_{2}(y)\over vwy\phi(y)}
= \sum_{z\leq R}{\mu^{2}(z)\over \phi(z)} = \log R + O(1),
\end{equation}
and
\begin{eqnarray}
& & \sumprime_{vwy \leq R}
{\mu^{2}(v)\mu^{2}(w)\mu(y)\phi_{2}(y)\over vwy\phi(y)}
\sum_{p \mid vwy}{\log p\over p}
= \sum_{z\leq R}{\mu^{2}(z)\over \phi(z)}\sum_{p \mid z}{\log p\over p}
 \nonumber \\
& & = \sum_{p\leq R}{\log p\over p}
\sum_{\stackrel{\scr z\leq R}{\scr p\mid z}}
{\mu^{2}(z)\over \phi(z)}
=  \sum_{p\leq R}{\log p\over p\phi(p)}
\sum_{\stackrel{\scr m\leq {R\over p}}{\scr (m,p)=1}}
{\mu^{2}(m)\over \phi(m)} \nonumber \\
& & \ll \sum_{p\leq R}{\log {2R\over p}\log p\over p^2}
\ll  \log R .
\end{eqnarray}
We can take the $\log\min(v,w)$-term as twice the summands with $w<v$, that is
\begin{equation}
2 \sumprime_{\stackrel{\scr vwy \leq R}{\scr w<v}}
{\mu^{2}(v)\mu^{2}(w)\mu(y)\phi_{2}(y)\over vwy\phi(y)}\log w.
\end{equation}
Now observe that in the sum of (7.5) the terms with $v=w$, which can only
be present if $v=w=1$, contribute
\begin{equation}
\sum_{y\leq R}{\mu(y)\phi_{2}(y)\over y\phi(y)} \ll 1,
\end{equation}
by Lemma 2. So the sum on the left-hand side of (7.5) is
\begin{displaymath}
2 \sumprime_{\stackrel{\scr vwy \leq R}{\scr w<v}}
{\mu^{2}(v)\mu^{2}(w)\mu(y)\phi_{2}(y)\over vwy\phi(y)} + O(1).
\end{displaymath}
Writing
\begin{equation}
\sum_{\stackrel{\scr y\leq {R\over vw}}{\scr (y,vw)=1}}
{\mu(y)\phi_{2}(y)\over y\phi(y)} = f_{vw}({R\over vw}); \qquad
\sum_{\stackrel{\scr w<v\leq {R\over w}}{\scr (v,w)=1}}
{\mu^{2}(v)\over v}
f_{vw}({R\over vw}) = g_{w}(w)
\end{equation}
(the subscripts of $f$ and $g$ refer to the coprimality conditions in
the sums),
(7.5) says that
\begin{equation}
2\sum_{w \leq \sqrt{R}}{\mu^{2}(w)\over w}g_{w}(w) = \log R + O(1).
\end{equation}
With this notation (7.7) can be rewritten as
\begin{equation}
2\sum_{w\leq \sqrt{R}}{\mu^{2}(w)\over w}g_{w}(w)\log w,
\end{equation}
and by partial summation on (7.10) we find
\begin{equation}
2\sum_{w\leq \sqrt{R}}{\mu^{2}(w)\over w}g_{w}(w)\log w = {1\over 4}\log^2 R
+ O(\log R).
\end{equation}
Thus the first line of (7.3) has been shown to be
$({3\over 4}\log^{2}R + O(\log 2R))$, for $R\geq 1$.

It remains to consider the error term of (7.3), which can be rewritten as
\begin{equation}
{2\over \sqrt{R}}\sumprime_{\stackrel{\scr vwy \leq R}{\scr w<v}}
{\mu^{2}(v)\mu^{2}(w)\mu^{2}(y)m(vwy)\phi_{2}(y)
\sqrt{vy}\over \phi(v)\phi(w)\phi^{2}(y)} + {1\over \sqrt{R}}
\sum_{y\leq R}{\mu^{2}(y)m(y)\phi_{2}(y)\sqrt{y}\over \phi^{2}(y)}.
\end{equation}
Let
\begin{equation}
h_{vw}({R\over vw}) =
\sum_{\stackrel{\scr y\leq {R\over vw}}{\scr (y,vw)=1}}
{\mu^{2}(y)m(y)\phi_{2}(y)\sqrt{y}\over \phi^{2}(y)} ,
\end{equation}
so that (7.13) is
\begin{eqnarray}
& & \qquad {2\over \sqrt{R}}
\sumprime_{\stackrel{\scr vw \leq R}{\scr w<v}}
{\mu^{2}(v)\mu^{2}(w)m(vw)\sqrt{v}\over \phi(v)\phi(w)}h_{vw}({R\over
vw}) + h_{1}(R)
\nonumber  \\
& = & {2\over \sqrt{R}}\sum_{w\leq \sqrt{R}}{\mu^{2}(w)m(w)\over \phi(w)}
\sum_{\stackrel{\scr w<v \leq {R\over w}}{\scr (v,w)=1}}
{\mu^{2}(v)m(v)\sqrt{v}\over \phi(v)}h_{vw}({R\over vw}) + h_{1}(R)
\nonumber \\
& = & {2\over \sqrt{R}}\sum_{w\leq \sqrt{R}}{\mu^{2}(w)m(w)\over
\phi(w)}k_{w}(w) + h_{1}(R) ,
\end{eqnarray}
say.
Note that $h_{1}(R)$ is the contribution of the terms with
$v=w=1$ in the error term of (7.3), and we know by (2.11) that
\begin{equation}
h_{1}(R) \leq {1\over \sqrt{R}}\sum_{y\leq R}{\mu^{2}(y)m(y)\over
\sqrt{y}} \ll 1 .
\end{equation}
Now consider the expression
\begin{equation}
{1\over \sqrt{R}}
\sumprime_{vwy \leq R}{\mu^{2}(v)\mu^{2}(w)\mu^{2}(y)m(vwy)\phi_{2}(y)
\sqrt{vwy}\over \phi(v)\phi(w)\phi^{2}(y)}.
\end{equation}
On one hand with the notation defined in (7.14) and (7.15),
(7.17) may be re-expressed as
\begin{equation}
{2\over \sqrt{R}}\sum_{w\leq \sqrt{R}}{\mu^{2}(w)m(w)\sqrt{w}\over
\phi(w)}k_{w}(w) + h_{1}(R) .
\end{equation}
On the other hand by putting $z=vwy$, the sum of (7.17) becomes
\begin{eqnarray}
& & \sum_{z\leq R}{\mu^{2}(z)m(z)\sqrt{z}\over \phi(z)}\sum_{y\mid
z} {\phi_{2}(y)\over \phi(y)}\sum_{w\mid {z\over y}}1  =
\sum_{z\leq R}{\mu^{2}(z)m(z)\sqrt{z}d(z)\over \phi(z)}\sum_{y\mid
z}{\phi_{2}(y)\over d(y)\phi(y)} \nonumber \\
& & = \sum_{z\leq R}{\mu^{2}(z)m(z)\sqrt{z}d(z)\over
\phi(z)}\prod_{p\mid z}(1+ {p-2\over 2(p-1)})  =  \sum_{z\leq
R}\mu^{2}(z) \prod_{p\mid z}{3p-4\over (p-1)(\sqrt{p}-1)}.
\nonumber \\
\mbox{} & & \mbox{}
\end{eqnarray}
The last sum has been calculated in Lemma 3, giving
\begin{equation}
\sum_{w\leq \sqrt{R}}{\mu^{2}(w)m(w)\sqrt{w}\over \phi(w)}k_{w}(w) =
2P(1)R^{{1\over 2}}\log^{2}\sqrt{R} + O(R^{{1\over 2}}\log R).
\end{equation}
From here we obtain by partial summation
\begin{equation}
\sum_{w\leq \sqrt{R}}{\mu^{2}(w)m(w)\over \phi(w)}k_{w}(w) =
P(1)R^{{1\over 4}}\log^{2}R + O(R^{{1\over 4}}\log R).
\end{equation}
By (7.21) and (7.16), the quantity in (7.15), i.e. the error
term of (7.3) is $O(1)$. This completes the evaluation.

\section{Pure triple correlations: The case $j_1 = j_2 = j \neq 0$}

In this case (6.12) assumes the form
\begin{equation}
\sumprime_{\stackrel{\stackrel{\scr uvy \leq R}{\scr uwy
\leq R}}{\scr vwy \leq R}}
{\mu^{2}(u)\mu(v)\mu\cdot\phi((v,j))
\mu(w)\mu\cdot\phi((w,j))
\mu^{2}(y)\mu\cdot\phi((y,j))\phi_{2}(y)
\over \phi(u)\phi^{2}(v)\phi^{2}(w)\phi^{3}(y)} ,
\end{equation}
and doing first the sum over $u$ as in (7.2), we get
\begin{eqnarray}
& & \sumprime_{vwy \leq R}
{\mu(v)\mu\cdot\phi((v,j))\mu(w)\mu\cdot\phi((w,j))
\mu^{2}(y)\mu\cdot\phi((y,j))\phi_{2}(y)
\over v\phi(v) w\phi(w) y\phi^{2}(y)} \nonumber \\
& & \qquad \qquad \times
[\log\min({R\over vy}, {R\over wy}) + D_{1}
+ \sum_{p \mid vwy}{\log p\over p}]  \\
& & +O(\sumprime_{vwy \leq
R}{\mu^{2}(v)\mu^{2}(w)\mu^{2}(y)\phi_{2}(y)\phi((vwy,j))
m(vwy)\over \phi^{2}(v)\phi^{2}(w)\phi^{3}(y)
\sqrt{\min({R\over vy},{R\over wy})}}).
\nonumber
\end{eqnarray}
The last error term is at most as large as the error term of (7.3) which was
shown to be $O(1)$. As in (7.4), letting $z=vwy$, the $\log$-term in the
brackets being $\log({R\over z}) + \log\min(v,w)$,
we view the main part of (8.2) except for the contribution of the last
$\log\min$-term in the form
\begin{equation}
\sum_{z\leq R}{\mu(z)\mu\cdot\phi((z,j))a(z)\over z\phi(z)}
\sum_{y\mid z}{\mu(y)\phi_{2}(y)\over \phi(y)}\sum_{w\mid {z\over
y}}1   = \sum_{z\leq R}{\mu(z)\mu\cdot\phi((z,j))a(z)\over
\phi^{2}(z)}.
\end{equation}
In the present situation $a(z) = \log({R\over z}) + D_{1} +
\sum_{p|z}{\log p\over p}$. By Lemma 4 and (2.21), we have
\begin{eqnarray}
& & \qquad \qquad \sum_{z\leq R}
{\mu(z)\mu\cdot\phi((z,j))\over \phi^{2}(z)}(\log {R\over z} + D_{1}) = \\
& & \gs_{2}(j)(\log R + D_{1})
+ O({j^{*}d(j^{*})\log R\over \phi(j^{*})R}) \nonumber \\
& & + [2\ndiv j]\gs_{2}(2j){\log 2\over 2} + [2|j]
\gs_{2}(j)[\sum_{p\ndiv j}{\log p\over p(p-2)} -
\sum_{p|j}{\log p\over p}] . \nonumber
\end{eqnarray}
Next we have, by Lemma 4,
\begin{eqnarray*}
& & \sum_{z\leq R}
{\mu(z)\mu\cdot\phi((z,j))\over \phi^{2}(z)}\sum_{p|z}{\log p\over p} =
\sum_{p\leq R}{\log p\over p}\sum_{\stackrel{\scr z\leq R}{\scr p|z}}
{\mu(z)\mu\cdot\phi((z,j))\over \phi^{2}(z)}  \\
& & = \sum_{p\leq R}{\log p\over p}
{\mu(p)\mu\cdot\phi((p,j))\over \phi^{2}(p)}
\sum_{\stackrel{\scr u\leq {R\over p}}{\scr (u,p)=1}}
{\mu(u)\mu\cdot\phi((u,j))\over \phi^{2}(u)}
\end{eqnarray*}
\medskip
\begin{eqnarray}
& & = -{\mu((2,j))\log 2\over 2} \sum_{\stackrel{\scr u\leq
{R\over 2}}{\scr (u,2)=1}} {\mu(u)\mu\cdot\phi((u,j))\over
\phi^{2}(u)} \nonumber \\
& & \qquad - \sum_{2<p \leq R}{\mu\cdot\phi((p,j))\log p\over
p\phi^{2}(p)} \sum_{\stackrel{\scr u\leq {R\over p}}{\scr
(u,p)=1}}
{\mu(u)\mu\cdot\phi((u,j))\over \phi^{2}(u)} \nonumber \\
& & = -{\mu((2,j))\log 2\over 4}\gs_{2}(2j) -
[2|j]\, 2C_{2}
\sum_{2<p \leq R}{\mu\cdot\phi((p,j))\log p\over p^{2}(p-2)}
\prod_{\stackrel{\stackrel{\scr q: \mathrm{ prime}}{\scr q|j}}{\scr
q\neq 2,\, p}}({q-1\over q-2}) \nonumber \\
& & \qquad \qquad + O({j^{*}d(j^{*})\over R\phi(j^{*})}\sum_{p\leq
R}{\phi((p,j))\log p\over (p-1)^{2}}) \nonumber \\
& & = -{\mu((2,j))\log 2\over 4}\gs_{2}(2j) + \gs_{2}(j)
[\sum_{\stackrel{\scr 2<p\leq R}{\scr p|j}}{\log p\over p^{2}} -
\sum_{\stackrel{\scr 2<p\leq R}{\scr p\ndiv j}}{\log p\over
p^{2}(p-2)}]
\nonumber \\
& & \qquad + O({j^{*}d(j^{*})\over R\phi(j^{*})}
 [\sum_{\stackrel{\scr
p\leq R}{\scr p|j}}{\log p\over p-1} + \sum_{\stackrel{\scr p\leq
R}{\scr p\ndiv j}}{\log p\over (p-1)^{2}}]),
\end{eqnarray}
where, by (2.1), the sums in the very last $O$-term are $\ll \log\log  3j^{*}$.

To finish off this case we deal with the $\log\min$-term for which we have
\begin{eqnarray}
& & |\sum_{z\leq R}{\mu(z)\mu\cdot\phi((z,j))\over z\phi(z)}
\sum_{y\mid z}{\mu(y)\phi_{2}(y)\over \phi(y)}\sum_{w\mid {z\over
y}}\log\min({z/y\over w},w)|  \nonumber \\
& & \leq \sum_{z\leq R}{\mu^{2}(z)\phi((z,j))\over z\phi(z)}
\sum_{y\mid z}{\mu^{2}(y)\phi_{2}(y)\over \phi(y)}d({z\over y})\log z 
\nonumber \\
& & = \sum_{z\leq R}{\mu^{2}(z)\phi((z,j))d(z)\log z\over z\phi(z)}
\prod_{p\mid z}(1+{p-2\over 2(p-1)}) \nonumber  \\
& & \leq \sum_{z\leq R}{\mu^{2}(z)\phi((z,j))3^{\omega(z)}\log z\over z\phi(z)}
\nonumber \\
& & = \sum_{d|j^{*}}\phi(d)\sum_{\stackrel{\scr z\leq R}{\scr (z,j)=d}}
{\mu^{2}(z)3^{\omega(z)}\log z\over z\phi(z)} \nonumber \\
& & = 
\sum_{d|j^{*}}{3^{\omega(d)}\over d}
\sum_{\stackrel{\scr t\leq {R\over d}}{\scr (t,j)=1}}
{\mu^{2}(t)3^{\omega(t)}\over t\phi(t)}\log dt  \\
& & = \sum_{d|j^{*}}{3^{\omega(d)}\log d\over d}
\sum_{\stackrel{\scr t\leq {R\over d}}{\scr (t,j)=1}}
{\mu^{2}(t)3^{\omega(t)}\over t\phi(t)} +
\sum_{d|j^{*}}{3^{\omega(d)}\over d}
\sum_{\stackrel{\scr t\leq {R\over d}}{\scr (t,j)=1}}
{\mu^{2}(t)3^{\omega(t)}\log t\over t\phi(t)}. \nonumber 
\end{eqnarray}
Now note that 
\begin{equation}
\sum_{\stackrel{\scr n=1}{\scr (n,j)=1}}^{\infty}
{\mu^{2}(n)3^{\omega(n)}\over \phi(n) n^{s}}
= \prod_{p\ndiv j}(1+{3\over (p-1)p^{s}}), \qquad (\Re s > 0)
\end{equation}
from which it follows that
\begin{eqnarray}
& & \sum_{\stackrel{\scr n=1}{\scr (n,j)=1}}^{\infty}
{\mu^{2}(n)3^{\omega(n)}\log n\over \phi(n) n} 
= -(\sum_{\stackrel{\scr n=1}{\scr (n,j)=1}}^{\infty}
{\mu^{2}(n)3^{\omega(n)}\over \phi(n) n^{s}})'|_{s=1}   \\
& & = \prod_{p}(1+{3\over p(p-1)}) \prod_{p|j}({p(p-1)\over p^2 -p+3})
\sum_{p\ndiv j}{3\log p\over p^2 -p+3} \ll 1. \nonumber
\end{eqnarray}
So the last line of (8.6) is
\begin{equation}
\ll 1+ \sum_{d|j^{*}}{3^{\omega(d)}\log d\over d} = 1 + \sum_{p|j^{*}}\log p
\sum_{\stackrel{\scr k|j^{*}}{\scr p|k}}{3^{\omega(k)}\over k} ,
\end{equation}
in which the inner sum over $k$ is
\begin{displaymath}
= {3\over p}\sum_{m|{j^{*}\over p}}{3^{\omega(m)}\over m}
= {3\over p}\prod_{\stackrel{\scr q|{j^{*}\over p}}{\scr q: \mathrm{ prime}}}
(1+{3\over q}) = {3\over p+3}
\prod_{\stackrel{q|j^{*}}{\scr q: \mathrm{ prime}}}(1+{3\over q}) .
\end{displaymath}
Hence (8.9) is
\begin{equation}
= 1 + 3\prod_{p|j^{*}}(1+{3\over p})\sum_{p|j^{*}}{\log p\over p+3} 
\ll (\log\log 3 j^{*})^{4}.
\end{equation}
The results of this section are combined to give the evaluation of
(8.1) when $R\to\infty$ as
\begin{equation}
\gs_{2}(j)\log R + O((\log\log 3 j^{*})^{4}) .
\end{equation}

\section{Pure triple correlations: The case $j_1 \neq j_2, j_{1}j_{2} \neq 0$}

We begin the evaluation of (6.12) by employing Lemma 4 to
sum over $a_{12}$ as
\begin{eqnarray}
\lefteqn{
\sum_{\stackrel{\scr a_{12}\leq \min ({R\over a_{13}a_{123}},
{R\over a_{23}a_{123}})}{\scr (a_{12},a_{13}a_{23}a_{123})=1}}
{\mu(a_{12})\mu\cdot\phi((a_{12},j_{1}-j_{2}))\over \phi^{2}(a_{12})} = }
 \nonumber \\
& & \{1-[2\ndiv a_{13}a_{23}a_{123}]\mu((2,j_{1}-j_{2}))\}
C_{2}\prod_{\stackrel{\scr p|a_{13}a_{23}a_{123}}{\scr p>2}}{(p-1)^{2}\over 
p(p-2)}
\prod_{\stackrel{\stackrel{\scr p|j_{1}-j_{2}}{\scr p\ndiv 
a_{13}a_{23}a_{123}}}{\scr p>2}}
({p-1\over p-2}) \nonumber \\
& & \qquad + O({d(j')j'\over \phi(j')
\min ({R\over a_{13}a_{123}},{R\over a_{23}a_{123}})}) ,
\end{eqnarray}
where 
\begin{equation}
j' = {(j_{1}-j_{2})^{*}\over ((j_1 - j_2)^* ,a_{13}a_{23}a_{123})}.
\end{equation}
Let us first consider the main term, which after (9.1) takes the form
\begin{eqnarray}
& & C_{2} \sumprime_{a_{13}a_{23}a_{123} \leq R}
\{1-[2\ndiv a_{13}a_{23}a_{123}]\mu((2,j_{1}-j_{2}))\}
{\mu(a_{13})\mu\cdot\phi((a_{13},j_{1}))\over \phi^{2}(a_{13})}
\nonumber \\
& & \qquad \times \, {\mu(a_{23})\mu\cdot\phi((a_{23},j_{2}))\over
\phi^{2}(a_{23})}
{\mu^{2}(a_{123})\mu\cdot\phi\cdot\phi_{2}((a_{123},j_{1},j_{2}))
\over \phi^{3}(a_{123})} \\
& & \qquad \times \, \phi_{2}({(a_{123},j_{1})\over
(a_{123},j_{1},j_{2})})
\phi_{2}({(a_{123},j_{2})\over (a_{123},j_{1},j_{2})})
\phi_{2}({(a_{123},j_{1}-j_{2})\over
(a_{123},j_{1}-j_{2},j_{1}j_{2})})
\nonumber \\
& & \qquad \times \, \mu\cdot d({a_{123}\over
(a_{123},j_{1}j_{2}(j_{1}-j_{2}))})
\prod_{\stackrel{\scr p|a_{13}a_{23}a_{123}}{\scr p>2}}{(p-1)^{2}\over
p(p-2)}
\prod_{\stackrel{\stackrel{\scr p|j_{1}-j_{2}}{\scr p\ndiv
a_{13}a_{23}a_{123}}}{\scr p>2}}{p-1\over p-2}. \nonumber
\end{eqnarray}
Writing $a_{13} = v,\, a_{23}=w,\, a_{123} = y,\, z=vwy$, turns this
expression into
\begin{eqnarray}
& & C_{2} \sum_{z \leq R}
\{1-[2\ndiv z]\mu((2,j_{1}-j_{2}))\}\mu(z)
\prod_{\stackrel{\scr p|z}{\scr p>2}}{1\over p(p-2)}
\prod_{\stackrel{\stackrel{\scr p|j_{1}-j_{2}}{\scr p\ndiv
z}}{\scr p>2}}({p-1\over p-2}) \nonumber \\
& & 
\times\sum_{y|z}{\mu(y)\mu\cdot\phi\cdot\phi_{2}((y,j_{1},j_{2}))
\phi_{2}({(y,j_{1})\over (y,j_{1},j_{2})})
\phi_{2}({(y,j_{2})\over (y,j_{1},j_{2})})
\phi_{2}({(y,j_{1}-j_{2})\over (y,j_{1}-j_{2},j_{1}j_{2})})\over 
\phi(y)} \nonumber  \\
& & \qquad \times
\prod_{\stackrel{\scr p|y}{\scr p\ndiv j_{1}j_{2}(j_{1}-j_{2})}}(-2)
\sum_{w|{z\over y}}\mu\cdot\phi((w,j_{2}))
\mu\cdot\phi(({z\over yw},j_{1})).
\end{eqnarray}
The innermost sum, over $w$, is
\begin{eqnarray}
& &  = \mu\cdot\phi(({z\over y},j_{1}))\sum_{w|{z\over y}}
{\mu\cdot\phi((w,j_{2}))\over \mu\cdot\phi((w,j_{1}))}  \nonumber \\
& & = \mu\cdot\phi(({z\over y},j_{1}))
\prod_{p|({z\over y},j_{1},j_{2})}2
\prod_{\stackrel{\scr p|({z\over y},j_{1})}{\scr p\ndiv j_{2}}}
(1-{1\over p-1})\prod_{\stackrel{\scr p|({z\over y},j_{2})}{p\ndiv j_{1}}}
(1-(p-1))\prod_{\stackrel{\scr p|{z\over y}}{\scr p\ndiv j_{1}j_{2}}}2 
\nonumber \\
& & = 
\mu\cdot\phi(({z\over y},j_{1}))d(({z\over y},j_{1},j_{2}))
d({{z\over y}\over ({z\over y},j_{1}j_{2})})
{\phi_{2}\over \phi}({({z\over y},j_{1})\over ({z\over
y},j_{1},j_{2})})
\mu\cdot\phi_{2}({({z\over y},j_{2})\over ({z\over y},j_{1},j_{2})}).
\end{eqnarray}
Notice that for even $z$ the contribution to (9.4) of an even $y$ 
is $0$ for any $j_{1}$ and $j_{2}$, and (9.5) shows that for even $z$ an odd
$y$ makes a nonzero contribution only if $j_{1}-j_{2}$ is even. Plugging
(9.5) into (9.4), upon effecting simplifications, (9.4) becomes
\begin{eqnarray}
& & \qquad \gs_{2}(j_{1}-j_{2})\,\times  \\
& & \!\!\!\!\!\!\sum_{z\leq R}\!
{\mu\cdot d(z)\mu\cdot\phi_{2}(({z\over (z,2)},j_{1}))
\mu\cdot\phi_{2}(({z\over (z,2)},j_{2}))\mu\!\cdot\!\phi\!\cdot\!
d((z,j_{1},j_{2}))
\phi_{2}(({z\over (z,2)},j_{1}\! -\! j_{2}))\over
z\phi_{2}({z\over (z,2)})\phi((z,j_{1}-\! j_{2}))
d((z,j_{1}j_{2}))\phi_{2}^{2}(({z\over (z,2)},j_{1},j_{2}))} 
\nonumber \\
& & \!\!\!\!\!\!\!\!\!
\sum_{y|{z\over (z,2)}}\!\!\!{\mu(\!(y,j_{1})\!)\mu(\!(y,j_{2})\!)
\mu(\!(y,j_{1}j_{2}(j_{1}-\! j_{2})\!)\!)\phi_{2}(\!(y,j_{1},j_{2})\!)
\phi_{2}(\!(y,j_{1}-j_{2})\!)d(\!(y,j_{1}j_{2})\!)\over
\phi(y)\phi_{2}(\!(y,j_{1}-\! j_{2},j_{1}j_{2})\!)d(\!(y,j_{1},j_{2})\!)
d(\!(y,j_{1}j_{2}(j_{1}-\! j_{2})\!)\!)}.  \nonumber 
\end{eqnarray}
Recalling the notation used in \S 5,
$J = [j_{1}j_{2}(j_{1}-j_{2})]^{*}$,
$j_{1}^{*} = (j_{1}^{*})'(j_{1}^{*},j_{2}^{*})$,
$j_{2}^{*} = (j_{2}^{*})'(j_{1}^{*},j_{2}^{*}),\, (j_{1} - j_{2})^{*} =
(j_{1}^{*},j_{2}^{*})j_{3}^{*}$, so that
$J=(j_{1}^{*})'(j_{2}^{*})'(j_{1}^{*},j_{2}^{*})j_{3}^{*}$ is a product
of relatively prime factors, the inner sum over $y$ is expressed as
\begin{eqnarray}
& & \sum_{y|{z\over (z,2)}}
{\mu\cdot\phi_{2}((y,j_{1}^{*},j_{2}^{*}))
\mu\cdot\phi_{2}((y,j_{3}^{*}))\over
\phi(y) d((y,j_{1}^{*},j_{2}^{*}))d((y,j_{3}^{*}))} = 
\sum_{y|{z\over (z,2)}}{\mu^{2}(y)\over \phi(y)}
{\mu\cdot\phi_{2}\over d}((y,(j_{1} - j_{2})))  \nonumber \\
& & = \prod_{\stackrel{\scr p|{z\over (z,2)}}{\scr p\ndiv j_{1}-j_{2}}}
(1+{1\over p-1}) \prod_{p|({z\over (z,2)},j_{1}-j_{2})} 
(1- {p-2\over 2(p-1)}) = {z\over \phi(z)d((z,j_{1}-j_{2}))},
\end{eqnarray}
since $j_{1}-j_{2}$ must be even for nonvanishing (9.6).
Hence (9.6) becomes, upon further simplifications,
\begin{eqnarray}
& & \qquad \gs_{2}(j_{1}-j_{2})\,\times \\
& & \!\!\!\sum_{z\leq R}\!
{\mu\!\cdot\! d(z)\mu\!\cdot\!\phi_{2}(({z\over (z,2)},j_{1}))
\mu\!\cdot\!\phi_{2}(({z\over (z,2)},j_{2}))\mu\!\cdot\!\phi\!\cdot\!
d((z,j_{1},j_{2}))
\phi_{2}(({z\over (z,2)},j_{1}- \! j_{2}))\over 
\phi(z)\phi_{2}({z\over (z,2)})\phi\!\cdot\! d((z,j_{1}-\! j_{2}))
d((z,j_{1}j_{2}))\phi_{2}^{2}(({z\over (z,2)},j_{1},j_{2}))} 
\nonumber  \\
& & = \gs_{2}(j_{1}-j_{2})\sum_{z\leq R}
{\mu(z) d(z)\over \phi(z)\phi_{2}({z\over (z,2)})}
{\mu\over d}((z,J)){\mu\over \phi}((z,j_{3}^{*}))
\phi_{2}(({z\over (z,2)},J)).  \nonumber
\end{eqnarray}
The last sum has been settled in Lemma 5, so that (9.4), that is
the main term, is evaluated as
\begin{eqnarray}
& & = \gs_{2}(j_{1}-j_{2})\, 2\, 
[2\ndiv j_{3}^{*}]
\prod_{p\ndiv J}{p(p-3)\over (p-1)(p-2)}
\prod_{\stackrel{\scr p>2}{\scr p|j_{1}j_{2}}}{p\over p-1}
\prod_{\stackrel{\scr p>2}{\scr p|j_{3}^{*}}}{p(p-2)\over (p-1)^{2}}
\nonumber \\
& & = [3|J][2|j_{1}-j_{2}][2\ndiv j_{3}^{*}]
\, 4C_{2} 
\prod_{p\ndiv J}{p(p-3)\over
(p-1)(p-2)} \prod_{\stackrel{\scr p|J}{\scr p>2}}{p\over p-1}
\prod_{\stackrel{\scr p|(j_{1}^{*},j_{2}^{*})}{\scr p>2}}
{p-1\over p-2} \nonumber \\
& & = [3|J][2|j_{1}-j_{2}][2\ndiv j_{3}^{*}]
\, ( 2C_{2} 
\prod_{\stackrel{\scr p|(j_{1}^{*},j_{2}^{*})}{\scr p>2}}
{p-1\over p-2} )\, C_{3}\, 3
\prod_{\stackrel{\scr p|J}{\scr p>3}}{p-2\over p-3} \nonumber \\
& & = \gs_{3}(J)\gs_{2}((j_{1}^{*},j_{2}^{*})) = \gs((0,j_{1},j_{2})),
\end{eqnarray}
and there is an error of $O(R^{-1+\epsilon})$. 

It remains to carry out the follow-up of the error term (9.1) in (6.12),
which is
\begin{eqnarray}
& \ll & {(j_{1}-j_{2})^{*}d((j_{1}-j_{2})^{*})\over
\phi((j_{1}-j_{2})^{*})} 
\sumprime_{z=vwy\leq
R}{\mu^{2}(z)\phi((v,j_{1}))\phi((w,j_{2}))\phi((y,j_{1},j_{2}))
\over \phi^{2}(z)\phi(y)\min({R\over vy},{R\over wy})} \nonumber \\
& & \times \; \phi_{2}((y,j_{1}))\phi_{2}({(y,j_{2})\over
(y,j_{1},j_{2})})\phi_{2}({(y,j_{1}-j_{2})\over
(y,j_{1}-j_{2},j_{1}j_{2})})d({y\over (y,(j_{1}-j_{2})j_{1}j_{2})})
\nonumber \\
& \ll & 
R^{-1+\epsilon}\sum_{z\leq R}
{\mu^{2}(z)z\phi((z,J))\phi_{2}(y)\over \phi^{2}(z)\phi(y)\min(w,v)}
\nonumber \\
& = & 
R^{-1+\epsilon}\sum_{z\leq R}{\mu^{2}(z)z\phi((z,J))\over \phi^{2}(z)}
\sum_{y|z}{\phi_{2}(y)\over \phi(y)}\sum_{w|{z\over y}}{1\over
\min(w,{z\over yw})} \nonumber \\
& \leq & 
R^{-1+\epsilon}\sum_{z\leq R}{\mu^{2}(z)z\phi((z,J))d(z)\over \phi^{2}(z)}
\sum_{y|z}{\phi_{2}(y)\over \phi(y)d(y)} \nonumber \\
& = & 
R^{-1+\epsilon}\sum_{z\leq R}\mu^{2}(z)\phi((z,J))\prod_{p|z}{p(3p-4)\over
(p-1)^{3}} \nonumber \\
& \ll & 
R^{-1+\epsilon}\sum_{z\leq R}\mu^{2}(z)\prod_{p|z}{p(3p-4)\over
(p-1)^{3}} \nonumber \\
& \ll & R^{-1+\epsilon}. \nonumber \\
&  & \mbox{} 
\end{eqnarray}
This completes the proof of Theorem 1 in the case considered in this
section.

\section{Preliminaries for proof of Theorem 3}
Consider the quantities ${\mathcal{M}}_{k}'(N,h,\psi_{R},\rho,C,A)$, say
$\calm_{k}'$ for brevity, defined as
\begin{equation}
\calm_{k}' = 
\sum_{n=N+1}^{2N}(\psi_R(n+h)-\psi_R(n)-h -CA)^{k-1}
(\psi(n+h)-\psi(n)-h-\rho A),
\end{equation}
where $\rho$ and $C$ are to be constants, and $A=o(h)$ will be chosen
appropriately. The $'$ on $\calm_{k}'$ (and on $M_{k}',\, 
\mathcal{S}_{k}'$) signifies that the sum over $n$ runs from $N+1$ to
$2N$. Our main interest here is the case $k=3$.
The Generalized Riemann Hypothesis will be assumed in
this section. 

For $k=1$ we have by (1.37),
\begin{eqnarray}
\calm_{1}' & = & \sum_{n=N+1}^{2N}(\psi(n+h)-\psi(n)-h-\rho A) \\
& = & M_{1}'(N,h,\psi) - (h+\rho A)N = -\rho AN + 
O(N^{{1\over 2}}h\log^{2}N), \nonumber
\end{eqnarray}
valid for $h\leq N$ and $N\to \infty$. The $k=2$ case is more
illuminating in that it helps us see the appropriate choice for $A$.
We have
\begin{eqnarray}
\calm_{2}' & = & \sum_{n=N+1}^{2N}(\psi_R(n+h)-\psi_R(n)-h -CA)
(\psi(n+h)-\psi(n)-h-\rho A) \nonumber \\
& = & \tilde{M}_{2}'(N,h,\psi_{R}) - (h+\rho
A)M_{1}'(N,h,\psi_{R}) \nonumber \\
&  & \mbox{} \;\; - (h+ CA)M_{1}'(N,h,\psi) + (h+\rho A)(h+CA)N . 
\end{eqnarray}
Eq.s (1.28) and (4.1) give
\begin{equation}
M_{1}'(N,h,\psi_{R}) = \sum_{1\leq j_{1} \leq h}
\mathcal{ S}_{1}'(N,(j_{1}),(1)) = Nh + O(Rh) .
\end{equation}
By (1.38) we have
\begin{equation}
\tilde{M}_{2}'(N,h,\psi_{R})\! = \!\mathcal{L}_{1}(R) \!\!
\sum_{1\leq j_{1} \leq h}\!\tilde{\mathcal{ S}}_{1}'
(N,(j_{1}),(1)) \! + \!
\sum_{\stackrel{\scr 1\le j_1 ,j_2 \le h}{ \scr \mathrm{ distinct}}}
\!\!\tilde{\mathcal{ S}}_{2}'(N,(j_{1},j_{2}),(1,1)) \!+ \! O(RN^{\e}).
\end{equation}
Here the first sum is $M_{1}'(N,h,\psi)$ which we know. In the second sum
we cannot use the evaluation (5.6) since the contribution of error term
coming from the Bombieri-Vinogradov theorem will be greater than the main
term of $\calm_{2}'$ when $h$ is a power of $N$. We begin by proceeding 
similar to (5.1)-(5.6), and we write $j= j_2 - j_1$, so that
\begin{eqnarray}
& & \sum_{\stackrel{\scr 1\le j_1 ,j_2 \le h}{ \scr \mathrm{ distinct}}}
\tilde{\mathcal{ S}}_{2}'(N,(j_{1},j_{2}),(1,1)) =
N\sum_{\stackrel{\scr 1\le j_1 ,j_2 \le h}{ \scr \mathrm{ distinct}}}
\sum_{r\le R}{\mu^{2}(r)\over \phi(r)}\sum_{\stackrel{\scr d|r}
{\scr (d,j_1 - j_2 )=1}}{d\mu(d)\over \phi(d)} \nonumber \\
& & \qquad \qquad \qquad + 
O\bigl(\sum_{\stackrel{\scr 1\le j_1 ,j_2 \le h}{ \scr \mathrm{distinct}}}
\sum_{r\leq R}{\mu^{2}(r)\over \phi(r)}\sum_{d|r}d\max_{x\leq 2N+h}
|E(x;d,j_2 - j_1)| \bigr) \nonumber \\
& & = N\sum_{\stackrel{\scr 1\le j_1 ,j_2 \le h}{ \scr \mathrm{ distinct}}}
[\gs_{2}(j_2 - j_1) + O\left({j^{*}d(j^{*})\over R\phi(j^{*})}\right)] 
\nonumber \\
& & \qquad + O\bigl(\sum_{1\leq j_1 \leq h}
\sum_{r\leq R}{\mu^{2}(r)\over \phi(r)}\sum_{d|r}d
\sum_{1\leq j \leq h-1}\max_{x\leq 2N+h}|E(x;d,j)| \bigr) 
\nonumber \\
& & = 2N\sum_{1\leq j \leq h-1}\sum_{1\leq j_1 \leq h-j}[\gs_{2}(j) +
O\left({jd(j)\over R\phi(j)}\right)] \nonumber \\
& & \qquad + O\bigl(\sum_{1\leq j_1 \leq h}
\sum_{r\leq R}{\mu^{2}(r)\over \phi(r)}\sum_{d|r}d
(\sum_{1\leq j \leq h-1}1)^{{1\over 2}}
(\sum_{1\leq j \leq h-1}\max_{x\leq 2N+h}|E(x;d,j)|^{2})^{{1\over
2}}\bigr) \nonumber \\
& & =  2N\sum_{1\leq j \leq h}(h-j)\gs_{2}(j)+O({Nh^{2}\log h\over R}) 
\nonumber \\
& & \qquad + O\bigl(h^{{3\over 2}}
\sum_{r\leq R}{\mu^{2}(r)\over \phi(r)}\sum_{d|r}d
(\sum_{1\leq j \leq h-1}\max_{x\leq 2N+h}|E(x;d,j)|^{2})^{{1\over
2}}\bigr) .
\end{eqnarray}
For the calculation of the main term we know from \cite{FG} that 
\begin{equation}
\sum_{1\leq j \leq h}(h-j)\gs_{2}(j) = {h^{2}\over 2}-{h\log h\over 2}
+ {1-\gamma-\log 2\pi\over 2}h + O(h^{{1\over 2}+\e}) .
\end{equation}
The sum
of the $E$'s is $O((1+{h\over d})N\log^{4}N)$ by Hooley's estimate (1.47) 
which depends on GRH, so that
\begin{eqnarray}
\sum_{\stackrel{\scr 1\le j_1 ,j_2 \le h}{ \scr \mathrm{ distinct}}}
\!\tilde{\mathcal{ S}}_{2}'(N,(j_{1},j_{2}),(1,1))\!\! & = & 
\! Nh^2 \!+\! Nh(-\log h +1-
\gamma -\!\log2\pi)\!+\!O(Nh^{{1\over 2}+\e}) \nonumber \\
+\, O({Nh^{2}\log h\over R})\!\!\!\! & & \!\!
+ O(N^{{1\over 2}}h^{{3\over 2}}R\log^{2}N)
+ O(N^{{1\over 2}}h^{2}R^{{1\over 2}}\log^{2}N).
\end{eqnarray}
Using (1.37) and (10.8) in (10.5) we have
\begin{eqnarray}
& & \!\!\!\!\tilde{M}_{2}'(N,h,\psi_{R}) = Nh^2 + 
Nh(\mathcal{L}_{1}(R) -\log h +1 -\gamma - \log2\pi) 
+ O(Nh^{{1\over 2}+\e}) \nonumber \\
& & \qquad + O({Nh^{2}\log h\over R}) 
+ O(N^{{1\over 2}}h\log^{3}N) 
+ O((N^{{1\over 2}}h^{{3\over 2}}R^{{1\over 2}}\log^{2}N)
(h^{{1\over 2}} + R^{{1\over 2}}))
\end{eqnarray}
and plugging (10.9), (10.4), and (1.37) in (10.3) we obtain
\begin{eqnarray}
\calm_{2}'&  = & Nh(\mathcal{L}_{1}(R) -\log h +1 -\gamma - \log2\pi) 
+ \rho C N A^2 + O(Rh^2 ) + O(Nh^{{1\over 2}+\e}) \nonumber \\
& & \qquad + O((N^{{1\over 2}}h^{{3\over 2}}R^{{1\over 2}}\log^{2}N)
(h^{{1\over 2}} + R^{{1\over 2}}))
+ O({Nh^{2}\log h\over R})\nonumber \\
& & \qquad + O((N^{{1\over 2}}h\log^{2}N (h+\log N)) 
+ O(CAN^{{1\over 2}}h\log^{2}N) + O(\rho ARh). 
\end{eqnarray}
Here the term $\rho CNA^2$ will be of the same order of magnitude with 
the main term if we choose 
\begin{equation}
A = (h\log N)^{{1\over 2}}.
\end{equation}
Since (10.1) would be meaningful for $A=o(h)$, we must have $\log N = o(h)$.
We need to have $R=o(N)$ so that (4.1) has an asymptotic interpretation,
and $R \gg N^{\e}$ so that $\mathcal{L}_{1}(R)$ is of the same order of
magnitude with $\log N$. We also require that the error term 
$\displaystyle O({Nh^{2}\log h\over R})$ to be smaller than $Nh$, and this
imposes $h\log h=o(R)$. With these in mind (10.10) reduces to
\begin{eqnarray}
\calm_{2}'&  = & Nh(\rho C\log N + \mathcal{L}_{1}(R) -\log h +1 -\gamma -
\log2\pi) + O(Nh^{{1\over 2}+\e}) \nonumber \\
& & \qquad + O({Nh^{2}\log h\over R})
+ O(N^{{1\over 2}}h^{{3\over 2}}R\log^{2}N). \nonumber \\
\end{eqnarray}
This will have asymptotic significance when
the very last error term is smaller than $Nh$, i.e. for
\begin{equation}
h^{{1\over 2}}R = o({N^{{1\over 2}}\over \log^{2}N}),
\end{equation}
which restricts us to
\begin{equation}
h = o({N^{{1\over 3}}\over \log^{2}N}).
\end{equation}
When all these conditions are met we have
\begin{equation}
\calm_{2}' \sim  Nh(\rho C\log N + \log {R\over h}).
\end{equation}
If it were that
$\psi(n+h)-\psi(n)-h = o((h\log n)^{{1\over 2}})$, then (10.15) 
with $C=0$ would imply
\begin{displaymath}
\sum_{n=N+1}^{2N}(\psi_R(n+h)-\psi_R(n) - h) \sim -{Nh^{{1\over 2}}\over
\rho}{\log {R\over h}\over (\log N)^{{1\over 2}}}
\end{displaymath}
for any fixed $\rho \neq 0$, which is absurd since the left-hand side
doesn't depend on $\rho$. Thus we obtain
\begin{equation}
\max_{x\leq y \leq 2x}|\psi(y+h)-\psi(y)-h| \gg_{\e} (h\log x)^{{1\over 2}}
\end{equation}
for $1 \leq h\leq x^{{1\over 3}-\e}$, a result which is already implicit in
previous works (e.g. \cite{GY1}) involving the lower bound method depending
on the second order correlations of $\lambda_{R}(n)$'s.

\section{The Third Mixed Moment for the Proof of Theorem 3}
We now turn our attention to 
\begin{eqnarray}
\calm_{3}' & = & \sum_{n=N+1}^{2N}(\psi_R(n+h)-\psi_R(n)-h -CA)^{2}
(\psi(n+h)-\psi(n)-h-\rho A) \nonumber \\
& = & \tilde{M}_{3}'(N,h,\psi_{R}) - (h+\rho A)M_{2}'(N,h,\psi_{R}) 
-2(h+CA)\tilde{M}_{2}'(N,h,\psi_{R}) \nonumber \\
& & \;\;\; +2(h+CA)(h+\rho A)M_{1}'(N,h,\psi_{R}) 
+ (h+ CA)^{2}M_{1}'(N,h,\psi) 
\nonumber \\
& & \;\;\;\; - (h+\rho A)(h+CA)^{2}N . 
\end{eqnarray}
We now need to evaluate $M_{2}'(N,h,\psi_{R})$ and
$\tilde{M}_{3}'(N,h,\psi_{R})$. By (1.28) we have
\begin{equation}
{M}_{2}'(N,h,\psi_{R}) = 
\sum_{1\leq j_{1} \leq h}\mathcal{ S}_{2}'(N,(j_{1}),(2)) + 
2\sum_{1\le j_1 < j_2 \le h}\mathcal{ S}_{2}'(N,(j_{1},j_{2}),(1,1)).
\end{equation}
The first sum on the right is evaluated by (4.7) as
\begin{equation}
\sum_{1\leq j_{1} \leq h}\mathcal{ S}_{2}'(N,(j_{1}),(2)) =
Nh\mathcal{L}_{1}(R) + O(hR^2 ) .
\end{equation}
For the second sum, letting $j= j_2 -j_1$, we have by (4.9)
\begin{equation}
2\sum_{1\le j_1 < j_2 \le h}\mathcal{ S}_{2}'(N,(j_{1},j_{2}),(1,1)) =
2N\sum_{1\leq j \leq h}(h-j)\gs_{2}(j) + O({Nh^{2}\log h\over R}) + 
O(h^2 R^2)
\end{equation}
(it is more convenient to keep the sum of $\gs_2$'s as is, not
using (10.7) until the end, and also view the first line of (10.9) except
for the term $Nh\mathcal{L}_{1}(R)$ in this manner). Hence we obtain
\begin{equation}
{M}_{2}'(N,h,\psi_{R}) =
2N\sum_{1\leq j \leq h}(h-j)\gs_{2}(j) + Nh\mathcal{L}_{1}(R) 
+ O({Nh^{2}\log h\over R}) + O(R^2 h^2).
\end{equation}

By (1.38) we have
\begin{eqnarray}
\tilde{M}_{3}'(N,h,\psi_{R}) & = & (\mathcal{L}_{1}(R))^2 
\sum_{1\leq j_{1} \leq h}\tilde{\mathcal{ S}}_{1}'(N,(j_{1}),(1)) 
+ \sum_{\stackrel{\scr 1\le j_1 ,j_2 \le h}{ \scr \mathrm{ distinct}}}
\tilde{\mathcal{ S}}_{3}'(N,(j_{1},j_{2}),(2,1)) \nonumber \\
& & \qquad + 2\mathcal{L}_{1}(R)
\sum_{\stackrel{\scr 1\le j_1 ,j_2 \le h}{ \scr \mathrm{ distinct}}}
\tilde{\mathcal{ S}}_{2}'(N,(j_{1},j_{2}),(1,1)) \nonumber \\
& & \qquad + 
\sum_{\stackrel{\scr 1\le j_1 ,j_2, j_3 \le h}{ \scr \mathrm{ distinct}}}
\tilde{\mathcal{ S}}_{3}'(N,(j_{1},j_{2},j_{3}),(1,1,1)) 
+ O(RN^{\e}). 
\end{eqnarray}
In calculating the sum of
$\tilde{\mathcal{ S}}_{3}'(N,(j_{1},j_{2}),(2,1))$'s over $j_1,\, j_2$
it is not suitable to use Theorem 2 because we will pick up an
error term $O(Nh^{2}\log\log^{2} 3j^{*})$ arising from the 
$\log\min$ term of (5.20), and this will be greater than the eventual
main term due to cancellations of larger terms
(note that $j_1,\, j_2$ are not the same as in \S 5). 
Therefore we start anew as
\begin{eqnarray}
& & \!\!\!\!
\sum_{\stackrel{\scr 1\le j_1 ,j_2 \le h}{ \scr \mathrm{ distinct}}}
\tilde{\mathcal{ S}}_{3}'(N,(j_{1},j_{2}),(2,1)) = 
\sum_{\stackrel{\scr 1\le j_1 ,j_2 \le h}{ \scr \mathrm{ distinct}}}
\sum_{n=N+1}^{2N}\lambda_{R}^{2}(n+j_1 )\Lambda(n+j_2 ) \nonumber \\
& = & \sum_{\stackrel{\scr 1\le j_1 ,j_2 \le h}{ \scr \mathrm{ distinct}}}
\sum_{r_1, r_2 \leq R}{\mu^{2}(r_1 )\mu^{2}(r_2 )\over \phi(r_1 )\phi(r_2 )}
\sum_{\stackrel{\scr d|r_1 }{\scr e|r_2 }}d\mu(d)e\mu(e)
\sum_{\stackrel{\scr n=N+1}{\scr n\equiv -j_{1}(\mathrm{mod}\, [d,e])}}
^{2N}\Lambda(n+j_2 ) \nonumber \\
& = & \sum_{\stackrel{\scr 1\le j_1 ,j_2 \le h}{ \scr \mathrm{ distinct}}}
\sum_{r_1, r_2 \leq R}{\mu^{2}(r_1 )\mu^{2}(r_2 )\over \phi(r_1 )\phi(r_2 )}
\;\; \bigl\{ \sum_{\stackrel{\stackrel{\scr d|r_1 }{\scr e|r_2 }}
{\scr ([d,e],j_2 -j_1)=1}}{N\over \phi([d,e])}d\mu(d)e\mu(e) \nonumber \\
& & \hspace{5.2cm} + O(\sum_{\stackrel{\scr d|r_1 }{\scr e|r_2 }}de
\max_{u\leq 2N+h}|E(u;[d,e],j_2 - j_1)|) \bigr\} \nonumber \\
& = & N\sum_{1\leq |j| \leq h}(h-|j|)
\sum_{r_1, r_2 \leq R}{\mu^{2}(r_1 )\mu^{2}(r_2 )\over \phi(r_1 )\phi(r_2 )}
\sum_{\stackrel{\stackrel{\scr d|r_1 }{\scr e|r_2 }}
{\scr ([d,e],j)=1}}{d\mu(d)e\mu(e)\over \phi([d,e])} \nonumber  \\
& & \qquad  + 
O(h\sum_{r_1, r_2 \leq R}{\mu^{2}(r_1 )\mu^{2}(r_2 )\over \phi(r_1 )\phi(r_2 )}
\sum_{\stackrel{\scr d|r_1 }{\scr e|r_2 }}de
\sum_{1\leq |j| \leq h}\max_{u\leq 2N+h}|E(u;[d,e],j)|). 
\end{eqnarray}
The last error term is dealt with in the same way as in 
(10.6), first by applying Cauchy-Schwarz inequality to
the sum over $j_1 $ and then using 
Hooley's GRH-dependent estimate (1.47), which makes it
\begin{eqnarray*}
& \ll & 
N^{{1\over 2}}h^{{3\over 2}}\log^{2}N 
\sum_{r_1, r_2 \leq R}{\mu^{2}(r_1)\mu^{2}(r_2 )\over \phi(r_1 )\phi(r_2 )}
\sum_{\stackrel{\scr d|r_1 }{\scr e|r_2 }}de
(1+({h\over [d,e]})^{{1\over 2}}) \\
& \ll & N^{{1\over 2}}h^{{3\over 2}}R^{2}\log^{2}N + 
N^{{1\over 2}}h^{2}\log^{2}N 
\sum_{r_1, r_2 \leq R}{\mu^{2}(r_1)\mu^{2}(r_2 )\over \phi(r_1 )\phi(r_2 )}
\sum_{\stackrel{\scr d|r_1 }{\scr e|r_2 }}{de\over ([d,e])^{{1\over 2}}}
\end{eqnarray*}
The last sum is treated similar to what follows (4.5), and with the notation 
used there the inner sums over $d$ and $e$ become
\begin{eqnarray*}
& & \sum_{\delta|(r_{1},r_{2})}\delta^{{3\over 2}}
\sum_{d'|{r{_1}\over \delta}}(d')^{{1\over 2}}
\sum_{\stackrel{\scr e'|{r_{2}\over \delta}}{\scr (e',d')=1}}(e')^{{1\over 2}}\\
& = & \prod_{p|r_2 }(1+\sqrt{p})\sum_{\delta|(r_{1},r_{2})}\prod_{p|\delta }
{p^{{3\over 2}}\over 1+\sqrt{p}}\sum_{d'|{r{_1}\over \delta}}(d')^{{1\over 2}}
\prod_{p|(r_2 , d')}{1\over 1+\sqrt{p}} \\
& = & \prod_{p|r_2 }(1+\sqrt{p})\sum_{\delta|(r_{1},r_{2})}\prod_{p|\delta }
{p^{{3\over 2}}\over 1+\sqrt{p}}
\prod_{p|{(r_1 , r_2 )\over \delta}}
{1+2\sqrt{p}\over 1+\sqrt{p}}\prod_{p|{r_{1}\over (r_1 ,r_2)}}(1+\sqrt{p}) \\
& = & \prod_{p|r_1 }(1+\sqrt{p})\prod_{p|r_2 }(1+\sqrt{p})
\prod_{p|(r_1 , r_2)}(1+2\sqrt{p} + p^{{3\over 2}}) .
\end{eqnarray*}
Hence the sum we are concerned with is now expressed as
\begin{displaymath}
\sum_{r_1, r_2 \leq R}{\mu^{2}(r_1)\mu^{2}(r_2 )\over \phi(r_1 )\phi(r_2)}
\prod_{p|r_1 }(1+\sqrt{p})\prod_{p|r_2 }(1+\sqrt{p})
\prod_{p|(r_1 , r_2)}(1+2\sqrt{p} + p^{{3\over 2}}) ,
\end{displaymath}
and writing $r_1 = a_1 a_{12}, \, r_2 = a_2 a_{12}$ with $a_{12}=(r_1 ,
r_2)$ this becomes
\begin{eqnarray*}
& & \sumprime_{\stackrel{\scr a_1 a_{12}\leq R}{\scr a_2 a_{12}\leq R}}
\mu^{2}(a_1 )\mu^{2}(a_2 )\mu^{2}(a_{12})\prod_{p|a_1 a_2 }{1\over \sqrt{p}-1}
\prod_{p|a_{12}}{1+2\sqrt{p}+p^{{3\over 2}}\over (p-1)^{2}} \\
& \ll & \sum_{a_{12}\leq R}
\prod_{p|a_{12}}{1+2\sqrt{p}+p^{{3\over 2}}\over (p-1)^{2}} 
(\sum_{a\leq {R\over a_{12}}}\prod_{p|a}{1\over \sqrt{p}-1})^2 \\
& \ll & R\sum_{a_{12}\leq R}
\prod_{p|a_{12}}{1+2\sqrt{p}+p^{{3\over 2}}\over p(p-1)^{2}} \ll R .
\end{eqnarray*}
Thus the error term of (11.7) has been shown to be
\begin{equation}
\ll N^{{1\over 2}}h^{{3\over 2}}R^{2}\log^{2}N + N^{{1\over 2}}h^{2}R\log^{2}N .
\end{equation}
We now treat the main term of (11.7), which is equal to
\begin{equation}
2N\sum_{r_1, r_2 \leq R}{\mu^{2}(r_1)\mu^{2}(r_2)\over \phi(r_1)\phi(r_2)}
\sum_{\stackrel{\scr d|r_1 }{\scr e|r_2 }}{d\mu(d)e\mu(e)\over \phi([d,e])} 
\sum_{\stackrel{\scr 1\leq j \leq h}{\scr (j,[d,e])=1}}(h-j) .
\end{equation}
Here the innermost sum is
\begin{eqnarray}
& & \sum_{\stackrel{\scr 1\leq j \leq h}{\scr (j,[d,e])=1}}(h-j) = 
\sum_{1\leq j \leq h}(h-j)\sum_{k|(j,[d,e])}\mu(k) 
= \sum_{k|[d,e]}\mu(k)\sum_{\stackrel{\scr 1\leq j \leq h}{\scr k|j}}(h-j) 
\nonumber \\
& = & \sum_{k|[d,e]}\mu(k)\{h({h\over k} + O(1)) - k\sum_{1\leq j' \leq 
{h\over k}}j' \} = \sum_{k|[d,e]}\mu(k)({h^{2}\over 2k} + O(h)) \nonumber \\
& & \qquad \qquad = {h^{2}\over 2}{\phi([d,e])\over [d,e]} + O(hd([d,e])) .
\end{eqnarray}
The main term of (11.10) plugged in (11.9) gives
\begin{equation}
Nh^{2}\sum_{r_1, r_2 \leq R}{\mu^{2}(r_1)\mu^{2}(r_2)\over \phi(r_1)\phi(r_2)}
\sum_{\stackrel{\scr d|r_1 }{\scr e|r_2 }}{d\mu(d)e\mu(e)\over [d,e]} =
Nh^{2}\mathcal{L}_{1}(R)
\end{equation}
where the last evaluation is by virtue of (4.5)--(4.7). It remains to
consider the contribution of the error term of (11.10) plugged in (11.9),
\begin{equation}
Nh\sum_{r_1, r_2 \leq R}{\mu^{2}(r_1)\mu^{2}(r_2)\over \phi(r_1)\phi(r_2)}
\sum_{\stackrel{\scr d|r_1 }{\scr e|r_2 }}{ded([d,e])\over \phi([d,e])} .
\end{equation}
Keeping the notation we have been using
the inner sums over $d$ and $e$ are rewritten as
\begin{eqnarray*}
&  & \sum_{\delta|(r_{1},r_{2})}{\delta^{2}d(\delta)\over \phi(\delta)}
\sum_{d'|{r{_1}\over \delta}}{d' d(d')\over \phi(d')}
\sum_{\stackrel{\scr e'|{r_{2}\over \delta}}{\scr (e',d')=1}}
{e' d(e')\over \phi(e')} \\
& = & \prod_{p|r_2 }{3p-1\over p-1}\sum_{\delta|(r_1 , r_2)}
\prod_{p|\delta}{2p^{2}\over 3p-1}
\sum_{d'|{r{_1}\over \delta}}{d' d(d')\over \phi(d')}
\prod_{p|(r_2 , d')}{p-1\over 3p-1} \\
& = & \prod_{p|r_2 }{3p-1\over p-1}\sum_{\delta|(r_1 , r_2)}
\prod_{p|\delta}{2p^{2}\over 3p-1}\prod_{p|{(r_1 , r_2 )\over \delta}}
{5p-1\over 3p-1}\prod_{p|{r_{1}\over (r_1 ,r_2)}}{3p-1\over p-1} \\
& = & \prod_{p|r_1 }{3p-1\over p-1}\prod_{p|r_2 }{3p-1\over p-1}
\prod_{p|(r_1 , r_2)}{2p^2 +5p-1\over p-1}.
\end{eqnarray*}
Hence (11.12) becomes
\begin{eqnarray}
& & \qquad
Nh\sum_{r_1, r_2 \leq R}{\mu^{2}(r_1)\mu^{2}(r_2)\over \phi(r_1)\phi(r_2)}
\prod_{p|r_1 }{3p-1\over p-1}\prod_{p|r_2 }{3p-1\over p-1}
\prod_{p|(r_1 , r_2)}{2p^2 +5p-1\over p-1}  \\
& = & Nh\sumprime_{\stackrel{\scr a_1 a_{12}\leq R}{\scr a_2 a_{12}\leq R}}
\mu^{2}(a_1 )\mu^{2}(a_2 )\mu^{2}(a_{12})
\prod_{p|a_1 a_2 }{3p-1\over (p-1)^2 }
\prod_{p|a_{12}}{2p^2 +5p-1\over (p-1)^{3}} \nonumber \\
& \ll & Nh\sum_{a_{12}\leq R}\prod_{p|a_{12}}{2p^2 +5p-1\over (p-1)^{3}} 
(\sum_{a\leq {R\over a_{12}}}\prod_{p|a}{3p-1\over (p-1)^2 })^2 \nonumber \\
& \ll & Nh\log^{6}N 
\sum_{a_{12}\leq R}\prod_{p|a_{12}}{2p^2 +5p-1\over (p-1)^{3}} 
\ll Nh\log^{8}N. \nonumber 
\end{eqnarray}
This completes the evaluation begun in (11.7), giving
\begin{eqnarray}
\sum_{\stackrel{\scr 1\le j_1 ,j_2 \le h}{ \scr \mathrm{ distinct}}}
\tilde{\mathcal{ S}}_{3}'(N,(j_{1},j_{2}),(2,1)) & = &
Nh^{2}\mathcal{L}_{1}(R) +O(Nh\log^{8}N) \nonumber \\
& & \qquad + O(N^{{1\over 2}}h^{{3\over 2}}R^{2}\log^{2}N) + 
O(N^{{1\over 2}}h^{2}R\log^{2}N).
\end{eqnarray}

\section{Completion of the Proof of Theorem 3}
The last quantity involving the correlations that remains to be considered
for the mixed third moment is
\begin{eqnarray}
\lefteqn{
\sum_{\stackrel{\scr 1\le j_1 ,j_2, j_3 \le h}{ \scr \mathrm{ distinct}}}
\tilde{\mathcal{ S}}_{3}'(N,(j_{1},j_{2},j_{3}),(1,1,1)) } \nonumber \\
& & 
= \sum_{\stackrel{\scr 1\le j_1 ,j_2 ,j_3 \le h}{ \scr \mathrm{ distinct}}}
\sum_{n=N+1}^{2N}\lambda_{R}(n+j_1 )\lambda_{R}(n+j_2)\Lambda(n+j_3) .
\end{eqnarray}
The inner sum here is the same as (5.7) except for a shift, 
and the innermost sum of (5.8) adapted to the present situation is
\begin{equation}
\sum_{\stackrel{\stackrel{\scr n=N+1+j_3 }{\scr n\equiv 
j_{3}-j_{1}(\mathrm{mod}\, d)}}
{\scr n\equiv j_{3}-j_{2}(\mathrm{mod}\, e)}}
^{2N+ j_3 }\Lambda(n) ,
\end{equation}
where the divisibility conditions, compatible only when $(d,e)|j_1 - j_2 $, 
can be combined as $n\equiv j(\mathrm{mod}\, [d,e])$ say. The last sum
is equal to
\begin{equation}
[([d,e],j)=1]{N\over \phi([d,e])} + E(2N+j_3 ;[d,e],j)- E(N+1+j_3 ;[d,e],j).
\end{equation}
We get the same main term as that of (5.11) except that the conditions now
read $(d,e)|j_1 -j_2 ,\, (d,j_3 -j_1 )=1,\, (e,j_3 -j_2 )=1$, and the
calculation carried out in \S 5 evaluates it as
\begin{equation}
N\gs((j_{1},j_{2},j_{3})) + O(NR^{-1+\e}) .
\end{equation}
We shall take up the contribution of (12.4) in (12.1) after considering
the contibution of the error terms of (12.3) to (12.1) which can be 
majorized as
\begin{equation}
\ll \sum_{\stackrel{\scr 1\le j_1 ,j_2 ,j_3 \le h}{ \scr \mathrm{ distinct}}}
\sum_{r_1, r_2 \leq R}{\mu^{2}(r_1)\mu^{2}(r_2)\over \phi(r_1)\phi(r_2)}
\sum_{\stackrel{\stackrel{\scr d|r_1 }{\scr e|r_2 }}{\scr (d,e)|j_1 -j_2 }}
de\max_{u\leq 2N+h}|E(u;[d,e],j)|.
\end{equation}
To simplify the expressions we may interchange the order of the two double
sums in (12.5) as in the first line of (5.10). This switching costs a factor
of $\log^{2}R$, which is unimportant in our application. (On two occasions
in \S 11 we didn't resort to this interchange). Since $j$ is a
function of $k_1 = j_3 -j_1 $ and $k_2 = j_3 - j_2 $ (and also of $d,\, e$),
we can rearrange the summations so as to see that (12.5) is 
\begin{equation}
\ll \log^{2}R \sum_{1\leq j_3 \leq h}
\sum_{d,\, e \leq R}{\mu^{2}(d)\mu^{2}(e)de\over \phi(d)\phi(e)} 
\sum_{\stackrel{\stackrel{\scr j_3 -h \leq k_1 ,\, k_2 \leq j_3 - 1}
{\scr (k_2 -k_1 )k_1 k_2 \neq 0}}{\scr (d,e)|k_2 -k_1 }}
\max_{u\leq 2N+h}|E(u;[d,e],j)|.
\end{equation}
Note that if $(d,e) \geq h$, then the innermost sum is void. The number
of permissible pairs of $k_{1},\, k_2 $ is 
$\displaystyle \ll h(\lceil {h\over (d,e)}\rceil - 1)$. Upon applying
Cauchy-Schwarz inequality and Hooley's GRH-dependent estimate (1.47)
we have that (12.6) is
\begin{eqnarray*}
\!\!\!& \!\!\ll \!\! & h\log^{2}R \!\!
\sum_{\stackrel{\scr d,\, e \leq R}{\scr (d,e)<h}}
{\mu^{2}(d)\mu^{2}(e)de\over \phi(d)\phi(e)} 
{h\over \sqrt{(d,e)}}[(1+{h^{2}\over de})\!\!\sum_{j(\bmod \, [d,e])}
\max_{u\leq 2N+h}|E(u;[d,e],j)|^{2}]^{{1\over 2}} \\
\!\!\!& \!\! \ll \!\! & N^{{1\over 2}}h^2 \log^{4}N \!\!
\sum_{\stackrel{\scr d,\, e \leq R}{\scr (d,e)<h}}
{\mu^{2}(d)\mu^{2}(e)de\over \phi(d)\phi(e)\sqrt{(d,e)}} 
+ N^{{1\over 2}}h^3 \log^{4}N \!\!
\sum_{\stackrel{\scr d,\, e \leq R}{\scr (d,e)<h}}
{\mu^{2}(d)\mu^{2}(e)\sqrt{de}\over \phi(d)\phi(e)\sqrt{(d,e)}}.
\end{eqnarray*}
We have 
\begin{eqnarray*}
\sum_{\stackrel{\scr d,\, e \leq R}{\scr (d,e)<h}}
{\mu^{2}(d)\mu^{2}(e)de\over \phi(d)\phi(e)\sqrt{(d,e)}} 
& \ll & \sum_{\delta < h}{\delta^{{3\over 2}}\over
\phi^{2}(\delta)}(\sum_{d'\leq {R\over \delta}}{d' \over \phi(d')})^2
\ll R^2 \sum_{\delta<h}{1\over \delta^{{1\over 2}}\phi^{2}(\delta)} 
\ll R^2 , \\
\sum_{\stackrel{\scr d,\, e \leq R}{\scr (d,e)<h}}
{\mu^{2}(d)\mu^{2}(e)\sqrt{de}\over \phi(d)\phi(e)\sqrt{(d,e)}} 
& \ll & \sum_{\delta < h}{\delta^{{1\over 2}}\over
\phi^{2}(\delta)}(\sum_{d'\leq {R\over \delta}}{\sqrt{d'}\over \phi(d')})^2
\ll R \sum_{\delta<h}{1\over \delta^{{1\over 2}}\phi^{2}(\delta)} \ll R.
\end{eqnarray*}
Hence the contribution of the error terms in (12.3) to (12.1) is
\begin{equation}
\ll N^{{1\over 2}}h^2 R^2 \log^{4}N 
+ N^{{1\over 2}}h^3 R \log^{4}N .
\end{equation}

Now we calculate the contribution of (12.4) to (12.1), that of the error
term being $O(Nh^3 R^{-1+\epsilon})$. Inverting (1.24) we have
\begin{equation}
\gs((j_{1},j_{2},j_{3})) =
\sum_{\mathcal{J}\subset\{j_1,j_2,j_3\}}\frU(\mathcal{J}),
\end{equation}
which implies by (1.23) that (since $\frU(\emptyset) =1,\, \frU((k))=0,\,
\frU((k,l))=\gs((k,l))-1$ for $k\neq l$)
\begin{eqnarray}
\sum_{\stackrel{\scr 1\le j_1 ,j_2, j_3 \le h}{ \scr \mathrm{ distinct}}}
\!\!\!\!\!\!\gs((j_{1},j_{2},j_{3})) \!\! & = &\!\! -2h(h-1)(h-2) + 3(h-2)
\!\!\!\sum_{\stackrel{\scr 1\le j_1 ,j_2 \le h}{ \scr \mathrm{ distinct}}}
\!\!\!\!\!\!\gs((j_{1},j_{2})) + R_{3}(h) \nonumber \\
& = & \!\!-2h(h-1)(h-2) + 6(h-2)\!\!\!\sum_{1\leq j \leq h}
\!\!\!(h-j)\gs_{2}(j) +\!R_{3}(h). 
\end{eqnarray}
From (10.7) and (1.26) we obtain
\begin{equation}
\sum_{\stackrel{\scr 1\le j_1 ,j_2, j_3 \le h}{ \scr \mathrm{ distinct}}}
\gs((j_{1},j_{2},j_{3})) = h^3 -3h^2 \log h + 3(1-\gamma-\log 2\pi)h^2 +
O(h^{{3\over 2}+\epsilon})
\end{equation} 
(note that (1.26) has not been used in full force, here we only need to know
$R_{3}(h) \ll h^{{3\over 2}+\epsilon}$).
Hence we have 
\begin{eqnarray}
\lefteqn{
\sum_{\stackrel{\scr 1\le j_1 ,j_2, j_3 \le h}{ \scr \mathrm{ distinct}}}
\tilde{\mathcal{ S}}_{3}'(N,(j_{1},j_{2},j_{3}),(1,1,1)) } \\
& & = Nh^3 -3Nh^2 \log h + 3(1-\gamma-\log 2\pi)Nh^2 
+ O(Nh^{{3\over 2}+\epsilon}) \nonumber \\
& & \qquad 
+ O(Nh^3 R^{-1+\epsilon}) + O((N^{{1\over 2}}h^2  R\log^{4}N)(R+h)) .\nonumber 
\end{eqnarray}
Combining (1.37), (10.8), (11.14) and (12.11) in (11.6) we obtain
\begin{eqnarray}
\tilde{M}_{3}'(N,h,\psi_{R}) & = & Nh^3 + 3Nh^2 (\mathcal{L}_{1}(R) -\log h)
+ 3Nh^2 (1-\gamma - \log2\pi) \nonumber \\
& & \qquad + O(Nh^{{3\over 2}+\epsilon}) 
+ O(Nh^3 R^{-1+\epsilon}) +O(Nh\log^{8}N) \nonumber \\
& & \qquad + O(N^{{1\over 2}}h^2 R^2 \log^{4}N) 
+ O(N^{{1\over 2}}h^3 R\log^{4}N) . 
\end{eqnarray}
Now we put our findings together in (11.1). In doing this, keeping the sums 
$\displaystyle \sum_{1\leq j \leq h}(h-j)\gs_{2}(j)$ (which appear in
several terms) unevaluated until the end not only 
facilitates the calculation, but also reveals the complete cancellation of 
the terms which contain
$\displaystyle Nh\sum_{1\leq j \leq h}(h-j)\gs_{2}(j)$. We get
\begin{eqnarray}
\calm_{3}' & = & (2C+\rho)ANh[\log h - \mathcal{L}_{1}(R) 
-(1-\gamma-\log 2\pi)] -\rho C^2 A^3 N \nonumber \\
&  & \qquad + O(ANh^{{1\over 2}+\epsilon}) 
+ O(N^{{1\over 2}}h^2 R^2 \log^{4}N) + O(N^{{1\over 2}}h^3 R \log^{4}N)
\nonumber \\
& & \qquad + O(Nh\log^{8}N) + O(h^3 R^2) + O(Nh^3 R^{-1+\epsilon}) + NR_{3}(h).
\end{eqnarray}
The main terms are at the same order of magnitude if $A=(h\log N)^{{1\over
2}}$ as before in (10.11) for $\calm_{2}'$. This choice of $A$ makes (12.13)
read as
\begin{eqnarray}
\calm_{3}' & = & Nh^{{3\over 2}}\log^{{1\over 2}}N[-\rho C^2 \log N
+(2C+\rho)(\log h - \mathcal{L}_{1}(R) -(1-\gamma-\log 2\pi)] \nonumber \\
& & \qquad + \; \mathrm{error \; terms \; of \; (12.13)} .
\end{eqnarray}
We are assuming that $R\gg N^{\epsilon}$, and the requirement that the error 
terms are smaller than the main term, i.e. 
$o(Nh^{{3\over 2}}\log^{{3\over 2}}N)$, brings the restrictions
\begin{equation}
h \ll R^{{2\over 3} -\epsilon}, \qquad \log^{13}N =o(h),\qquad
h^{{1\over 2}}R^2 = o(N^{{1\over 2}}\log^{-{5\over 2}}N) .
\end{equation}
Note that
the cancellation mentioned before (12.13) is essential in reaching a result,
for if (12.12) and (10.9) which depend on the
evaluation (10.7) had been used, then we would have acquired an error term
$O(Nh^{{3\over 2}+\epsilon})$ that is larger than the main term.
Upon this we need Montgomery and Soundararajan's estimate (1.26) for
$R_{3}(h)$.
Thus for
\begin{equation}
\log^{14}N \ll h \ll N^{{1\over 7} -\epsilon}
\end{equation}
we have the asymptotic result
\begin{equation}
\calm_{3}' \sim -Nh^{{3\over 2}}\log^{{1\over 2}}N(\rho C^2 \log N +
(2C+\rho )\log {R\over h}) .
\end{equation}
(The factor $N^{\epsilon}$ in (12.16) can be replaced by a small power
of $\log N$ if one bounds (5.38) more precisely as was remarked).
From (12.17) we can get the result (10.16), only this time for the smaller
range (12.16). The significance of (12.17) is that it allows us also
to get a result of the type (10.16) without the absolute value.
It is convenient to write $R=N^{\theta}, \, h=N^{\alpha}$. For a fixed
$\rho$ satisfying $0 < |\rho| < \sqrt{\theta - \alpha}$, with the choice
$C = -{\theta - \alpha\over \rho}$, (12.17) reads
\begin{equation}
\calm_{3}' \sim -\rho (\theta -\alpha)(1-{\theta - \alpha\over \rho^2 })
Nh^{{3\over 2}}\log^{{3\over 2}}N .
\end{equation}
We see that $\calm_{3}'$ is positive for $0 < \rho < \sqrt{\theta - \alpha}$, 
and negative for $-\sqrt{\theta -\alpha} < \rho < 0$,
and $\gg N(h\log N)^{{3\over 2}}$ in either case. This means that given an
arbitrarily small but fixed $\eta > 0$, for all
sufficiently large $N$ and $h$ subject to (12.16), there
exist $n_{1},\, n_{2} \in [N+1,2N]$ such that
\begin{eqnarray}
\psi(n_1 +h) -\psi(n_1) - h & > & ({\sqrt{1-5\alpha}\over 2}-\eta) 
(h\log N)^{{1\over 2}} \nonumber \\
\psi(n_2 +h) -\psi(n_2) - h & < & -({\sqrt{1-5\alpha}\over 2}-\eta) 
(h\log N)^{{1\over 2}} .
\end{eqnarray}


\begin{thebibliography}{10}

\bibitem{D} H. Davenport, \textit{Multiplicative Number Theory}, Third edition,
Revised and with a preface by H. L. Montgomery, Springer-Verlag, New York,
2000.

\bibitem {GA} P. X. Gallagher, \textit{On the distribution of primes in short 
intervals}, Mathematika \textbf{23} (1976), 4--9.

\bibitem{FG} D. A. Goldston, \textit{Linnik's theorem on Goldbach numbers
in short intervals}, Glasgow Math. J. \textbf{32} (1990), 285--297.

\bibitem{G1} D. A. Goldston, \textit{ A lower bound for the second moment of 
primes in short intervals},  Expo. Math. \textbf{13} (1995), 366--376.

\bibitem{G2} D. A. Goldston, \textit{A lower bound method for binary
additive problems involving primes}, Unpublished manuscript

\bibitem{GM} D. A. Goldston and H. L. Montgomery, \textit{Pair correlation
of zeros  and primes in short intervals,} Analytic Number Theory and
Diophantine Problems, Birkha\"user, Boston, Mass., 1987, 183--203.  

\bibitem{GY1}  D. A. Goldston and C. Y. Y{\i}ld{\i}r{\i}m, \textit{On the 
second moment for primes in an arithmetic progression}, Acta 
Arithmetica, \textbf{C.1} (2001), 85--104.

\bibitem{GY2}  D. A. Goldston and C. Y. Y{\i}ld{\i}r{\i}m, \textit{Higher
correlations of divisor sums related to primes I: Triple correlations},
Integers {\bf 3} (2003), A5, 66 pp. (electronic).

\bibitem{HL} G. H. Hardy and J. E. Littlewood, \textit{ Some problems of 
`Partitio Numerorum': III On the expression of a number as a sum of primes}, 
Acta Math. \textbf{ 44} (1923), 1--70.

\bibitem{Hi} A. Hildebrand, \textit{ \"{U}ber die punktweise Konvergenz von 
Ramanujan-Entwicklungen zahlentheoretischer Funktionen}, Acta Arithmetica 
\textbf{XLIV} (1984),  109--140.

\bibitem{HR} H. Halberstam and H. -E. Richert, \textit{On a result of R. R. 
Hall}, J. Number Theory (1) \textbf{11} (1979), 76--89.

\bibitem{HB} D. R. Heath-Brown,  \textit{The ternary Goldbach problem}, Rev. 
Mat. Iberoamericana \textbf{1} (1985), no. 1, 45--59.

\bibitem{H6} C. Hooley, \textit{On the Barban-Davenport-Halberstam theorem 
VI}, J. London Math. Soc. (2) \textbf{13} (1976), 57--64. 

\bibitem{H12} C. Hooley, \textit{On the Barban-Davenport-Halberstam theorem 
XII}, Number Theory in Progress (Zakopane 1997), Vol. II, Walter de Gruyter,
1999, 893--910.

\bibitem{H13} C. Hooley, \textit{On the Barban-Davenport-Halberstam theorem 
XIII}, Acta Arith. \textbf{94} (2000), no. 1, 53--86. 

\bibitem{KV} A. A. Karatsuba and S. M. Voronin, \textit{The Riemann
Zeta-function} (translated from the Russian by N. Koblitz), 
Walter de Gruyter, Berlin, 1992.

\bibitem {MS} H. L. Montgomery and K. Soundararajan, \textit{Primes in
short intervals}, Comm. Math. Phys. (2005), 29 pages (to appear).

\bibitem{S1} A. Selberg, \textit{ On an elementary method in the theory of
primes}, Norske Vid. Selsk. Forh. Trondhjem \textbf{19} (1947), 64--67.

\bibitem{T} E. C. Titchmarsh, \textit{The theory of the Riemann 
zeta-function}, Second edition. Edited and with a preface by D. R. 
Heath-Brown. The Clarendon Press, Oxford University Press, New York, 1986.



\end{thebibliography}
\end{document}